\documentclass[12pt,a4paper,reqno]{article}

\pdfoutput=1

\usepackage{amsmath}
\usepackage{amsfonts}
\usepackage{amssymb}

\usepackage[utf8]{inputenc}
\usepackage[T1]{fontenc}
\usepackage{lmodern}

\usepackage{enumerate}
\usepackage{geometry}
\usepackage{graphicx}
\usepackage[unicode,bookmarks]{hyperref}
\usepackage{mathbbol}
\usepackage[numbers,square]{natbib}
\usepackage{subfig}

\newcommand{\bracket}[1]{\ensuremath{\left < #1 \right >}}

\newcommand{\basis}[1]{\ensuremath{\mathrm{#1}}}
\newcommand{\coord}[1]{\ensuremath{#1}}
\newcommand{\field}[1]{\ensuremath{#1}}

\newcommand{\bx}{\ensuremath{\basis{x}}}
\newcommand{\by}{\ensuremath{\basis{y}}}
\newcommand{\bz}{\ensuremath{\basis{z}}}
\newcommand{\bw}{\ensuremath{\basis{w}}}

\newcommand{\cx}{\ensuremath{\coord{x}}}
\newcommand{\cy}{\ensuremath{\coord{y}}}
\newcommand{\cz}{\ensuremath{\coord{z}}}
\newcommand{\ct}{\ensuremath{\coord{t}}}

\newcommand{\fu}{\ensuremath{\field{u}}}
\newcommand{\fv}{\ensuremath{\field{v}}}
\newcommand{\fw}{\ensuremath{\field{w}}}

\newcommand{\jb}{\ensuremath{\mathrm{J}}}
\newcommand{\jp}{\ensuremath{\mathrm{j}}}
\newcommand{\jpd}{\ensuremath{\mathrm{j}}}

\newcommand{\eps}{\ensuremath{\epsilon}}
\newcommand{\phy}{\ensuremath{\varphi}}

\newcommand{\del}{\ensuremath{\mathrm{D_{EL}}}}
\newcommand{\id}{\ensuremath{\mathrm{id}}}
\newcommand{\fr}{\ensuremath{D}}

\let\div\undefined

\DeclareMathOperator{\div}{div}

\providecommand{\keywords}[1]{%
\noindent\rule{3cm}{0.4pt}\linebreak%
\small{\textit{Keywords:} #1}%
}

\begin{document}

\title{\vspace{-1em}Variational Integrators for Nonvariational\\ Partial Differential Equations}

\author{
\large{Michael Kraus}\\
\small{(michael.kraus@ipp.mpg.de)}
\vspace{.5em}\\
\normalsize{Max-Planck-Institut f\"ur Plasmaphysik}\\
\normalsize{Boltzmannstra\ss{}e 2, 85748 Garching, Deutschland}%
\vspace{.5em}\\
\normalsize{Technische Universit\"at M\"unchen, Zentrum Mathematik}\\
\normalsize{Boltzmannstra\ss{}e 3, 85748 Garching, Deutschland}%
\vspace{1em}\\
\and
\large{Omar Maj}\\
\small{(omar.maj@ipp.mpg.de)}
\vspace{.5em}\\
\normalsize{Max-Planck-Institut f\"ur Plasmaphysik}\\
\normalsize{Boltzmannstra\ss{}e 2, 85748 Garching, Deutschland}%
\vspace{1em}\\
}

\date{\today}

\maketitle

\begin{abstract}
Variational integrators for Lagrangian dynamical systems provide a systematic way to derive geometric numerical methods.
These methods preserve a discrete multisymplectic form as well as momenta associated to symmetries of the Lagrangian via Noether's theorem.
An inevitable prerequisite for the derivation of variational integrators is the existence of
a variational formulation for the considered problem.
Even though for a large class of systems this requirement is fulfilled, there are many interesting examples which do not belong to this class, e.g., equations of advection-diffusion type frequently encountered in fluid dynamics or plasma physics.

On the other hand, it is always possible to embed an arbitrary dynamical system into a larger Lagrangian system using the method of formal (or adjoint) Lagrangians.
We investigate the application of the variational integrator method to formal Lagrangians, and thereby extend the application domain of variational integrators to include potentially all dynamical systems.

The theory is supported by physically relevant examples, such as the advection equation and the vorticity equation, and numerically verified.
Remarkably, the integrator for the vorticity equation combines Arakawa's discretisation of the Poisson brackets with a symplectic time stepping scheme in a fully covariant way such that the discrete energy is exactly preserved.
In the presentation of the results, we try to make the geometric framework of variational integrators accessible to non specialists.
\end{abstract}


\newpage

\tableofcontents

\vspace{3em}

\section{Introduction}
\label{sec:intro}

In recent years, the field of structure-preserving or geometric discretisation 
\cite{Christiansen:2011, HairerLubichWanner:2006, BuddPiggot:2000} has become a flourishing discipline of numerical analysis and scientific computing. One particular family of geometric discretisation methods is that of variational integrators \cite{WendlandtMarsden:1997, MarsdenWendlandt:1997, MarsdenPatrick:1998, KouranbaevaShkoller:2000, MarsdenWest:2001, Lew:2004b, Lew:2004a}, which are based on the discretisation of Hamilton's principle of stationary action \cite{JoseSaletan:1998, Arnold:1989, Holm:2009, MarsdenRatiu:2002, AbrahamMarsden:1978}.
Variational integrators preserve a discrete multisymplectic form and have good longtime energy behaviour. As we will see, they can be designed to preserve energy even exactly, which in practice means up to machine precision.
Furthermore, they preserve momenta associated to symmetries of the discrete equations of motion via a discrete version of Noether's theorem \cite{Noether:1918, KosmannSchwarzbach:2010}.

While in most standard discretisation techniques for dynamical systems the equations of motion are directly discretised,
the basic idea of variational integrators is to construct a discrete counterpart to the considered system.
This means that the fundamental building blocks of classical mechanics and field theory, namely the action functional, the Lagrangian,
the variational principle, and the Noether theorem, all have discrete equivalents.
The application of the discrete variational principle to the discrete action then leads to discrete Euler-Lagrange equations.
The evolution map that corresponds to the discrete Euler-Lagrange equations is what is called a variational integrator.
The discrete Noether theorem can be used to relate symmetries of the discretised system to discrete momenta
that are exactly preserved by this integrator.
Whereas most standard techniques put emphasis on the minimisation of local errors, for variational
integrators the focus is rather on the preservation of global or geometric properties of the system.

An obvious limitation of the variational integrator method is its applicability to Lagrangian systems only.
This excludes a large class of interesting systems, for example the problems of advection-diffusion type often found in fluid dynamics and plasma physics.
We propose here that the method of formal (or adjoint) Lagrangians \cite{AthertonHomsy:1975} can be used as an expedient to avoid this limitation.
More specifically, formal Lagrangians allow us to embed any given system into a larger system which, in turn, admits a Lagrangian formulation.
To obtain a formal Lagrangian $L$, the equation at hand, say $\mathcal{F} [\fu] = 0$, is multiplied by an adjoint variable $\fv$, giving $L = \fv \cdot \mathcal{F} [\fu]$.
Variation of the resulting action functional, $\mathcal{A} = \int L \, d^{n+1} \bx$, with respect to $\fv$ gives the original equation $\mathcal{F} [\fu] = 0$.
Variation of the action functional with respect to the physical variable $\fu$ gives an additional equation that determines the evolution of the adjoint variable $\fv$.

At first sight one might be tempted to regard the formal Lagrangian formalism as merely a method for obtaining a weak formulation of the problem at hand. 
Then, if our goal is to obtain an integrator, the details of the dynamics of the adjoint variable $\fv$ would seem irrelevant.
However, it turns out that the dynamics of $\fv$ play a key role in relating symmetries of the formal Lagrangian to conservation laws satisfied by $\fu$.
\citeauthor{Ibragimov:2006} \cite{Ibragimov:2006, Ibragimov:2007a, Ibragimov:2007b} developed a theory for the analysis of conservation laws of arbitrary differential equations by extending the Noether theorem to formal Lagrangians. This leads to conservation laws for the extended system $(\fu, \fv)$, which can be restricted to the original system provided that it is possible to express the solution of the adjoint variable $\fv$ in terms of $\fu$.

In this work, we propose the combination of the discrete variational principle with Ibragimov's theory in order to derive variational integrators for systems without a natural Lagrangian formulation and to determine the associated discrete conservation laws.
Thereby we extend significantly the applicability of the variational integrator method.
The goal of this approach is to design numerical schemes which respect certain conservation laws of a given system in a rather systematic way.

We proceed as follows. In section \ref{sec:vi}, we present the theory of variational integrators in simple terminology. To set the stage and fix notation we review the continuous action principle for field theories and the corresponding Noether theorem before passing over to the discrete theory, which is extended to account for discrete divergence symmetries. The style of presentation is chosen to make the theory accessible to a wide audience without extensive background in modern differential geometry. This implies some loss of generality, but hopefully not too much of the geometric beauty of the original work is lost.
In section \ref{sec:adjoint}, we recall the inverse problem of the calculus of variations, review the theory of formal Lagrangians and explain the derivation of conserved quantities in this setting. We also provide a geometric formulation of the theory, which to our knowledge has not been presented, yet.
Finally, in section \ref{sec:applications}, we apply the method to some prototypical examples, including the advection equation and the vorticity equation, and verify the theoretical properties in numerical experiments.
More elaborate numerical examples for the Vlasov-Poisson system as well as ideal and reduced magnetohydrodynamics will be presented elsewhere \cite{KrausMajSonnendruecker:2015, KrausMaj:2015, KrausTassi:2015}.
The examples we provide here can be seen as building blocks for these more complicated systems.

\section{Variational Integrators}
\label{sec:vi}

\subsection{Geometry and Notation}

In this work, we are concerned with the discretisation of partial differential equations (PDEs) of evolution type.
A field is a map $\fu : X \rightarrow F$ from a bounded domain $X \subset \mathbb{R}^{n+1}$ taking values in an open set $F \subseteq \mathbb{R}^{m}$.
Most often, $X$ corresponds to some region of spacetime with coordinates
\begin{align*}
\bx
&= (\bx^{\mu})
= (\bx^{0}, \bx^{i})
= (\ct, \cx, \cy, \cz) &
& \text{with} &
0 \leq \mu \leq n , \;
1 \leq i \leq n ,
\end{align*}
and $n = \dim{X} - 1$ being the number of space-like dimensions.
Most of the theoretic results rely neither on this space-plus-time splitting nor on $X$ being a subset of an Euclidean space. Thus $X$ can be replaced by a differentiable manifold, in which case the results are valid locally in a coordinate chart.
Points on $F$ are denoted by $\by = (\by^{a})$ with $1 \leq a \leq m$ and $m = \dim{F}$ being the number of field components.
We make use of the Einstein summation convention on repeated indices, both for coordinates and fields.

The configuration of a field $\fu$ is geometrically represented by its graph
\begin{align*}
\mathrm{graph} (\fu) = \big\{ (\bx, \by) \; \big\vert \; \by = \fu(\bx) \big\} ,
\end{align*}
which is a subset of the Cartesian product
\begin{align*}
Y = X \times F = \big\{ (\bx, \by) \; \big\vert \; \bx \in X, \; \by \in F \big\} .
\end{align*}
With the projection onto the first factor,
\begin{align*}
\arraycolsep=2pt
\begin{array}{rcccc}
\pi & : & Y & \rightarrow & X , \\
& & (\bx, \by) & \mapsto & \bx ,
\end{array}
\end{align*}
we can define a geometrical structure $(Y, X, \pi)$ called (trivial) fibre bundle.
Here, $X$ is called the base space, $Y$ the configuration space, and $F$ the fibre.
Local coordinates on $Y$ are given by
\begin{align*}
& ( \bx^{\mu}, \by^{a} ) 
& \text{with} &
& 0 \leq \mu \leq n , \; 1 \leq a \leq m .
\end{align*}
A field $\fu$ can be identified with a section of the bundle, i.e., a map $\phy : X \rightarrow Y$, satisfying the condition
\begin{align*}
\pi \circ \phy = \id_{X} ,
\end{align*}
where $\id_{X}$ is the identity map on $X$.
This condition ensures that the image of a section $\phy$ corresponds to the graph of a field $\fu$,
\begin{align*}
\phy (X)
= \big\{ (\bx, \by) = \phy (\bx) \; \big\vert \; \bx \in X \big\}
= \big\{ (\bx, \fu (\bx)) \; \big\vert \; \bx \in X \big\} .
\end{align*}
In local coordinates, a section $\phy : X \rightarrow Y$ can be written as $\phy (\bx) = (\bx^{\mu}, \phy^{a} (\bx) )$, that is $\phy (\bx)$ is given by functions $\by^{a} = \phy^{a} (\bx) = \fu^{a} (\bx)$.

Partial derivatives of a field component $\fu^{a}$ with respect to $\bx^{\mu}$ or $\bx^{i}$ are denoted $\fu^{a}_{\mu}$ or $\fu^{a}_{i}$, respectively.
The collection of all partial derivatives of a given order $k$ is denoted $\fu_{(k)}$.
The appropriate geometric setting for partial derivatives of a field $\fu$ is the theory of jet bundles \cite{Saunders:1989, Kolar:1993, Olver:1995, Ivancevic:2007, Campos:2010}.
Jets combine into a single object the values of a field at a point and the values of its derivatives.
More specifically, the jet prolongation of a section $\phy$ is a map
\begin{align*}
\jp^{1} \phy : \bx \mapsto \big( \bx^{\mu} , \phy^{a} (\bx) , \phy^{a}_{\mu} (\bx) \big) ,
\end{align*}
taking values in $Y \times \mathbb{R}^{m (n+1)}$. This space is identified with the first jet bundle of $Y$, denoted $\jb^{1} Y$, so that the first jet prolongation of a field is a section of $\jb^{1} Y$.
Local coordinates on $\jb^{1} Y$ are given by
\begin{align*}
& ( \bx^{\mu}, \by^{a}, \bz^{a}_{\mu} )
& \text{with} &
& 0 \leq \mu \leq n , \; 1 \leq a \leq m ,
\end{align*}
where $\bz^{a}_{\mu}$ represents the possible values of partial derivatives.
The Lagrangian of a first order field theory for example will be a function defined on the first jet bundle $\jb^{1} Y$.
The jet bundle $\jb^{1} Y$ has a natural projection on both $Y$ and $X$, denoted by $\pi_{Y}$ and $\pi_{X}$, respectively. The latter defines a fibre bundle $(\jb^{1} Y, X, \pi_{X})$. 
It is worth noting that not every section $\psi$ of the bundle $\jb^{1} Y \rightarrow X$ is the jet prolongation of a section $\phy : X \rightarrow Y$. When that happens we say that $\psi$ is holonomic.

Higher order jet bundles $\jb^{k} Y$ of $Y$ can be defined analogously by considering the $k$th jet prolongation of a section $\phy$,
\begin{align*}
\jp^{k} \phy : \bx \mapsto \big( \bx^{\mu} , \phy^{a} (\bx) , \phy^{a}_{\mu} (\bx) , \phy^{a}_{\mu \nu} (\bx) , ... \big) ,
\end{align*}
which includes derivatives of $\phy$ up to $k$th order.
Higher-order jet bundles $\jb^{k} Y$ can also be defined as submanifolds of iterated first-order jet bundles,
 e.g., $\jb^{2} Y$ is a submanifold of $\jb^{1} ( \jb^{1} Y )$.
If coordinates on $\jb^{1} ( \jb^{1} Y)$ are denoted by
\begin{align*}
& ( \bx^{\mu}, \by^{a}, \bz^{a}_{\mu}, \bar{\bz}^{a}_{\mu}, \bw^{a}_{\mu \nu} )
& \text{with} &
& 0 \leq \mu, \nu \leq n , \; 1 \leq a \leq m ,
\end{align*}
then $\jb^{2} Y$ is obtained by requiring that $\bz^{a}_{\mu} = \bar{\bz}^{a}_{\mu}$ for all $a$ and all $\mu$,
and accordingly for larger $k$ (for more details see e.g. \cite[Chapter 32]{Kolar:1993} or \cite[Section 4.1.4]{Campos:2010}).

Jet bundles have become the standard framework for PDE analysis and variational calculus.
They provide a natural setting for the formulation of field theories and the analysis of conservation laws which are central to this work.
This framework generalises in a standard way to non-trivial cases in which $X$
is a manifold and $Y$ cannot be written globally as a Cartesian product.
For a thorough introduction into the geometric framework we refer to \citet{GotayMarsden:1998}.

\subsection{Continuous Action Principle}
\label{sec:vi_continuous_action_principle}

For definiteness, let us consider only first order Lagrangian field theories,
i.e., the Lagrangian $L$ shall depend only on the coordinates, the fields and their first derivatives.
Therefore, the Lagrangian is a function defined on the first jet bundle,
\begin{align}\label{eq:variational_continuous_lagrangian}
L : \jb^{1} Y \rightarrow \mathbb{R} .
\end{align}
The corresponding action functional is
\begin{align}\label{eq:variational_continuous_action}
\mathcal{A} [\phy] = \int \limits_{X} L \big( \jp^{1} \phy \big) \, d^{n+1} \bx .
\end{align}
Hamilton's principle of stationary action \cite{JoseSaletan:1998, Arnold:1989, AbrahamMarsden:1978} states that among all possible field configurations $\phy : X \rightarrow Y$, the one chosen by nature makes the action functional (\ref{eq:variational_continuous_action}) stationary.
As usual in geometric mechanics \cite{Holm:2009, MarsdenRatiu:2002, GotayMarsden:1998, MarsdenPatrick:1998}, stationary points of the action are meant in the formal sense.
The variations of a field $\phy$ are defined in terms of a geometric transformation of the underlying bundle $Y$, that is a sufficiently regular one-parameter group of transformations $\sigma^\eps (\bx, \by)$ of $Y$ into itself, defined for $\eps$ in a neighbourhood of zero.
In order to have near identity transformations, it is required that $\sigma^\eps$ reduces to the identity map at $\eps = 0$, i.e., $\sigma^\eps \vert_{\eps = 0} = \id$.
Furthermore, variations should vanish at boundary points, so that $\sigma^\eps = \id$ for $\bx \in \partial X$.
The map $\sigma^\eps (\bx, \by)$ is the flow of the vector field $V$ over $Y$ given by
\begin{align}\label{eq:variational_continuous_vector_field}
V (\bx, \by) = \dfrac{d}{d\eps} \sigma^\eps (\bx, \by) \bigg\vert_{\eps=0} ,
\end{align}
and written in components as
\begin{align}\label{eq:variational_continuous_vector_field_components}
V (\bx, \by) = \eta^{a} (\bx, \by) \, \dfrac{\partial}{\partial \by^{a}} ,
\end{align}
where we have $\eta^{a} (\bx, \by) = 0$ for $\bx \in \partial X$.
Here, vector fields are identified with first-order differential operators.
The variation of a field $\phy$ in the direction of $V$ is then defined by
\begin{align}
\phy^{\eps} = (\sigma^\eps \circ \phy) (\bx) = \sigma^\eps ( \phy(\bx) ) ,
\end{align}
and we have $\phy^{\eps} \vert_{\eps=0} = \phy$ and $\phy^{\eps} = \phy$ on boundary points $\bx \in \partial X$.
Loosely speaking, $\phy^{\eps}$ is obtained by dragging $\phy$ along the flow of the vector field $V$.
For sake of simplicity, we consider transformations in $\sigma^\eps$ that leave the point $\bx$ unchanged, thus, preserving the fibres of the bundle $Y$.
Those are referred to as vertical transformations\footnote{
As noted in \cite{MarsdenPatrick:1998, MarsdenPekarsky:2001}, one really should consider general transformations that might affect the points of the base space as well as the fields. Otherwise, the Cartan form and thus the multisymplectic form will not be correctly obtained from the variational principle.
However, in this work we are only concerned with the Euler-Lagrange equations and the Noether theorem. Both will be obtained correctly even if only transformations of the fields (vertical transformations in the language of jet bundles) are considered.
Transformations of the base space (horizontal transformations) would substantially complicate the derivations and are neglected for sake of simplicity.
}.
We say that $\phy$ is a stationary point of $\mathcal{A} [\phy]$ if for every $\sigma^\eps$ we find that $\eps=0$ is a stationary point of $\mathcal{A} [\phy^{\eps}]$ viewed as a function of $\eps$ in an open interval about $\eps = 0$. Then Hamilton's principle amounts to
\begin{align}\label{eq:variational_continuous_variation_1}
\dfrac{d}{d\eps} \mathcal{A} [ \sigma^\eps \circ \phy ] \bigg\vert_{\eps=0} = 0
\quad
\text{for every $\sigma^\eps$} .
\end{align}
Explicitly,
\begin{align}
\dfrac{d}{d\eps} \mathcal{A} [ \sigma^\eps \circ \phy ] \bigg\vert_{\eps=0}
\label{eq:variational_continuous_variation_2}
&= \int \limits_{X} \dfrac{d}{d\eps} \Big[ L \big( \jp^{1} ( \sigma^\eps \circ \phy) \big) \Big] \bigg\vert_{\eps=0} \, d^{n+1} \bx \\
\label{eq:variational_continuous_variation_3}
&= \int \limits_{X} \bigg[ \dfrac{\partial L}{\partial \by^{a}} \big( \jp^{1} \phy \big) \cdot \dfrac{d (\sigma^\eps)^{a}}{d \eps} + \dfrac{\partial L}{\partial \bz^{a}_{\mu}} \big( \jp^{1} \phy \big) \cdot \dfrac{d}{d \eps} \dfrac{\partial (\sigma^\eps)^{a}}{\partial \bx^{\mu}} \bigg] \bigg\vert_{\eps=0} \, d^{n+1} \bx .
\end{align}
The second term in (\ref{eq:variational_continuous_variation_3}) can be integrated by parts,
\begin{align}\label{eq:variational_continuous_variation_4}
\dfrac{d}{d\eps} \mathcal{A} [ \sigma^\eps \circ \phy ] \bigg\vert_{\eps=0}
&= \int \limits_{X} \bigg[ \dfrac{\partial L}{\partial \by^{a}} \big( \jp^{1} \phy \big) - \dfrac{\partial}{\partial \bx^{\mu}} \left( \dfrac{\partial L}{\partial \bz^{a}_{\mu}} \big( \jp^{1} \phy \big) \right) \bigg] \cdot \eta^{a} (\phy) \, d^{n+1} \bx ,
\end{align}
where we used (\ref{eq:variational_continuous_vector_field}) and (\ref{eq:variational_continuous_vector_field_components}).
As the variation of the action has to vanish for arbitrary transformations $\sigma^\eps$ and therefore for arbitrary vector fields $V$, the term in square brackets has to vanish identically, that is
\begin{align}\label{eq:variational_continuous_euler_lagrange}
\big( \jp^{2} \phy \big)^{*} \big( \del(L) \big)^{a} =
\dfrac{\partial L}{\partial \by^{a}} \big( \jp^{1} \phy \big) - \dfrac{\partial}{\partial \bx^{\mu}} \left( \dfrac{\partial L}{\partial \bz^{a}_{\mu}} \big( \jp^{1} \phy \big) \right) = 0 .
\end{align}
These are the Euler-Lagrange field equations, i.e., the equations of motion for a first order Lagrangian field theory, which for a non-degenerate Lagrangian $L$ live in the second-order jet bundle $\jb^{2} Y$.
By $\del$ we denote the Euler-Lagrange operator acting on the Lagrangian $L$. The pull-back notation $(\jp^{2} \phy)^{*}$ states that the result is evaluated on the second jet prolongation $\jp^{2} \phy$.
The theory is fully covariant, that is for time-dependent problems, time is regarded as a component of $\bx$.

\subsection{Continuous Noether Theorem}
\label{sec:vi_continous_noether_theorem}

The Noether theorem \cite{Noether:1918, KosmannSchwarzbach:2010}
states that each Lie point symmetry of a Lagrangian corresponds to a conservation law of the associated Euler-Lagrange equations.

We restrict our attention to conservation laws that are generated by vertical transformations of the configuration bundle $Y$, i.e., transformations which leave the base space $X$ invariant.
In the framework of formal Lagrangians addressed below, this is often sufficient to uncover interesting conservation laws (including conservation of momentum and energy).
Our derivation of the Noether theorem is essentially based on reference \citep{MarsdenPatrick:1998}.

As in section \ref{sec:vi_continuous_action_principle}, a transformation is generated by a map $\sigma^\eps (\bx, \by)$, that is
\begin{align}\label{eq:noether_continuous_transformation}
\phy^{\eps} &= \sigma^\eps \circ \phy &
& \text{with} &
\sigma^\eps \vert_{\eps=0} &= \mathrm{id} &
& \text{and} &
V = \dfrac{d\sigma^\eps}{d\eps} \bigg\vert_{\eps=0} ,
\end{align}
but it is not required that $\sigma^\eps$ reduces to the identity at boundary points.
In our analysis we will usually just prescribe the generating vector field $V$ instead of the actual transformation $\sigma^\eps$.
In components, $V$ can be written as
\begin{align}\label{eq:noether_continuous_vector_field}
V (\bx, \by) &= \eta^{a} (\bx, \by) \, \dfrac{\partial}{\partial \by^{a}} &
& \text{with} &
\eta^{a} (\bx, \by) &= \dfrac{d}{d\eps} (\sigma^\eps)^{a} (\bx, \by) \bigg\vert_{\eps=0} ,
\end{align}
where $(\sigma^\eps)^{a}$ is the $a$-th component of the transformation map. 
Since Lagrangians are functions defined on $\jb^{1} Y$, we need to compute the first prolongation of the generating vector field.
The prolongation of $V$ is defined via the jet prolongation of the transformation map $\sigma^\eps$,
\begin{align}
\jp^{1} V = \dfrac{d}{d\eps} \Big[ \jp^{1} \sigma^\eps (\bx, \by) \Big] \bigg\vert_{\eps=0} ,
\end{align}
and is given by
\begin{align}\label{eq:noether_continuous_vector_field_prolongation}
\jp^{1} V
= \eta^{a} \, \dfrac{\partial}{\partial \by^{a}} + \bigg( \dfrac{\partial \eta^{a}}{\partial \bx^{\mu}} + \bz^{b}_{\mu} \dfrac{\partial \eta^{a}}{\partial \by^{b}} \bigg) \dfrac{\partial}{\partial \bz^{a}_{\mu}}
= \eta^{a} \, \dfrac{\partial}{\partial \by^{a}} + \eta^{a}_{\mu} \, \dfrac{\partial}{\partial \bz^{a}_{\mu}}
.
\end{align}
A vertical transformation $\sigma^\eps$ is a symmetry transformation for the Lagrangian (\ref{eq:variational_continuous_lagrangian}) if the invariance condition,
\begin{align}\label{eq:noether_continuous_symmetry_1}
L \big( \jp^{1} (\sigma^\eps \circ \phy) \big) = L \big( \jp^{1} \phy \big) ,
\end{align}
is satisfied\footnote{
Here, invariance of the Lagrangian suffices, but for general transformations, i.e., transformations that transform the base space in addition to the configuration space, invariance of the Lagrangian does not guarantee that $\phy^{\eps}$ is a solution whenever $\phy$ is a solution.
Invariance of solutions requires invariance of the action, and invariance of the action requires equivariance of the Lagrangian, which also takes the deformation of the integration domain $X$ into account.
For a detailed discussion see, e.g., \citet{Lew:2003}, section 6.5.}.
Taking the $\eps$ derivative of (\ref{eq:noether_continuous_symmetry_1}), we obtain an infinitesimal invariance condition,
\begin{align}\label{eq:noether_continuous_symmetry_2}
\dfrac{d}{d\eps} L \big( \jp^{1} (\sigma^\eps \circ \phy) \big) \bigg\vert_{\eps=0}
= 0 ,
\end{align}
which is equivalent to (\ref{eq:noether_continuous_symmetry_1}).
Explicitly computing the $\eps$ derivative, we obtain
\begin{align}\label{eq:noether_continuous_symmetry_3}
\jp^{1} V \big( L \big) \big( \jp^{1} \phy \big)
&= \dfrac{\partial L}{\partial \by^{a}} \big( \jp^{1} \phy \big) \cdot \eta^{a} (\phy) + \dfrac{\partial L}{\partial \bz^{a}_{\mu}} \big( \jp^{1} \phy \big) \cdot \bigg[ \dfrac{\partial \eta^{a}}{\partial \bx^{\mu}} + \dfrac{\partial \eta^{a}}{\partial \by^{b}} \dfrac{\partial \phy^{b}}{\partial \bx^{\mu}} \bigg]
= 0 .
\end{align}
If $\phy$ solves the Euler-Lagrange field equations (\ref{eq:variational_continuous_euler_lagrange}), we can replace the first term on the right-hand side of (\ref{eq:noether_continuous_symmetry_3}) to obtain
\begin{align}\label{eq:noether_continuous_symmetry_4}
\bigg[ \dfrac{\partial}{\partial \bx^{\mu}} \dfrac{\partial L}{\partial \bz^{a}_{\mu}} \big( \jp^{1} \phy \big) \bigg] \cdot \eta^{a} (\phy) + \dfrac{\partial L}{\partial \bz^{a}_{\mu}} \big( \jp^{1} \phy \big) \cdot \bigg[ \dfrac{\partial \eta^{a}}{\partial \bx^{\mu}} + \dfrac{\partial \eta^{a}}{\partial \by^{b}} \dfrac{\partial \phy^{b}}{\partial \bx^{\mu}} \bigg]
= 0 .
\end{align}
This, at last, amounts to a total divergence,
\begin{align}\label{eq:noether_continuous_noether_theorem}
\div \big[J (\jp^{1} \phy)\big] = 0 ,
\end{align}
with the Noether current $J$ given by
\begin{align}\label{eq:noether_continuous_current_definition}
J^{\mu} \big( \jp^{1} \phy \big) = \dfrac{\partial L}{\partial \bz^{a}_{\mu}} \big( \jp^{1} \phy \big) \cdot \eta^{a} (\phy) .
\end{align}
The fact that the Noether current is divergence-free, expresses the conservation law satisfied by solutions $\phy$ of the Euler-Lagrange field equations (\ref{eq:variational_continuous_euler_lagrange}).
The flux of $J$ through the boundary of any domain $\Omega \subseteq X$ is zero.

\subsubsection*{Global Form of Conservation Laws}

Let us consider an alternative point of view that is better suited for actual computations on the discrete level.
If the Lagrangian is invariant under the vertical transformation (\ref{eq:noether_continuous_transformation}), the action is also invariant under this transformation,
\begin{align}\label{eq:noether_continuous_global_invariance_condition}
\dfrac{d}{d\eps} \mathcal{A} [\sigma^\eps \circ \phy] \bigg\vert_{\eps=0}
= \int \limits_{X} \jp^{1} V \big( L \big) \big( \jp^{1} \phy \big) \, d^{n+1} \bx
= 0 .
\end{align}
Repeating the steps that lead from (\ref{eq:noether_continuous_symmetry_3}) to (\ref{eq:noether_continuous_noether_theorem}) this becomes
\begin{align}\label{eq:noether_continuous_conservation_law_1}
\dfrac{d}{d\eps} \mathcal{A} [\sigma^\eps \circ \phy] \bigg\vert_{\eps=0}
= \int \limits_{X} \dfrac{\partial}{\partial \bx^{\mu}} \bigg[ \dfrac{\partial L}{\partial \bz^{a}_{\mu}} \big( \jp^{1} \phy \big) \cdot \eta^{a} (\phy) \bigg] \, d^{n+1} \bx
= 0 .
\end{align}
With appropriate boundary conditions, one has
\begin{align}\label{eq:noether_continuous_current_spatial_integral}
\int \dfrac{\partial}{\partial \bx^{i}} \bigg[ \dfrac{\partial L}{\partial \bz^{a}_{i}} \big( \jp^{1} \phy \big) \cdot \eta^{a} (\phy) \bigg] \, d^{n} \bx = 0 ,
\end{align}
and thus, integrating (\ref{eq:noether_continuous_conservation_law_1}) for $\bx^{0} = \ct \in [t_{0}, t_{1}]$,
\begin{align}\label{eq:noether_continuous_conservation_law_2}
\dfrac{d}{d\eps} \mathcal{A} [\sigma^\eps \circ \phy] \bigg\vert_{\eps=0}
= \bigg[ \int \dfrac{\partial L}{\partial \bz^{a}_{t}} \big( \jp^{1} \phy \big) \cdot \eta^{a} (\phy) \, d^{n} \bx \bigg]_{t_{0}}^{t_{1}}
= 0 .
\end{align}
This is equivalent to integrating the divergence of the Noether current (\ref{eq:noether_continuous_noether_theorem}) over the spatial dimensions and using (\ref{eq:noether_continuous_current_spatial_integral}).
Since $t_{0}$ and $t_{1}$ are arbitrary, this implies the conservation of
\begin{align}\label{eq:noether_continuous_charge}
\mathcal{J}
= \int \dfrac{\partial L}{\partial \bz^{a}_{t}} \big( \jp^{1} \phy \big) \cdot \eta^{a} (\phy) \, d^{n} \bx
= \int J^{t} (\jp^{1} \phy) \, d^{n} \bx ,
\end{align}
which is called Noether charge.

\subsubsection*{Divergence Symmetries}
\label{sec:noether_divergence}

Sometimes it is necessary to consider a slightly more general version of Noether's theorem \cite{Olver:1993, Ibragimov:2007a}.
The invariance condition (\ref{eq:noether_continuous_symmetry_3}) can be weakened to
\begin{align}\label{eq:divergence_symmetries_symmetry_condition}
\jp^{1} V (L) (\jp^{1} \phy) = \div B (\jp^{1} \phy) ,
\end{align}
with $B$ being a vector field over $Y$ of the form $B = B^{\mu} (\bx, \by) \, \partial / \partial \bx^{\mu}$,
and $\div B$ is a scalar field defined over $\jb^1 Y$ given by
\begin{align}\label{eq:divergence_symmetries_div}
\div B (\bx, \by, \bz) = \dfrac{\partial B^{\mu} (\bx, \by)}{\partial \bx^{\mu}} + \dfrac{\partial B^{\mu} (\bx, \by)}{\partial \by^{a}} \, \bz^{a}_{\mu} .
\end{align}
This form of the invariance condition is tightly related to the gauge freedom of the Lagrangian. We can add any divergence to the Lagrangian without changing the equations of motion as
\begin{align}
\mathcal{A} [\phy]
= \int \limits_{X} L' (\jp^{1} \phy) \, d^{n+1} \bx 
= \int \limits_{X} \big( L + \div H \big) (\jp^{1} \phy) \, d^{n+1} \bx
= \int \limits_{X} L (\jp^{1} \phy) \, d^{n+1} \bx ,
\end{align}
assuming appropriate boundary conditions.
Here, $H$ is a vector field over $Y$ of the form $H = H^{\mu} (\bx,
\by)\,\partial / \partial \bx^{\mu}$, $\div H$ is defined as in (\ref{eq:divergence_symmetries_div}), and
\begin{align}\label{eq:divergence_symmetries_modified_lagrangian}
L' (\bx, \by, \bz)
= \big( L + \div H \big) (\bx, \by, \bz) .
\end{align}
One can check that the Euler-Lagrange equations for $L$ and $L'$ are the same.
Requiring the invariance condition (\ref{eq:noether_continuous_symmetry_3}) to be satisfied for $L'$,
\begin{align}
\jp^{1} V (L') (\jp^{1} \phy) = \jp^{1} V (L) (\jp^{1} \phy) + \div \tilde{H} (\jp^{1} \phy) = 0 ,
\quad \text{with} \quad
\tilde{H} = V(H^\mu) \frac{\partial}{\partial x^\mu},
\end{align}
we find the following relation between the vector fields $B$ and $H$,
\begin{align}\label{eq:divergence_symmetries_vector_field_H}
B = - \tilde{H}.
\end{align}
In summary, if $V$ is a divergence symmetry for $L$ in the sense of
(\ref{eq:divergence_symmetries_symmetry_condition}) and there exists a vector
field $H$ satisfying (\ref{eq:divergence_symmetries_vector_field_H}), then $V$ is a Lie-point symmetry (\ref{eq:noether_continuous_symmetry_3}) of the equivalent Lagrangian $L'$ defined in (\ref{eq:divergence_symmetries_modified_lagrangian}).
We can therefore apply the Noether theorem to $L'$, obtaining the Noether current 
\begin{align}\label{eq:divergence_symmetries_noether_continuous_current}
J' (\jp^{1} \phy) = J (\jp^{1} \phy) - B (\phy) ,
\end{align}
with $J'$ the conserved Noether current associated to the divergence symmetry (\ref{eq:divergence_symmetries_symmetry_condition}) and $J$ as defined in (\ref{eq:noether_continuous_current_definition}), which is not conserved in this case.

\subsection{Discrete Jet Bundles}
\label{sec:vi_discrete_jet_bundles}

In order to derive a discrete version of Hamilton's action principle, we proceed in three steps.
In this section, we define discrete analogues of the base space $X$, the configuration space $Y$, and the jet bundles $\jb^{1} Y$ and  $\jb^{2} Y$.
In the next section, we discretise the Lagrangian and the action functional, and finally, we work out the discrete action principle.

\begin{figure}[tb]
	\centering
	\subfloat[Rectangular Grid]{\label{fig:grid_types_rectangular}
		\includegraphics[height=3.5cm]{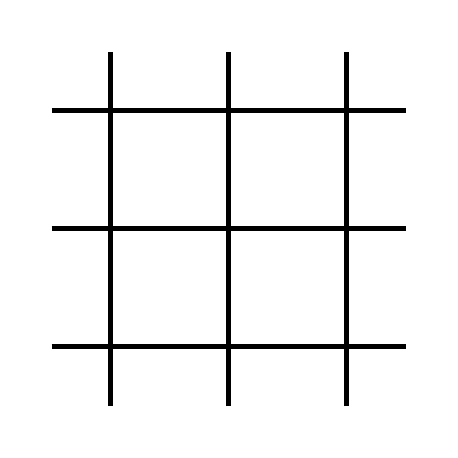}
	}
	\subfloat[Triangular Grid]{\label{fig:grid_types_triangular}
		\includegraphics[height=3.5cm]{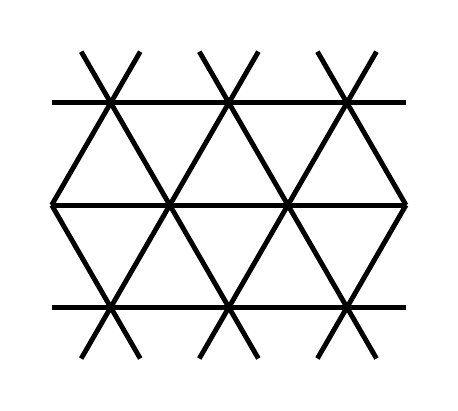}
	}
	\subfloat[Hexagonal Grid]{\label{fig:grid_types_hexagonal}
		\includegraphics[height=3.5cm]{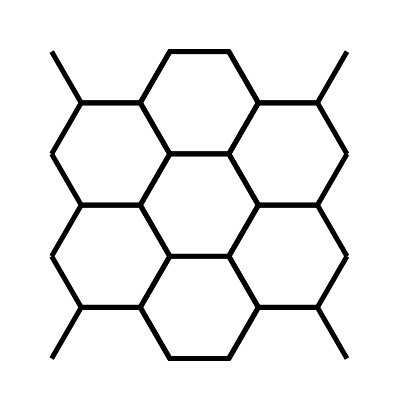}
	}
	\caption{Different regular grids.}
	\label{fig:grid_types}
\end{figure}

The continuous base space $X$ is replaced by its discrete analogue $X_{d}$, a bounded subset $X_{d} \subset \mathbb{Z}^{n+1}$ which corresponds to a grid of points in $X$.
In what follows we assume $X$ to be two-dimensional and choose an equidistant rectangular discretisation like it is depicted in Figure~\subref*{fig:grid_types_rectangular}. The theory is easily applicable to other types of grids as well, e.g., triangular as in Figure~\subref*{fig:grid_types_triangular}, hexagonal as in Figure~\subref*{fig:grid_types_hexagonal}, or even staggered or irregular ones.

\begin{figure}[t]
	\centering
	\subfloat[Grid Points]{\label{fig:gridbox_indices}
		\includegraphics[height=3.5cm]{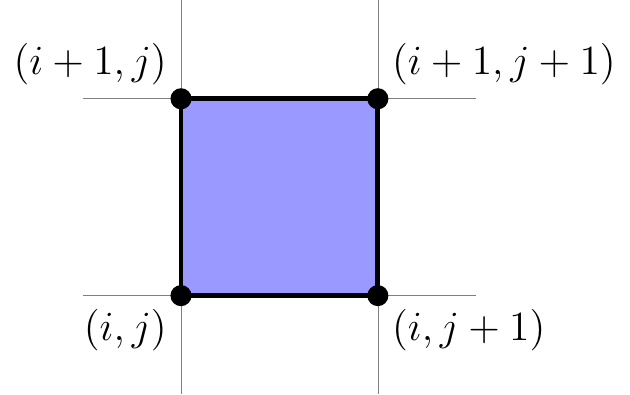}
	}
	\hspace{1cm}
	\subfloat[Field Components]{\label{fig:gridbox_fields}
		\includegraphics[height=3.5cm]{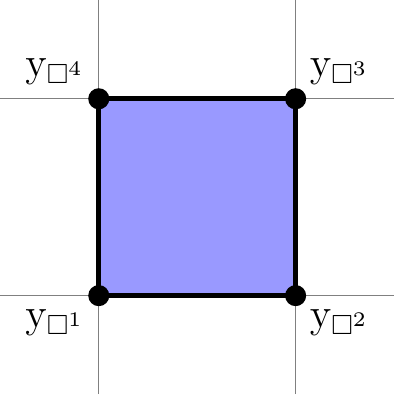}
	}
	\caption{Vertices of a primal grid cell in a two-dimensional rectangular grid and field components at those vertices.}
	\label{fig:gridbox}
\end{figure}

The coordinates on the discrete base space $X_{d} \subset \mathbb{Z} \times \mathbb{Z}$ are denoted $(i,j)$ and assumed to take values $i \in \{ 1, ..., N_{0} \}$ and $j \in \{ 1, ..., N_{1} \}$, respectively, where $N_{0}$ and $N_{1}$ are the number of points for each dimension.
This corresponds to a grid of points $\bx_{i,j} = (ih, jh)$ in $X$ with, for simplicity, the same step size $h$ in both directions.
Then
\begin{align}\label{eq:vi_discrete_base_space}
X_{d} &\cong \big\{ (ih, jh) \in X \; \big\vert \; i=1, ..., N_{0}, \, j=1, ..., N_{1} \big\} .
\end{align}
The discrete configuration space is defined as the cartesian product
\begin{align}\label{eq:vi_discrete_configuration_space}
Y_{d} = X_{d} \times F ,
\end{align}
where $F$ is the same as in the continuous case and we have an analogous projection
\begin{align}
\arraycolsep=2pt
\begin{array}{rcccc}
\pi_{d} & : & Y_{d} & \rightarrow & X_{d} , \\
& & (\by_{i,j}) & \mapsto & (i,j) .
\end{array}
\end{align}
Coordinates of $Y_{d}$ are denoted $\by^{a}_{i,j}$ with $1 \leq a \leq \dim F$. The coordinates of the base point $(i,j)$ are already implied and therefore not specified separately.
While coordinates on $X_{d}$ and $Y_{d}$ are defined point-wise, coordinates on $\jb^{1} Y_{d}$ will be defined grid-cell-wise.
We therefore introduce the following abstract but convenient notation.
A square $\square$ on $X_{d}$ is an ordered quadruplet
\begin{align}\label{eq:vi_discrete_field_primal_grid_cell}
\square = \big( (i, j), (i, j+1), (i+1, j+1), (i+1, j) \big) ,
\end{align}
defining a primal grid cell.
The set of such cells on $X_{d}$ is denoted $X^{\square}$.
Vertices $\square^{l}$ of a square $\square$ with $1 \leq l \leq 4$ are counted counter-clockwise from the bottom left (c.f. Figure~\subref*{fig:gridbox_indices}), namely,
\begin{align}\label{eq:vi_discrete_vertices}
\square^{1} &= (i,   j  ) , &
\square^{2} &= (i,   j+1) , &
\square^{3} &= (i+1, j+1) , &
\square^{4} &= (i+1, j  ) .
\end{align}
On a quadrilateral grid, the first jet bundle of $Y_{d}$ becomes \cite{MarsdenPatrick:1998}
\begin{align}\label{eq:vi_discrete_jet_space}
\jb^{1} Y_{d}
&= X^{\square} \times F^{4} ,
\end{align}
with coordinates given by
\begin{align}\label{eq:vi_discrete_jet_bundle_coordinates}
\big\{ ( \square, \, \by_{\square^{l}}^{a} ) \; \big\vert \; \square \in X^{\square} , \; \by_{\square^{l}}^{a} \in \mathbb{R} , \; 1 \leq l \leq 4 , \; 1 \leq a \leq \dim F \big\} .
\end{align}
This implies that if $\square^{1} = (i,j)$, then we have
\begin{align}\label{eq:vi_discrete_coordinates_vertices}
\by_{\square^{1}} &= \by_{i,  j  } , &
\by_{\square^{2}} &= \by_{i,  j+1} , &
\by_{\square^{3}} &= \by_{i+1,j+1} , &
\by_{\square^{4}} &= \by_{i+1,j  } .
\end{align}
A section of $Y_{d}$, representing a discrete field, is a map $\phy_{d} : X_{d} \rightarrow Y_{d}$ such that $\pi \circ \phy_{d} = \id_{X_{d}}$.
The components of a discrete field $\phy_{d}$ at the vertices (\ref{eq:vi_discrete_vertices}) of a specific square are
\begin{align}\label{eq:vi_discrete_field_vertices}
\phy_{\square^{1}} &= \phy_{d} (\square^{1}) , &
\phy_{\square^{2}} &= \phy_{d} (\square^{2}) , &
\phy_{\square^{3}} &= \phy_{d} (\square^{3}) , &
\phy_{\square^{4}} &= \phy_{d} (\square^{4}) ,
\end{align}
as depicted in Figure~\subref*{fig:gridbox_fields}.
The discrete first jet prolongation of a discrete field $\phy_{d}$ is a map $\jpd^{1} \phy_{d} : X^{\square} \rightarrow \jb^{1} Y_{d}$, defined as
\begin{align}\label{eq:vi_discrete_prolongation_field}
\jpd^{1} \phy_{d} (\square) = \big( \square, \phy_{\square^{1}}, \phy_{\square^{2}}, \phy_{\square^{3}}, \phy_{\square^{4}} \big) .
\end{align}
This contains all the information necessary to define discrete first-order derivatives.

We will see that the discrete Euler-Lagrange equations are living on $\jb^{2} Y_{d}$, i.e., the second jet bundle of $Y_{d}$, similar to the continuous case (see \cite{KouranbaevaShkoller:2000} for another discussion of discrete second order jet bundles).
Coordinates on $\jb^{2} Y_{d}$ will be defined on a quadruple of cells $\boxplus$, namely those cells $\square$ in $X^{\square}$ which share a common vertex.
The cells $\square_{m}$ of $\boxplus$ with $1 \leq m \leq 4$ are counted counter-clockwise from the bottom left, so that
\begin{align}\label{eq:eq:vi_discrete_definition_quadruple}
\boxplus = \big\{ \square_{1}, \square_{2}, \square_{3}, \square_{4} \in X^{\square} \; \big\vert \; \square_{1}^{3} = \square_{2}^{4} = \square_{3}^{1} = \square_{4}^{2} \big\} . 
\end{align}
The set of all such quadruples on $X_{d}$ is denoted by $X^{\boxplus}$.
Vertices $\boxplus^{l}$ of a quadruple $\boxplus$ with $1 \leq l \leq 9$ are counted row-wise from bottom left to top right, namely,
\begin{align*}
\boxplus^{1} &= (i-1, j-1) , &
\boxplus^{2} &= (i-1, j  ) , &
\boxplus^{3} &= (i-1, j+1) , \\
\boxplus^{4} &= (i,   j-1) , &
\boxplus^{5} &= (i,   j  ) , &
\boxplus^{6} &= (i,   j+1) , \\
\boxplus^{7} &= (i+1, j-1) , &
\boxplus^{8} &= (i+1, j  ) , &
\boxplus^{9} &= (i+1, j+1) .
\end{align*}
On a quadrilateral grid, the second jet bundle of $Y_{d}$ can be identified with
\begin{align}\label{eq:vi_discrete_jet_space2}
\jb^{2} Y_{d}
&= X^{\boxplus} \times F^{9} ,
\end{align}
with coordinates given by
\begin{align}
\big\{ ( \boxplus, \, \by_{\boxplus^{l}}^{a} ) \; \big\vert \; \boxplus \in X^{\boxplus} , \; \by_{\boxplus^{l}}^{a} \in \mathbb{R} , \; 1 \leq l \leq 9 , \; 1 \leq a \leq \dim F \big\} ,
\end{align}
in analogy with (\ref{eq:vi_discrete_jet_bundle_coordinates}).
The second jet prolongation of a discrete field $\phy_{d}$ is a map $\jpd^{2} \phy_{d} : X^{\boxplus} \rightarrow \jb^{2} Y_{d}$, defined as
\begin{align}\label{eq:vi_discrete_prolongation_field2}
\jpd^{2} \phy_{d} (\boxplus) = \big( \boxplus, \phy_{\boxplus^{1}} , \phy_{\boxplus^{2}}, \phy_{\boxplus^{3}}, \phy_{\boxplus^{4}}, \phy_{\boxplus^{5}}, \phy_{\boxplus^{6}}, \phy_{\boxplus^{7}}, \phy_{\boxplus^{8}}, \phy_{\boxplus^{9}} \big) .
\end{align}
Similar to the continuous case, $\jb^{2} Y_{d}$ can also be defined via iteration as the first jet bundle of $\jb^{1} Y_{d}$. If we write coordinates on $\jb^{1} ( \jb^{1} Y_{d} )$ as
\begin{align}
\big\{ ( \boxplus , \, \square_{m} , \, \by_{\square_{m}^{l}}^{a} ) \; \big\vert \; \boxplus \in X^{\boxplus} , \; \square_{m} \in \boxplus , \; \by_{\square_{m}^{l}}^{a} \in \mathbb{R} , \; 1 \leq l, m \leq 4 , \; 1 \leq a \leq \dim F \big\} ,
\end{align}
then $\jb^{2} Y_{d}$ consists of all elements of $\jb^{1} ( \jb^{1} Y_{d} )$, which satisfy
\begin{align}
   \by_{\square_{1}^{3}}^{a}
&= \by_{\square_{2}^{4}}^{a}
 = \by_{\square_{3}^{1}}^{a}
 = \by_{\square_{4}^{2}}^{a} , &
\by_{\square_{1}^{2}}^{a} &= \by_{\square_{2}^{1}}^{a} , &
\by_{\square_{2}^{3}}^{a} &= \by_{\square_{3}^{2}}^{a} , &
\by_{\square_{3}^{4}}^{a} &= \by_{\square_{4}^{3}}^{a} , &
\by_{\square_{4}^{1}}^{a} &= \by_{\square_{1}^{4}}^{a} , 
\end{align}
for all $a$, which in our construction is always guaranteed due to (\ref{eq:vi_discrete_coordinates_vertices}) and (\ref{eq:eq:vi_discrete_definition_quadruple}).
Therefore, $\jb^{2} Y_{d}$ can be identified with $\jb^{1} ( \jb^{1} Y_{d} )$, unlike the continuous case, where $\jb^{2} Y$ is strictly embedded into $\jb^{1} ( \jb^{1} Y )$.

\subsection{Discrete Action Principle}
\label{sec:vi_discrete_action_principle}

The discretisation of the Lagrangian is based on the observation that the continuous action functional can be written as
\begin{align}\label{eq:vi_discrete_action_sum}
\mathcal{A} [\phy]
= \sum \limits_{\square \in X^{\square}} \int \limits_{\text{Vol}(\square)} L \big( \jp^{1} \phy \big) \, d^{n+1} \bx ,
\end{align}
where $\text{Vol}(\square) \subset X$ is the physical domain enclosed by $\square$.
The integral in (\ref{eq:vi_discrete_action_sum}) is approximated by a function of values of $\phy_{d}$ in four different points in the spacetime grid, which corresponds to the discrete jet prolongation defined in (\ref{eq:vi_discrete_prolongation_field}). This is the discrete Lagrangian $L_{d}$, i.e.,
\begin{align}\label{eq:vi_discrete_replace_lagrangian}
\int \limits_{\text{Vol}(\square)} L \big( \jp^{1} \phy \big) \, d^{n+1} \bx
\approx
L_{d} \big( \jpd^{1} \phy_{d} (\square) \big) =
L_{d} \big( \square, \phy_{\square^1}, \phy_{\square^2}, \phy_{\square^3}, \phy_{\square^4} \big) .
\end{align}
The action functional (\ref{eq:vi_discrete_action_sum}) is then approximated by
\begin{align}\label{eq:vi_discrete_replace_action}
\mathcal{A}_{d} [\phy_{d}] = \sum \limits_{\square \in X^{\square}} L_{d} \big( \jpd^{1} \phy_{d} (\square) \big) ,
\end{align}
which can also be written explicitly as
\begin{align}
\mathcal{A}_{d} [\phy_{d}] = \sum \limits_{i=1}^{N_{0}-1} \sum \limits_{j=1}^{N_{1}-1} L_{d} \big( \phy_{i,j}, \phy_{i,j+1}, \phy_{i+1,j+1}, \phy_{i+1,j} \big) ,
\end{align}
but for most of our derivations the abstract notation is more practical.

Specifically, the discrete Lagrangian is obtained by introducing a quadrature rule (e.g., trapezoidal, midpoint, Simpson) to approximate the integral in (\ref{eq:vi_discrete_replace_lagrangian}) as well as approximations of the fields and their derivatives.
In the spirit of the Veselov discretisation \cite{Veselov:1988, Veselov:1991, MoserVeselov:1991, MarsdenPatrick:1998}, we will use the midpoint rule and first-order finite differences. This entails that the continuous fields in the Lagrangian are replaced with averages of field values at the four vertices of the grid cell $\square$, i.e.,
\begin{align}\label{eq:vi_discrete_average}
\by^{a} \rightarrow \dfrac{1}{4} \Big( \by_{\square^1}^{a} + \by_{\square^2}^{a} + \by_{\square^3}^{a} + \by_{\square^4}^{a} \Big) \equiv \overline{\by}^{a} (\square).
\end{align}
There are two possibilities for the definition of each of the derivatives.
With reference to Figure~\subref*{fig:gridbox_fields}, let the vertical dimension correspond to $\bx^{0}$ (time) and the horizontal dimension to $\bx^{1}$ (space).
Then $\partial_{0}$ can be defined along the left as well as along the right edge of the grid cell.
Similarly, $\partial_{1}$ can be defined along the upper as well as along the lower edge. 
Best results are usually obtained for the most symmetric discretisation of the Lagrangian. Therefore we use the average of the respective options and replace the derivatives in the Lagrangian according to
\begin{align}
\label{eq:vi_discrete_derivative_0}
\bz^{a}_{0} &\rightarrow \dfrac{1}{2} \bigg( \dfrac{\by_{\square^4}^{a} - \by_{\square^1}^{a}}{h} + \dfrac{\by_{\square^3}^{a} - \by_{\square^2}^{a}}{h} \bigg) \equiv \overline{\bz}_{0}^{a} (\square) , \\
\label{eq:vi_discrete_derivative_1}
\bz^{a}_{1} &\rightarrow \dfrac{1}{2} \bigg( \dfrac{\by_{\square^2}^{a} - \by_{\square^1}^{a}}{h} + \dfrac{\by_{\square^3}^{a} - \by_{\square^4}^{a}}{h} \bigg) \equiv \overline{\bz}_{1}^{a} (\square) .
\end{align}

At last, in order to obtain the discrete equations of motion, we only need to apply a discrete version of Hamilton's principle of stationary action.
The vertical transformation $\sigma^\eps$ is discretised in the same way as the
fields $\phy$, by considering its values over points $(i,j)$ of the discrete
base space $X_{d}$, i.e., we replace $\sigma^\eps (\bx, \by)$ with
$\sigma_d^\eps = \big(\sigma^\eps_{i,j} (\by_{i,j}) \big)$.
The generating vector field is computed as in (\ref{eq:variational_continuous_vector_field}), but on grid points only,
\begin{align}\label{eq:vi_discrete_vector_field}
\eta_{i,j}^{a} = \dfrac{d(\sigma^\eps_{i,j})^{a}}{d\eps} \bigg\vert_{\eps=0} .
\end{align}
The discrete version of the action principle (\ref{eq:variational_continuous_variation_1}) then reads
\begin{align}\label{eq:vi_discrete_variation_1}
\dfrac{d}{d\eps} \mathcal{A}_{d} [ \sigma_{d}^\eps \circ \phy_{d} ] \bigg\vert_{\eps=0} = 0
\quad
\text{for every $\sigma_{d}^\eps$} .
\end{align}
The explicit computation of (\ref{eq:vi_discrete_variation_1}) leads to
\begin{align}
\nonumber
\dfrac{d}{d\eps} \mathcal{A}_{d} [ \sigma^\eps_{d} \circ \phy_{d} ] \bigg\vert_{\eps=0}
&= \sum \limits_{\square \in X^{\square}} \dfrac{d}{d\eps} L_{d} \big( \jpd^{1} \phy_{d}^{\eps} (\square) \big) \bigg\vert_{\eps=0} \\
\label{eq:vi_discrete_variation_3}
&= \sum \limits_{\square \in X^{\square}} \sum \limits_{l=1}^{4} \dfrac{\partial L_{d}}{\partial \by_{\square^{l}}^{a}} \big( \jpd^{1} \phy_{d} (\square) \big) \cdot \eta_{\square^{l}}^{a}  (\phy_{\square^{l}}) .
\end{align}
As the variation of the action has to vanish for each $\eta_{i,j}^{a}$ on the spacetime grid independently, it is sufficient to consider only those contributions that are multiplied with the vector field at a fixed grid point $(i,j)$.
\begin{figure}[tb]
	\centering
	\includegraphics[width=.5\textwidth]{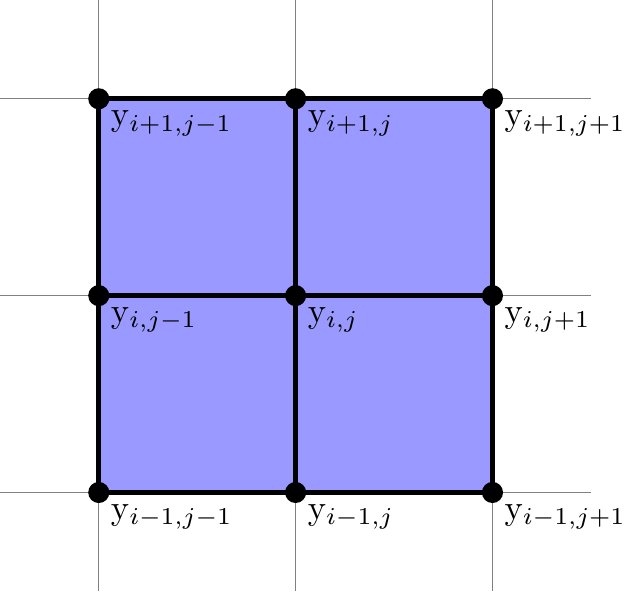}
	\caption{Contributions to the discrete Euler-Lagrange field equations. The derivatives of the discrete Lagrangian at the highlighted grid cells with respect to $\by_{i,j}$ correspond to the terms in equation (\ref{eq:vi_discrete_euler_lagrange_equations}).}
	\label{fig:variation}
\end{figure}
In total there are four such contributions (see Figure~\ref{fig:variation}),
\begin{align*}
& \dfrac{d}{d\eps} \mathcal{A}_{d} [ \sigma^\eps_{d} \circ \phy_{d} ] \bigg\vert_{\eps=0} = \\
& \hspace{2em}
   ... + \dfrac{\partial L_d}{\partial \by_{\square^1}^{a}} \Big( \phy_{i,  j  }, \phy_{i,  j+1}, \phy_{i+1,j+1}, \phy_{i+1,j  } \Big) \cdot \eta_{i,j}^{a} (\phy_{i,j}) + ... \\
& \hspace{3.5em}
   ... + \dfrac{\partial L_d}{\partial \by_{\square^2}^{a}} \Big( \phy_{i,  j-1}, \phy_{i,  j  }, \phy_{i+1,j  }, \phy_{i+1,j-1} \Big) \cdot \eta_{i,j}^{a} (\phy_{i,j}) + ... \\
& \hspace{5em}
   ... + \dfrac{\partial L_d}{\partial \by_{\square^3}^{a}} \Big( \phy_{i-1,j-1}, \phy_{i-1,j  }, \phy_{i,  j  }, \phy_{i,  j-1} \Big) \cdot \eta_{i,j}^{a} (\phy_{i,j}) + ... \\
& \hspace{6.5em}
   ... + \dfrac{\partial L_d}{\partial \by_{\square^4}^{a}} \Big( \phy_{i-1,j  }, \phy_{i-1,j+1}, \phy_{i,  j+1}, \phy_{i,  j  } \Big) \cdot \eta_{i,j}^{a} (\phy_{i,j}) + ...
= 0 .
\end{align*}
The variation of the discrete action vanishes, if the sum of all expressions multiplying $\eta_{i,j}^{a}$ vanishes identically for all $(i,j)$ (since the $(\sigma^\eps_{i,j})^{a}$ and therefore the $\eta_{i,j}^{a}$ are arbitrary). For each $a$, this requirement yields the discrete Euler-Lagrange field equations at $(i,j)$,
\begin{align}\label{eq:vi_discrete_euler_lagrange_equations}
\nonumber
0= & {\dfrac{\partial L_d}{\partial \by_{\square^1}^{a}} \Big( \phy_{i,  j  }, \phy_{i,  j+1}, \phy_{i+1,j+1}, \phy_{i+1,j  } \Big)}
 +   {\dfrac{\partial L_d}{\partial \by_{\square^2}^{a}} \Big( \phy_{i,  j-1}, \phy_{i,  j  }, \phy_{i+1,j  }, \phy_{i+1,j-1} \Big)} \\
 + & {\dfrac{\partial L_d}{\partial \by_{\square^3}^{a}} \Big( \phy_{i-1,j-1}, \phy_{i-1,j  }, \phy_{i,  j  }, \phy_{i,  j-1} \Big)}
 +   {\dfrac{\partial L_d}{\partial \by_{\square^4}^{a}} \Big( \phy_{i-1,j  }, \phy_{i-1,j+1}, \phy_{i,  j+1}, \phy_{i,  j  } \Big)} ,
\end{align}
which can be compactly written as an equation on the discrete second jet bundle $\jb^{2} Y_{d}$, namely
\begin{align}\label{eq:vi_discrete_euler_lagrange_equations_compact}
\big( \jpd^{2} \phy_{d} (\boxplus) \big)^{*} \big( \del(L_{d}) \big)^{a} =
\sum \limits_{\substack{l=1\\ \square^{l} = \boxplus^{5}}}^{4}
\dfrac{\partial L_{d}}{\partial \by_{\square^l}^{a}} \big( \jpd^{1} \phy_{d} (\square) \big) &= 0 &
& \text{for all $a$ and all $\boxplus$} .
\end{align}
Here, $\del(L_{d})$ is the discrete Euler-Lagrange operator acting on the discrete Lagrangian $L_{d}$, which is evaluated on the prolongation $\jpd^{2} \phy_{d} (\boxplus)$ as indicated by the pull-back notation.
These relations define the variational integrator for a first-order Lagrangian field theory according to the Veselov discretisation of the Lagrangian as it was described above.

\subsection{Discrete Noether Theorem}
\label{sec:vi_discrete_noether_theorem}

Following the derivation of the continuous theory from the previous section, we consider vertical transformations $\sigma^\eps_{i,j} (\by_{i,j})$ and define
\begin{align}\label{eq:noether_discrete_transformation}
\phy_{i,j}^{\eps} &= \sigma^\eps_{i,j} \circ \phy_{i,j} &
& \text{with} &
\sigma^\eps_{i,j} \vert_{\eps=0} &= \mathrm{id} &
& \text{and} &
\eta_{i,j}^{a} &= \dfrac{d (\sigma^\eps_{i,j})^{a}}{d\eps} \bigg\vert_{\eps = 0} .
\end{align}
The transformation $\sigma^\eps_{i,j}$ is a symmetry for the discrete Lagrangian (\ref{eq:vi_discrete_replace_lagrangian}) if
\begin{align}\label{eq:noether_discrete_symmetry_condition_1}
L_{d} \big( \jpd^{1} \phy^{\eps}_{d} (\square) \big) &= L_{d} \big( \jpd^{1} \phy_{d} (\square) \big) ,
\end{align}
which is equivalent to the infinitesimal symmetry condition
\begin{align}\label{eq:noether_discrete_symmetry_condition_2}
\dfrac{d}{d \eps} L_{d} \big( \jpd^{1} \phy^{\eps}_{d} (\square) \big) \bigg\vert_{\eps = 0}
= \sum \limits_{l=1}^{4} \dfrac{\partial L_{d}}{\partial \by^{a}_{\square^l}} \big( \jpd^{1} \phy_{d} (\square) \big) \cdot \eta^{a}_{\square^{l}} \big( \phy_{\square^{l}} \big) 
= 0 .
\end{align}
On each grid cell $\square$ we can define four discrete momentum maps in analogy to (\ref{eq:noether_continuous_current_definition}),
\begin{align}\label{eq:noether_discrete_current_abstract}
\big( \jpd^{1} \phy_{d} (\square) \big)^{*} J_{\square^{l}}
&= \dfrac{\partial L_d}{\partial \by^{a}_{\square^l}} \big( \jpd^{1} \phy_{d} (\square) \big) \cdot \eta^{a}_{\square^{l}} \big( \phy_{\square^{l}} \big) , &
l &\in \{ 1, 2, 3, 4 \} .
\end{align}
Instead of two components of $J$, corresponding to the coordinates of the base space, we now have four contributions, corresponding to the vertices of a grid cell.
With that, we can write the discrete symmetry condition (\ref{eq:noether_discrete_symmetry_condition_2}) as
\begin{align}\label{eq:noether_discrete_symmetry_condition_3}
J_{\square^{1}} + J_{\square^{2}} + J_{\square^{3}} + J_{\square^{4}} = 0, \quad 
\text{for every $\square \in X^\square$.}
\end{align}
In addition to (\ref{eq:noether_discrete_current_abstract}), we define some explicit shorthand notation that will become useful in the following derivation,
\begingroup
\allowdisplaybreaks
\begin{subequations}\label{eq:noether_discrete_current_shorthand}
\begin{align}
J_{i,j}^{1} &= \dfrac{\partial L_d}{\partial \by^{a}_{\square^1}} \Big( \phy_{i,  j  }, \phy_{i,  j+1}, \phy_{i+1,j+1}, \phy_{i+1,j  } \Big) \cdot \eta^{a}_{i,  j  } ( \phy_{i,  j  } ),
\\
J_{i,j}^{2} &= \dfrac{\partial L_d}{\partial \by^{a}_{\square^2}} \Big( \phy_{i,  j  }, \phy_{i,  j+1}, \phy_{i+1,j+1}, \phy_{i+1,j  } \Big) \cdot \eta^{a}_{i,  j+1} ( \phy_{i,  j+1} ) ,
\\
J_{i,j}^{3} &= \dfrac{\partial L_d}{\partial \by^{a}_{\square^3}} \Big( \phy_{i,  j  }, \phy_{i,  j+1}, \phy_{i+1,j+1}, \phy_{i+1,j  } \Big) \cdot \eta^{a}_{i+1,j+1} ( \phy_{i+1,j+1} ) ,
\\
J_{i,j}^{4} &= \dfrac{\partial L_d}{\partial \by^{a}_{\square^4}} \Big( \phy_{i,  j  }, \phy_{i,  j+1}, \phy_{i+1,j+1}, \phy_{i+1,j  } \Big) \cdot \eta^{a}_{i+1,j  } ( \phy_{i+1,j  } ).
\end{align}
\end{subequations}
\endgroup

\begin{figure}[tb]
	\centering
	\includegraphics[width=\textwidth]{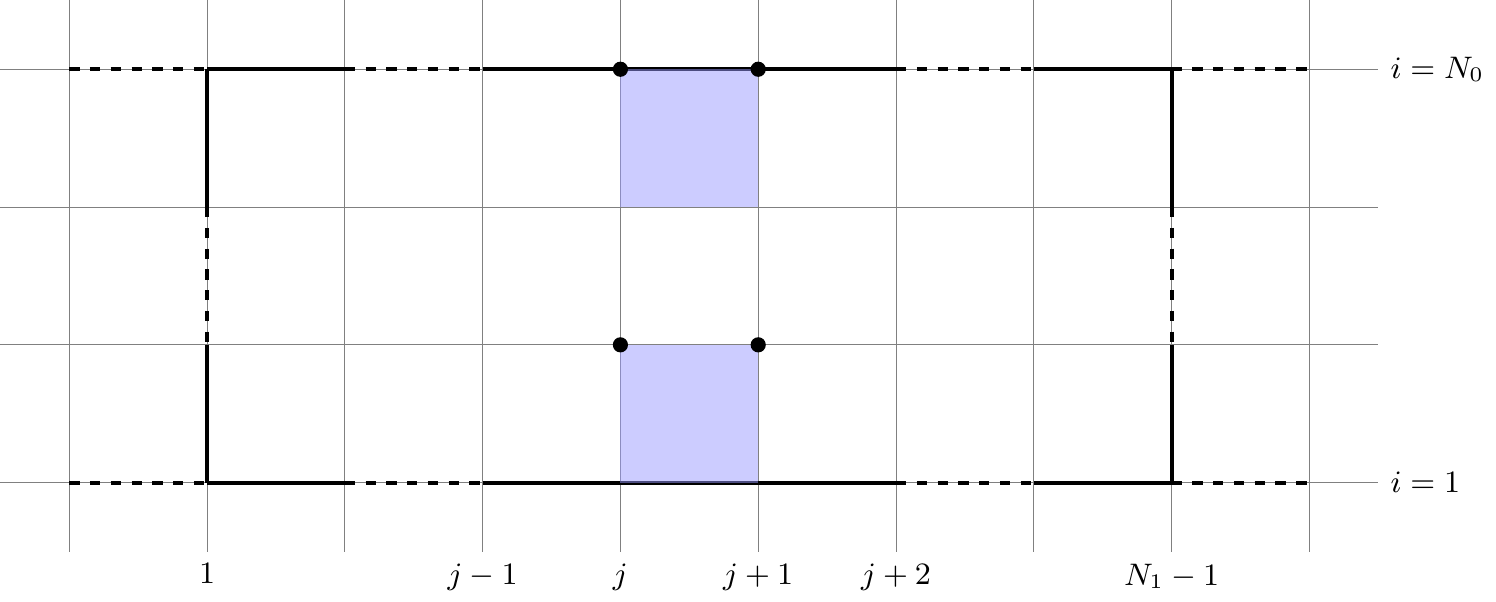}
	\caption{Contributions to the discrete Noether theorem.}
	\label{fig:discrete_noether_theorem}
\end{figure}

In order to derive the discrete conservation law, we follow the argument from
section \ref{sec:vi_continous_noether_theorem}.
If the Lagrangian is invariant under a given vertical transformation, then the
action is also invariant under this transformation, namely,
\begin{align}
\dfrac{d}{d \eps} \mathcal{A}_{d} \big[ \phy_{d}^{\eps} \big] \bigg\vert_{\eps = 0} = 0 ,
\end{align}
which becomes
\begin{align}\label{eq:noether_discrete_symmetry_condition_global}
\dfrac{d}{d \eps} \mathcal{A}_{d} \big[ \phy_{d}^{\eps} \big] \bigg\vert_{\eps = 0}
\nonumber
&= \sum \limits_{\square \in X^{\square}} \sum \limits_{l=1}^{4} \dfrac{\partial L_{d}}{\partial \by^{a}_{\square^l}} \big( \jpd^{1} \phy_{d} (\square) \big) \cdot \eta^{a}_{\square^{l}} \big( \phy_{\square^{l}} \big) \\
&= \sum \limits_{\square \in X^{\square}} \sum \limits_{l=1}^{4} \big( \jpd^{1} \phy_{d} (\square) \big)^{*} J_{\square^{l}}
 = 0 .
\end{align}
Those contributions to the sum that originate from interior points vanish in virtue of the discrete Euler-Lagrange field equations (\ref{eq:vi_discrete_euler_lagrange_equations}).
Specifically,
\begin{align}
  J^{1}_{i,  j  }
+ J^{2}_{i,  j-1}
+ J^{3}_{i-1,j-1}
+ J^{4}_{i-1,j  }
= \big( \jpd^{2} \phy_{d} (\boxplus) \big)^{*} \big( \del(L_{d}) \big)_{i,j}^{a} \cdot \eta^{a}_{i,j} (\phy_{i,j})
= 0 ,
\end{align}
so that only boundary cells contribute to (\ref{eq:noether_discrete_symmetry_condition_global}), c.f. Figure~(\ref{fig:discrete_noether_theorem}).
For simplicity, we will assume periodic
boundary conditions in the space-like variable, then the only contribution comes from the boundary of the time-like variable,
\begin{align}\label{eq:noether_discrete_conservation_law_1}
0 = \sum \limits_{j=1}^{N_{1}-1} \Big[
  J^{1}_{1,j  }
+ J^{2}_{1,j-1}
+ J^{3}_{N_{0}-1,j-1}
+ J^{4}_{N_{0}-1,j  }
\Big] .
\end{align}
Using the discrete symmetry condition (\ref{eq:noether_discrete_symmetry_condition_3}) to replace the first two terms as well as the periodicity of the spatial domain, we obtain the discrete counterpart of (\ref{eq:noether_continuous_conservation_law_2}), that is the discrete conservation law
\begin{align}\label{eq:noether_discrete_conservation_law_2}
  \sum \limits_{j=1}^{N_{1}-1} \Big[ J^{3}_{1,j}       + J^{4}_{1,j}       \Big]
= \sum \limits_{j=1}^{N_{1}-1} \Big[ J^{3}_{N_{0}-1,j} + J^{4}_{N_{0}-1,j} \Big] ,
\end{align}
see Figure~\ref{fig:discrete_noether_theorem}.
Since the number of timesteps $N_{0}$ is arbitrary, this is equivalent to
\begin{align}\label{eq:noether_discrete_conservation_law_3}
\mathcal{J}_{d} &= \sum \limits_{\substack{j=1\\ \square^{1} = (i,j)}}^{N_{1}-1} \Big[ J_{\square^{3}} \big( \jpd^{1} \phy_{d} (\square) \big) + J_{\square^{4}} \big( \jpd^{1} \phy_{d} (\square) \big) \Big] = \text{const. for all $i$} ,
\end{align}
in analogy to (\ref{eq:noether_continuous_charge}).
This proves that variational integrators preserve discrete invariants $\mathcal{J}_{d}$ to machine accuracy\footnote{
In practice, many variational integrators require the solution of nonlinear algebraic equations. In that case, the accuracy to which the invariants are preserved will generally depend on the accuracy of the nonlinear solver.
}.

\subsubsection*{Discrete Divergence Symmetries}

We can generalise the discrete Noether theorem in a similar way as we did with the continuous Noether theorem for divergence symmetries.
In analogy to (\ref{eq:divergence_symmetries_symmetry_condition}), we rewrite the discrete symmetry condition (\ref{eq:noether_discrete_symmetry_condition_3}) as
\begin{align}\label{eq:noether_discrete_divergence_symmetry_condition}
  J_{\square^{1}}
+ J_{\square^{2}}
+ J_{\square^{3}} 
+ J_{\square^{4}}
= B_{\square^{1}}
+ B_{\square^{2}}
+ B_{\square^{3}} 
+ B_{\square^{4}} .
\end{align}
Making no assumption on the $B_{\square^{l}}$, equation (\ref{eq:noether_discrete_conservation_law_1}) takes the form
\begin{multline}\label{eq:noether_discrete_divergence_conservation_law_1}
\sum \limits_{j=1}^{N_{1}-1} \Big[
   J^{1}_{1,j  }
 + J^{2}_{1,j-1}
 - B^{1}_{1,j  }
 - B^{2}_{1,j-1}
 + J^{3}_{N_{0}-1,j-1}
 + J^{4}_{N_{0}-1,j  }
 - B^{3}_{N_{0}-1,j-1}
 - B^{4}_{N_{0}-1,j  }
\Big] = \\
= \sum \limits_{i=2}^{N_{0}-1} \sum \limits_{j=1}^{N_{1}-1} \Big[
   B^{1}_{i,  j  }
 + B^{2}_{i,  j-1}
 + B^{3}_{i-1,j-1}
 + B^{4}_{i-1,j  }
\Big] .
\end{multline}
We conclude, that we have the discrete equivalent of a divergence symmetry, if the sum of all $B^{l}_{i,j}$ in the interior of $X_{d}$ vanishes,
\begin{align}
\sum \limits_{i=2}^{N_{0}-1} \sum \limits_{j=1}^{N_{1}-1} \Big[
   B^{1}_{i,  j  }
 + B^{2}_{i,  j-1}
 + B^{3}_{i-1,j-1}
 + B^{4}_{i-1,j  }
\Big] = 0 ,
\end{align}
such that only the fluxes at the boundaries of $X_{d}$ contribute to (\ref{eq:noether_discrete_divergence_conservation_law_1}).
Then, the generalised Noether charge (\ref{eq:noether_discrete_conservation_law_3}) becomes
\begin{align}\label{eq:noether_discrete_divergence_conservation_law_3}
\mathcal{J}_{d} &= \sum \limits_{\substack{j=1\\ \square^{1} = (i,j)}}^{N_{1}-1} \Big[ J_{\square^{3}} \big( \jpd^{1} \phy_{d} (\square) \big) + J_{\square^{4}} \big( \jpd^{1} \phy_{d} (\square) \big) -  B_{\square^{3}} \big( \jpd^{1} \phy_{d} (\square) \big) - B_{\square^{4}} \big( \jpd^{1} \phy_{d} (\square) \big) \Big] = \text{const. for all $i$} ,
\end{align}
in correspondence with (\ref{eq:divergence_symmetries_noether_continuous_current}).

\section{Formal Lagrangians}
\label{sec:adjoint}

In this section we review the idea of formal (or adjoint) Lagrangians \cite{AthertonHomsy:1975, Ibragimov:2006} and the corresponding Noether theorem as introduced by \citet{Ibragimov:2007a}. 
We provide a geometric view of this theory in line with the formalism adopted in the rest of the paper,
so that the application of the results of section \ref{sec:vi} to formal
Lagrangians becomes natural. At last, we propose a discrete version of
Ibragimov's extension of Noether's theorem, which allows us to prove discrete
conservation laws for variational integrators obtained from formal Lagrangians. 

\subsection{Inverse Problem of the Calculus of Variations}
\label{sec:inverse_problem}

Consider a generic system of nonlinear partial differential equations,
\begin{align}\label{eq:formal_lagrangians_system_of_pdes}
\mathcal{F} (\bx, \fu, \fu_{(1)}, ..., \fu_{(k)}) = 0 ,
\end{align}
for a field $\fu : X \rightarrow F$, where $\fu_{(k)}$ denotes all derivatives of order $k$.
Here, $\mathcal{F}$ is a (sufficiently regular) function of the point $\bx \in X$, the field $\fu (\bx)$, and its derivatives up to order $k$, taking values in $\mathbb{R}^{m}$, where $m = \dim F$, so that we have as many equations as variables. The components of $\mathcal{F}$ are denoted $\mathcal{F}_{a}$ with $1 \leq a \leq m$. 

Let $B$ be the space of functions $\fu : X \rightarrow F$ where the solution of (\ref{eq:formal_lagrangians_system_of_pdes}) is sought.
We assume that $B$ is a Banach space and denote by $B^{*}$ its topological dual, i.e., the space of linear continuous functionals from $B$ to $\mathbb{R}$.
The specific choice of the space $B$ depends on the problem at hand.
For definiteness, we require $B \hookrightarrow \big[L^{2} (X)\big]^m$, i.e., $B$ is continuously embedded in $\big[L^{2} (X)\big]^m$. If the function $\bx \mapsto \mathcal{F} (\bx, \fu, \fu_{(1)}, ..., \fu_{(k)})$ defined in (\ref{eq:formal_lagrangians_system_of_pdes}) belongs to $\big[L^{2} (X)\big]^m$, then it can be identified with the mapping $\mathcal{F} : B \rightarrow B^{*}$ defined by
\begin{align}\label{eq:formal_lagrangians_duality_pairing}
\bracket{ \mathcal{F} [\fu] , \fv } = \int \limits_{X} \fv (\bx) \cdot \mathcal{F} (\bx, \fu, \fu_{(1)}, ..., \fu_{(k)}) \, d^{n+1} \bx ,
\end{align}
where $\bracket{ \cdot , \cdot } : B^{*} \times B \rightarrow \mathbb{R}$ is the
duality pairing between elements of $B^{*}$ and $B$. With some abuse of notation
we denote by $\mathcal{F}$ both the function in
(\ref{eq:formal_lagrangians_system_of_pdes}), which is defined on a
finite-dimensional space, and the mapping $\mathcal{F} : B \rightarrow B^{*}$,
which is defined on a Banach space.

The nonlinear operator $\mathcal{F}$ admits a natural variational formulation if there exists a functional $\mathcal{A} : B \rightarrow \mathbb{R}$ of class $C^{1} (B)$, such that $\mathcal{F} = \fr \mathcal{A}$, where $\fr$ denotes the Fr\'{e}chet derivative operator.
We recall that the Fr\'{e}chet derivative of $\mathcal{A} : B \rightarrow \mathbb{R}$ is a map $\fr \mathcal{A} : B \rightarrow B^{*}$, such that for every $\fu, \fv \in B$, $\bracket{ \fr \mathcal{A}[\fu] , \fv }$ coincides with the Gateaux derivative of $\mathcal{A}$ along $\fv$ \cite{Olver:1993, Morrison:1998}, that is
\begin{align}\label{eq:formal_lagrangians_frechet_derivative}
\bracket{ \fr \mathcal{A}[\fu] , \fv }
= \lim \limits_{\eps \rightarrow 0} \dfrac{\mathcal{A} [\fu + \eps \fv] - \mathcal{A} [\fu]}{\eps}
= \dfrac{d}{d\eps} \mathcal{A} [\fu + \eps \fv] \bigg\vert_{\eps=0} .
\end{align}
In this framework, the functional derivative $\delta \mathcal{A} [\fu] / \delta \fu$, if it exists, is defined as the $L^{2}$-rep\-re\-sen\-ta\-tive of the Fr\'echet derivative.
Specifically, $\delta \mathcal{A} [\fu] / \delta \fu$ is the unique element of $\big[L^{2} (X)\big]^m$, if it exists, such that
\begin{align}
\bracket{ \fr \mathcal{A} [\fu] , \fv }
= \int \limits_{X} \dfrac{\delta \mathcal{A} [\fu]}{\delta \fu} \,\cdot\, \fv \, d\bx ,
\end{align}
for all $\fv \in B$ \cite{AbrahamMarsdenRatiu:1988}.
Here, we are implicitly assuming that boundary conditions on the field $\fu$ are
encoded in the definition of the linear space $B$; this is possible in
particular for both homogeneous Dirichlet and periodic boundary conditions.

Analogously to (\ref{eq:formal_lagrangians_frechet_derivative}), the Fr\'{e}chet
derivative at $\fu \in B$ of a functional $\mathcal{F} : B \rightarrow B^{*}$ of
class $C^1$ is the linear continuous operator $\fr \mathcal{F} [\fu] : B \rightarrow B^{*}$, such that for every $\fv \in B$, $\fr \mathcal{F} [\fu]\fv$ coincides with the Gateaux derivative of $\mathcal{F}$ along $\fv$, i.e.,
\begin{align}\label{eq:formal_lagrangians_frechet_derivative2}
\fr \mathcal{F} [\fu] \fv = \dfrac{d}{d\eps} \mathcal{F} [\fu + \eps \fv] \bigg\vert_{\eps=0} .
\end{align}
Therefore, $\fr \mathcal{F} : B \rightarrow \mathcal{L} (B, B^{*}) \cong B^{*}\otimes B^{*}$, where $\mathcal{L} (B, B^{*})$ is the space of linear continuous operators from $B$ to $B^{*}$.
If $\mathcal{F}$ admits a variational formulation, then the second-order Fr\'{e}chet derivative $\fr \mathcal{F} = \fr^{2} \mathcal{A} : B \rightarrow B^{*} \otimes B^{*}$ should be symmetric at each $\fu$, hence a necessary condition for the existence of the action $\mathcal{A}$ is
\begin{align}\label{eq:formal_lagrangians_self_adjointness}
\bracket{ \fr \mathcal{F} [\fu] \fv , \fw } &= \bracket{ \fr \mathcal{F} [\fu] \fw , \fv }
\hspace{1em}
\forall \, \fu, \fv, \fw \in B .
\end{align}
It turns out that this condition is also sufficient \cite{Olver:1993, BampiMorro:1982, Tonti:1969}.

Unfortunately, for many interesting systems this condition is not fulfilled.
Some of these systems admit a variational principle after a variable transformation.
This is the case for many equations from fluid dynamics \cite{SeligerWhitham:1968}.
But such transformations can be inconvenient since the new variables might suffer from problems with respect to non-uniqueness, boundary conditions or regularity.
In addition, the resulting variational principle might be subject to constraints on the variations which are not easily dealt with at the discrete level.
Instead, they call for extensions of the theory that complicate its application \cite{Pavlov:2011, Gawlik:2011}.
For a still larger class of systems, a variational formulation is not known, even after a change of coordinates.
Nevertheless, we can derive variational integrators for these systems by considering the following construction.

\subsection{Formal Lagrangians}
\label{sec:adjoint_formal_lagrangians}

A generic system of differential equations (\ref{eq:formal_lagrangians_system_of_pdes}) can be treated as part of a Lagrangian system, constructed by doubling the number of variables.
The action of such a Lagrangian system is
\begin{align}\label{eq:formal_lagrangians_action_fa}
\mathcal{A} [\fu, \fv] &= \bracket{ \mathcal{F} [\fu] , \fv } , 
\end{align}
with $u , v \in B$ and the pairing as defined in (\ref{eq:formal_lagrangians_duality_pairing}).
Then, the Lagrangian is
\begin{align}\label{eq:formal_lagrangians_lagrangian_fa}
L (\bx, \fu, \fu_{(1)}, ..., \fu_{(k)}, \fv, \fv_{(1)}, ..., \fv_{(k)}) = \fv (\bx) \cdot \mathcal{F} (\bx, \fu, \fu_{(1)}, ..., \fu_{(k)}) .
\end{align}
This is referred to as formal (or adjoint) Lagrangian \cite{AthertonHomsy:1975, Ibragimov:2006}.
The variational principle (\ref{eq:variational_continuous_variation_1}) applied to (\ref{eq:formal_lagrangians_lagrangian_fa})
gives the original equation,
\begin{align}\label{eq:formal_lagrangians_original_system}
\dfrac{\delta \mathcal{A} [u,v]}{\delta \fv} = \mathcal{F} [\fu] = 0 ,
\end{align}
as well as the \emph{adjoint equation}, which is defined by
\begin{align}\label{eq:formal_lagrangians_adjoint_system}
\dfrac{\delta \mathcal{A} [u,v]}{\delta \fu} = \mathcal{F}^{*} [\fu, \fv] = 0 ,
\end{align}
assuming that the functional derivative exists.

Summarising, if we consider a dynamical system described by the variables $(\fu, \fv)$, whose dynamics are governed by the equations $\mathcal{F} [\fu] = 0$ and $\mathcal{F}^{*} [\fu, \fv] = 0$, then this system has a variational formulation and contains the original system (\ref{eq:formal_lagrangians_system_of_pdes}) as a subsystem.
In the following, we refer to the system (\ref{eq:formal_lagrangians_original_system}-\ref{eq:formal_lagrangians_adjoint_system}) as \emph{extended system} describing the dynamics of the variables $(\fu, \fv)$.

\subsubsection*{Geometry of Formal Lagrangians}

As before, we identify a field $\fu$ with a section $\phy : X \rightarrow Y$ of the configuration bundle $Y$.
With $\jp^{k} \phy$ the $k$th jet prolongation of $\phy$, we can write the system of equations (\ref{eq:formal_lagrangians_system_of_pdes}) as a function on the jet bundle $\jb^{k} Y$,
\begin{align}
\mathcal{F} (\jp^{k} \phy) = 0 .
\end{align}
The configuration bundle $\tilde{Y}$ of the extended system is given by the Cartesian product
\begin{align}
\tilde{Y} &
= X \times F \times F ,
\end{align}
together with the projection
\begin{align}
\arraycolsep=2pt
\begin{array}{rcccc}
\tilde{\pi} & : & \tilde{Y} & \rightarrow & X , \\
& & (\bx, \by, \tilde{\by}) & \mapsto & \bx .
\end{array}
\end{align}
Local coordinates on $\tilde{Y}$ are denoted by $( \bx^{\mu}, \by^{a}, \tilde{\by}^{a} )$ with $1 \leq a \leq m = \dim F$.
Sections $\tilde{\phy} : X \rightarrow \tilde{Y}$ are at the same time associated to the graph of $\fu$ and $\fv$,
\begin{align}
\tilde{\phy} : \bx \mapsto \big( \bx^{\mu} , \phy^{a} (\bx) , \tilde{\phy}^{a} (\bx) \big)
= \big( \bx^{\mu} , \fu^{a} (\bx) , \fv^{a} (\bx) \big) ,
\end{align}
where the $\phy^{a}$ are associated with the physical variables $\fu^{a}$ and the $\tilde{\phy}^{a}$ are associated with the adjoint variables $\fv^{a}$.
With $\jp^{k} \tilde{\phy}$ the $k$th jet prolongation of $\tilde{\phy}$, we can write the formal Lagrangian (\ref{eq:formal_lagrangians_lagrangian_fa}) as a function on the jet bundle $\jb^{k} \tilde{Y}$,
\begin{align}\label{eq:formal_lagrangians_geometric_lagrangian}
L \big( \jp^{k} \tilde{\phy} \big) = \fv \cdot \mathcal{F} (\bx, \fu, \fu_{(1)}, ..., \fu_{(k)}) ,
\end{align}
and the formal action functional (\ref{eq:formal_lagrangians_action_fa}) as
\begin{align}\label{eq:adjoint_action}
\mathcal{A} [\tilde{\phy}] = \int \limits_{X} L \big( \jp^{k} \tilde{\phy} \big) \, d^{n+1} \bx .
\end{align}
In the following we shall consider first-order Lagrangians only, i.e., $k=1$.

\subsection{Noether Theorem}
\label{sec:formal_lagrangians_noether_theorem}

\citet{Ibragimov:2007a} has shown that the adjoint equations (\ref{eq:formal_lagrangians_adjoint_system}) inherit all symmetries of the original equations (\ref{eq:formal_lagrangians_original_system}).
Therefore it is possible to determine conservation laws of any system of differential equations, even without a natural Lagrangian, by application of the Noether theorem from section \ref{sec:vi_continous_noether_theorem} to the corresponding formal Lagrangian (\ref{eq:formal_lagrangians_geometric_lagrangian}).
Again, we restrict our attention to the case of first order field theories.

More specifically, if the the original system (\ref{eq:formal_lagrangians_system_of_pdes}) admits a symmetry related to the infinitesimal vector field
\begin{align}\label{eq:formal_lagrangians_noether_vector_field}
V &= \eta^{a} (\bx, \by) \dfrac{\partial}{\partial \by^{a}} ,
\end{align}
then the extended system (\ref{eq:formal_lagrangians_original_system}-\ref{eq:formal_lagrangians_adjoint_system}) admits a symmetry related to the vector field
\begin{align}\label{eq:formal_lagrangians_noether_extended_vector_field}
\tilde{V}
&= \eta^{a} (\bx, \by) \dfrac{\partial}{\partial \by^{a}}
 + \tilde{\eta}^{a} (\bx, \by, \tilde{\by}) \dfrac{\partial}{\partial \tilde{\by}^{a}} ,
\end{align}
with appropriately chosen coefficients $\tilde{\eta}^{a}$ \cite{Ibragimov:2007a}.
By definition, the vector field
(\ref{eq:formal_lagrangians_noether_vector_field}) describes a Lie point
symmetry of (\ref{eq:formal_lagrangians_system_of_pdes}), if there exists a
matrix-valued function $\lambda = (\lambda_{a}^{b} (\bx, \by, \bz) )$, such that
\begin{align}\label{eq:formal_lagrangians_noether_symmetry_equation}
\jp^{1} V (\mathcal{F}_{a}) = \lambda_{a}^{b} \mathcal{F}_{b} ,
\end{align}
where $\mathcal{F} (\bx, \by, \bz)$ is the function defining the system of first order equations $\mathcal{F} (\cx, \fu, \fu_{(1)})$.
We apply the jet prolongation of the extended vector field (\ref{eq:formal_lagrangians_noether_extended_vector_field}) to the Lagrangian (\ref{eq:formal_lagrangians_geometric_lagrangian}) of the extended system (\ref{eq:formal_lagrangians_original_system}-\ref{eq:formal_lagrangians_adjoint_system}), treating $\tilde{\eta}^{a}$ as unknown coefficients,
\begin{align}\label{eq:formal_lagrangians_noether_symmetry_condition}
\jp^{1} \tilde{V} (L) = \tilde{\eta}^{a} \mathcal{F}_{a} + \tilde{\by}^{a} \, \jp^{1} V (\mathcal{F}_{a}) = \tilde{\eta}^{a} \mathcal{F}_{a} + \tilde{\by}^{a} \lambda_{a}^{b} \mathcal{F}_{b} .
\end{align}
We want to choose the $\tilde{\eta}^{a}$ so that the symmetry condition (\ref{eq:noether_continuous_symmetry_3}) is satisfied.
We see that it is sufficient to set
\begin{align}\label{eq:formal_lagrangians_noether_extended_vector_field_component}
\tilde{\eta}^{a} = - \lambda^{a}_{b} \, \tilde{\by}^{b}
\end{align}
for (\ref{eq:formal_lagrangians_noether_symmetry_condition}) to vanish.
The extended vector field (\ref{eq:formal_lagrangians_noether_extended_vector_field}) therefore becomes
\begin{align}\label{eq:formal_lagrangians_noether_extended_vector_field_2}
\tilde{V} &= \eta^{a} \, \dfrac{\partial}{\partial \by^{a}} - \lambda^{a}_{b} \, \tilde{\by}^{b} \, \dfrac{\partial}{\partial \tilde{\by}^{a}} ,
\end{align}
with $\lambda$ given in the symmetry condition (\ref{eq:formal_lagrangians_noether_symmetry_equation}).
In the line of the standard Noether theorem, c.f. section \ref{sec:vi_continous_noether_theorem}, if (\ref{eq:formal_lagrangians_noether_extended_vector_field}) is a symmetry of the formal Lagrangian (\ref{eq:formal_lagrangians_geometric_lagrangian}), we obtain the local conservation law
\begin{align}\label{eq:formal_lagrangians_noether_theorem}
\div \tilde{J} (\jp^{1} \tilde{\phy}) = 0 ,
\end{align}
for all solutions $\tilde{\phy}$ of the Euler-Lagrange equations on $\tilde{Y}$, with the Noether current $\tilde{J}$ given by
\begin{align}\label{eq:formal_lagrangians_noether_current}
\tilde{J}^{\mu} \big( \jp^{1} \tilde{\phy} \big)
= \dfrac{\partial L}{\partial \bz_{\mu}^{a}} \big( \jp^{1} \tilde{\phy} \big) \cdot \eta^{a} (\tilde{\phy}) + \dfrac{\partial L}{\partial \tilde{\bz}_{\mu}^{a}} \big( \jp^{1} \tilde{\phy} \big) \cdot \tilde{\eta}^{a} (\tilde{\phy}) .
\end{align}
For formal Lagrangians of the form (\ref{eq:formal_lagrangians_geometric_lagrangian}), the second terms on the right-hand side are zero. Sometimes, however, we might want to symmetrise the Lagrangian (see e.g. section \ref{sec:vorticity}), in which case these terms do contribute to the Noether current.
From equation (\ref{eq:formal_lagrangians_noether_current}), we obtain the (global) Noether charge (\ref{eq:noether_continuous_charge}),
\begin{align}\label{eq:formal_lagrangians_noether_charge_full}
\tilde{\mathcal{J}} = \int \bigg[ \dfrac{\partial L}{\partial \bz_{t}^{a}} \big( \jp^{1} \tilde{\phy} \big) \cdot \eta^{a} (\tilde{\phy}) + \dfrac{\partial L}{\partial \tilde{\bz}_{t}^{a}} \big( \jp^{1} \tilde{\phy} \big) \cdot \tilde{\eta}^{a} (\tilde{\phy}) \bigg] \, d^{n} \bx ,
\end{align}
which is constant for all times $\ct$.
Usually, formal Lagrangians of the form (\ref{eq:formal_lagrangians_geometric_lagrangian}) do not feature a time derivative of the adjoint variables, so that the Noether charge simplifies accordingly.

\subsubsection*{Restriction of Conservation Laws}

The conservation law thus obtained is of course a conservation law of the extended system (\ref{eq:formal_lagrangians_original_system}-\ref{eq:formal_lagrangians_adjoint_system}), and therefore generally depends on both $\fu$ and $\fv$.
In order to restrict conservation laws of the extended systems to conservation laws of the physical system, we have to find a suitable way of restricting solutions of the extended system to solutions of the physical system.
For that purpose, the following two definitions will become useful.

If the adjoint equation (\ref{eq:formal_lagrangians_adjoint_system}), restricted to $\fv = \fu$, becomes equivalent to the original equation (\ref{eq:formal_lagrangians_system_of_pdes}), i.e.,
\begin{align}\label{eq:adjoint_system_self_adjoint}
\mathcal{F}^{*} [\fu, \fu] = \lambda \mathcal{F} [\fu],
\end{align}
for some matrix $\lambda$, possibly depending on the fields and their derivatives, the system (\ref{eq:formal_lagrangians_system_of_pdes}) is called self-adjoint in the sense of Ibragimov \cite{Ibragimov:2006}.
If $\mathcal{F} [\fu]$ is a linear operator and $\lambda$ the identity matrix $\mathbb{1}_{m \times m}$, the above definition coincides with the standard definition of formally self-adjoint operators.
In general, however, the systems we are considering are not self-adjoint.
\citeauthor{Ibragimov:2007b} relaxed the requirement of self-adjointness by introducing the concept of quasi-self-adjointness \cite{Ibragimov:2007b, Ibragimov:2010}.
This is a generalisation of self-adjointness where $v=u$ is generalised to $v = \phi(u)$ for some function $\phi : F \rightarrow F$.
The advantage of self-adjointness or quasi-self-adjointness in the sense of Ibragimov is the possibility to build a solution of the full extended system (\ref{eq:formal_lagrangians_original_system}-\ref{eq:formal_lagrangians_adjoint_system}) given a solution of the original problem (\ref{eq:formal_lagrangians_original_system}).

Given a diffeomorphism $\phi : F \rightarrow F$ of $F$ into itself, we can build the embedding
\begin{align}\label{eq:formal_lagrangians_embedding_geometric}
\arraycolsep=2pt
\begin{array}{rcccl}
\Phi & : & Y & \hookrightarrow & \quad \tilde{Y} , \\
& & \big( \bx, \by \big) & \mapsto & \big( \bx, \by, \phi (\by) \big) .
\end{array}
\end{align}
In order to use this embedding to restrict the Noether current to the physical
system, we need to lift $\Phi$ to a map from $\jb^{1} Y$ to $\jb^{1} \tilde{Y}$.
With this aim, let us consider a generic section $\phy$ of $Y$ for which the composed map
\begin{align}
\tilde{\phy} = \Phi \circ \phy : X \rightarrow \tilde{Y}
\end{align}
amounts to
\begin{align}
\tilde{\phy} : \bx \mapsto \big( \bx , \fu (\bx) , \phi ( \fu (\bx) ) \big) = \big( \bx , \fu (\bx) , \fv (\bx) \big) ,
\end{align}
where $\fu(\bx)$ is the field component corresponding to $\phy$, and we have defined the second field $\fv = \phi \circ \fu : X \rightarrow F$.
It follows that the condition $\tilde{\phy} \circ \tilde{\pi} = \id_{X}$ is satisfied and $\tilde{\phy} = \Phi \circ \phy$ is indeed a section of $\tilde{Y}$, i.e., the composition with $\Phi$ maps sections of $Y$ into sections of $\tilde{Y}$.
By the chain rule, we compute the first jet of $\tilde{\phy}$, that is
\begin{align}
\jp^{1} \tilde{\phy} (\bx)
\nonumber
&= \big( \bx , \, \fu (\bx) , \, \fv (\bx)          , \, D \fu (\bx) , \, D \fv (\bx) \big) \\
&= \big( \bx , \, \fu (\bx) , \, \phi ( \fu (\bx) ) , \, D \fu (\bx) , \, D \phi ( \fu (\bx) ) \cdot D \fu (\bx) \big) ,
\end{align}
where $Df$ denotes the Jacobian matrix of a function $f$.
We now can define the lift
\begin{align}\label{eq:formal_lagrangians_embedding_lifted}
\jp^{1} \Phi : \jb^{1} Y \rightarrow \jb^{1} \tilde{Y}
\end{align}
by
\begin{align}\label{eq:formal_lagrangians_embedding_lifted_2}
\jp^{1} \Phi  \left( \bx, \by, \bz \right) = \big( \bx , \, \by , \, \phi ( \by ) , \, \bz , \, D \phi (\by) \cdot \bz \big) ,
\end{align}
which by construction satisfies the identity
\begin{align}
\jp^{1} \big( \Phi \circ \phy \big) = \jp^{1} \Phi \circ \jp^{1} \phy ,
\end{align}
as desired.
With an abuse of notation, the symbol $\jp^{1} \Phi$ is not used here in the usual sense of jet prolongation, instead $\jp^{1} \Phi$ is a lift of $\Phi$ up to the first jet bundle.
We can now use (\ref{eq:formal_lagrangians_embedding_lifted_2}) to pull back the Noether current in (\ref{eq:formal_lagrangians_noether_theorem}) if we assume that $\tilde{\phy}$ can be realised in the form $\tilde{\phy} = \Phi \circ \phy$,
\begin{align}
\tilde{J} \big( \jp^{1} \tilde{\phy} \big)
&= \tilde{J} \big( \jp^{1} ( \Phi \circ \phy ) \big)
 = \tilde{J} \big( \jp^{1} \Phi \circ \jp^{1} \phy \big)
 = \big( \tilde{J} \circ \jp^{1} \Phi \big) \big( \jp^{1} \phy \big)
.
\end{align}
If $\tilde{\phy} = \Phi \circ \phy$ solves the equations of the extended system (\ref{eq:formal_lagrangians_original_system}-\ref{eq:formal_lagrangians_adjoint_system}), then $\div \tilde{J} \big( \jp^{1} \tilde{\phy} \big) = 0$.
Upon defining the restricted Noether current by
\begin{align}\label{eq:formal_lagrangians_noether_current_restricted}
J = \tilde{J} \circ \jp^{1} \Phi ,
\end{align}
the Noether theorem (\ref{eq:formal_lagrangians_noether_theorem}) takes the form
\begin{align}\label{eq:formal_lagrangians_noether_theorem_restricted}
\div J (\jp^{1} \phy) = 0 ,
\end{align}
which expresses the local conservation law for the physical system.

This result is crucial for the application of variational integrators to formal Lagrangians. 
Without the construction of a solution $(\fu, \phi(\fu))$ of the extended system from a solution $\fu$ of the physical system, it is in general not possible to determine the discrete momenta that are conserved by the variational integrator due to symmetries of the physical system.

\subsubsection*{Generalisations}

To simplify the derivations and the analysis of the conservation laws, it is sometimes useful to add a term $G (\jp^{k} \phy)$ to the Lagrangian, that is a function of the coordinates, the physical variables and their derivatives, but not of the adjoint variables and their derivatives.
The modified Lagrangian,
\begin{align}\label{eq:formal_lagrangians_generalisation}
L' \big( \jp^{k} \tilde{\phy} \big) &= L \big( \jp^{k} \tilde{\phy} \big) + G \big( \jp^{k} \tilde{\phy} \big) &
& \text{with} &
\dfrac{\partial G}{\partial \tilde{\by}^{a}} \big( \jp^{k} \tilde{\phy} \big) &= 0 &
& \text{and} &
\dfrac{\partial G}{\partial \tilde{\bz}_{\mu}^{a}} \big( \jp^{k} \tilde{\phy} \big) &= 0 &
& \text{for all} \; \mu, a ,
\end{align}
yields the same physical equations of motion as $L$, but in general will lead to different adjoint equations.
This freedom can be used to simplify the search for the embedding $\Phi$ used to restrict the conservation laws of the extended system to the physical system.

For example, for equations of advection type,
\begin{align}
\fu_{\ct} + A(\fu, \fu_{\cx}) = 0 ,
\end{align}
where $A$ is some possibly nonlinear function of the field $\fu$ and its spatial derivatives,
it is possible to construct an extended system which is guaranteed to be always self-adjoint.
It suffices to choose $G=-u \cdot A$, so that the Lagrangian becomes
\begin{align}
L' \big( \jp^{k} \tilde{\phy} \big) = \fv \cdot \fu_{\ct} + \fv \cdot A(\fu, \fu_{\cx}) - \fu \cdot A(\fu, \fu_{\cx}) .
\end{align}
While this simplifies the construction of the embedding $\Phi$, it complicates the extension of symmetries. In general it will not be possible to determine the components of the generating vector field in an algorithmic way as in (\ref{eq:formal_lagrangians_noether_extended_vector_field_component}).

\subsubsection*{Discrete Embedding}

At this point, the discretisation of formal Lagrangians can be carried out
straightforwardly as an application of the theory reviewed in section
\ref{sec:vi}. The only missing element is the restriction of discrete
conservation laws of the extended system to the physical system.  

Since the fibres of $Y_{d}$ are copies of $F$, we can consider a diffeomorphism $\phi : F \rightarrow F$, as in the continuous case and define the discrete embedding
\begin{align}\label{eq:formal_lagrangians_embedding_discrete}
\arraycolsep=2pt
\begin{array}{rcccl}
\Phi_{d} & : & Y_{d} & \hookrightarrow & \quad \tilde{Y}_{d} , \\
& & \big( \by_{i,j} \big) & \mapsto & \big( \by_{i,j}, \, \phi (\by_{i,j}) \big) ,
\end{array}
\end{align}
and its lift to the discrete first jet bundle,
\begin{align}
\jp^{1} \Phi_{d} : \jb^{1} Y_{d} \rightarrow \jb^{1} \tilde{Y}_{d},
\end{align}
which is given by
\begin{align}
\jpd^{1} \Phi_{d} (\square) = \big( \square, \; \by_{\square^{1}}, \by_{\square^{2}}, \by_{\square^{3}}, \by_{\square^{4}}, \; \phi( \by_{\square^{1}} ), \phi( \by_{\square^{2}} ), \phi( \by_{\square^{3}} ), \phi( \by_{\square^{4}} ) \big) .
\end{align}
With this we can pullback discrete conservation laws (\ref{eq:noether_discrete_conservation_law_3}) of the extended system,
\begin{align}\label{eq:formal_lagrangians_discrete_noether_charge_extended}
\tilde{\mathcal{J}}_{d} &= \sum \limits_{\substack{j=1\\ \square^{1} = (i,j)}}^{N_{1}-1} \Big[ \tilde{J}_{\square^{3}} \big( \jpd^{1} \tilde{\phy}_{d} (\square) \big) + \tilde{J}_{\square^{4}} \big( \jpd^{1} \tilde{\phy}_{d} (\square) \big) \Big] = \text{const. for all $i$} ,
\end{align}
to the physical system by defining the restricted discrete momentum map
\begin{align}
J_{\square^{l}} = \tilde{J}_{\square^{l}} \circ \jpd^{1} \Phi_{d} (\square) ,
\end{align}
such that the restricted Noether charge becomes
\begin{align}\label{eq:formal_lagrangians_discrete_noether_charge_restricted}
\mathcal{J}_{d} &= \sum \limits_{\substack{j=1\\ \square^{1} = (i,j)}}^{N_{1}-1} \Big[ J_{\square^{3}} \big( \jpd^{1} \phy_{d} (\square) \big) +  J_{\square^{4}} \big( \jpd^{1} \phy_{d} (\square) \big) \Big] = \text{const. for all $i$} .
\end{align}
It is worth noticing that we do not pull back the Noether charge $\tilde{\mathcal{J}}$ but the momentum maps $\tilde{J}_{\square^{l}}$, just as in the continuous case  (\ref{eq:formal_lagrangians_noether_current_restricted}), where we pull back the Noether current $\tilde{J}$.

\section{Applications}
\label{sec:applications}

In this section we present two applications, the linear advection equation and the vorticity equation.
The linear advection equation in one spatial dimension is a prototypical example that shares many characteristics with more complicated systems from plasma physics and fluid dynamics.
The vorticity equation in two spatial dimensions is an important equation of fluid dynamics which is widely used for example in atmospheric sciences but also in plasma physics. It consists of an advection equation which is coupled to a Poisson equation.

On the continuous level, we construct the formal Lagrangian and compute variations of the corresponding action functional. We check the resulting equations of motion for self-adjointness or quasi-self-adjointness in order to define the restriction map $\phi$. Then we check the equations for symmetries and compute the extended generating vector field. This vector field is applied to the Lagrangian in order to compute the corresponding conservation laws.
After discretising the formal Lagrangian, we derive the actual variational integrator and check the discrete equations of motion for discrete self-adjointness or quasi-self-adjointness. Then we check the discrete Lagrangian for symmetries and compute the discrete conservation laws.
We will discuss different discretisations of the Lagrangian, especially with respect to different quadrature rules, and their consequences for properties of the resulting variational integrator (e.g., explicit vs. implicit) and the corresponding conservation laws.

\subsubsection*{Discretisation of Formal Lagrangians}

Before we proceed, we want to give a short comment on the discretisation of formal Lagrangians.
A straight forward application of the Veselov discretisation, i.e., using the midpoint quadrature rule for both time and space integration, as it was presented in section \ref{sec:vi_discrete_action_principle} and as it is used in many works on variational integrators, does not appear to be suitable for formal Lagrangians due to their particular structure.
As we will see in section \ref{sec:advection} on the linear advection equation, it results in numerical schemes, which are prone to unphysical oscillations in space. 
Instead, when discretising the derivatives with first-order finite differences, the spatial integral can be approximated by the trapezoidal or Simpson rule. Another option to avoid this issue is to move to a slightly more abstract Galerkin framework and use e.g. Lagrange polynomials \cite{Chen:2008} or splines \cite{Kraus:2015:Splines} as basis functions.

\subsection{Some Definitions}

So far, we have only considered the midpoint rule (Veselov discretisation) for approximating the Lagrangian $L(\fu, \fu_{t}, \fu_{x})$, namely,
\begin{align}\label{eq:applications_discrete_lagrangian_midpoint}
L_{d}^{\mathrm{mp}} \big( \jpd^{1} \phy_{d} (\square) \big) = h_{t} h_{x} \; L \bigg(
\nonumber
& \tfrac{1}{4} \Big( \phy_{\square^{1}} + \phy_{\square^{2}} + \phy_{\square^{3}} + \phy_{\square^{4}} \Big), \\
\nonumber
& \tfrac{1}{2} \Big( \tfrac{\phy_{\square^{4}} - \phy_{\square^{1}}}{h_{t}} + \tfrac{\phy_{\square^{3}} - \phy_{\square^{2}}}{h_{t}} \Big), \\
& \tfrac{1}{2} \Big( \tfrac{\phy_{\square^{2}} - \phy_{\square^{1}}}{h_{x}} + \tfrac{\phy_{\square^{3}} - \phy_{\square^{4}}}{h_{x}} \Big)
\bigg) .
\end{align}
As will be discussed below, for formal Lagrangians the midpoint rule does not appear to lead to stable variational integrators. Therefore we will also consider some alternative discretisations of the Lagrangian.
First, the trapezoidal rule,
\begin{align}\label{eq:applications_discrete_lagrangian_leapfrog}
L_{d}^{\mathrm{tr}} \big( \jpd^{1} \phy_{d} (\square) \big) = \dfrac{h_{t} h_{x}}{4} \, \bigg(
\nonumber
  & L \Big( \phy_{\square^{1}}, \tfrac{\phy_{\square^{4}} - \phy_{\square^{1}}}{h_{t}}, \tfrac{\phy_{\square^{2}} - \phy_{\square^{1}}}{h_{x}} \Big) \\
\nonumber
+ & L \Big( \phy_{\square^{2}}, \tfrac{\phy_{\square^{3}} - \phy_{\square^{2}}}{h_{t}}, \tfrac{\phy_{\square^{2}} - \phy_{\square^{1}}}{h_{x}} \Big) \\
\nonumber
+ & L \Big( \phy_{\square^{3}}, \tfrac{\phy_{\square^{3}} - \phy_{\square^{2}}}{h_{t}}, \tfrac{\phy_{\square^{3}} - \phy_{\square^{4}}}{h_{x}} \Big) \\
+ & L \Big( \phy_{\square^{4}}, \tfrac{\phy_{\square^{4}} - \phy_{\square^{1}}}{h_{t}}, \tfrac{\phy_{\square^{3}} - \phy_{\square^{4}}}{h_{x}} \Big)
\bigg) ,
\end{align}
which leads to the (explicit) Leapfrog scheme.
Second, a combination of both, the midpoint rule for time integration and the trapezoidal rule for spatial integration,
\begin{align}\label{eq:applications_discrete_lagrangian_simplified}
L_{d}^{\mathrm{mt}} \big( \jpd^{1} \phy_{d} (\square) \big)
\nonumber
&= \dfrac{h_{t} h_{x}}{2} L \bigg(
\tfrac{\phy_{\square^{1}} + \phy_{\square^{4}}}{2},
\tfrac{\phy_{\square^{4}} - \phy_{\square^{1}}}{h_{t}},
\tfrac{1}{2} \Big(
   \tfrac{\phy_{\square^{2}} - \phy_{\square^{1}}}{h_{x}}
 + \tfrac{\phy_{\square^{3}} - \phy_{\square^{4}}}{h_{x}}
\Big)
\bigg)  \\
&+ \dfrac{h_{t} h_{x}}{2} L \bigg(
\tfrac{\phy_{\square^{2}} + \phy_{\square^{3}}}{2},
\tfrac{\phy_{\square^{3}} - \phy_{\square^{2}}}{h_{t}},
\tfrac{1}{2} \Big(
   \tfrac{\phy_{\square^{2}} - \phy_{\square^{1}}}{h_{x}}
 + \tfrac{\phy_{\square^{3}} - \phy_{\square^{4}}}{h_{x}}
\Big)
\bigg) .
\end{align}
This leads to a simplified implicit scheme with the same conservation properties as the midpoint scheme derived from (\ref{eq:applications_discrete_lagrangian_midpoint}) but without the aforementioned instability.

We want to define some shorthand notation similar to (\ref{eq:vi_discrete_average}-\ref{eq:vi_discrete_derivative_1}),
\begin{align}
\label{eq:applications_definitions_discrete_average}
u &\rightarrow \dfrac{1}{4} \Big( u_{\square^1} + u_{\square^2} + u_{\square^3} + u_{\square^4} \Big) \equiv \overline{u} (\square) \\
\label{eq:applications_definitions_discrete_derivative_t}
\dfrac{\partial u}{\partial t} &\rightarrow \dfrac{1}{2} \bigg( \dfrac{u_{\square^4} - u_{\square^1}}{h_{t}} + \dfrac{u_{\square^3} - u_{\square^2}}{h_{t}} \bigg) \equiv \overline{u}_{t} (\square) , \\
\label{eq:applications_definitions_discrete_derivative_x}
\dfrac{\partial u}{\partial x} &\rightarrow \dfrac{1}{2} \bigg( \dfrac{u_{\square^2} - u_{\square^1}}{h_{x}} + \dfrac{u_{\square^3} - u_{\square^4}}{h_{x}} \bigg) \equiv \overline{u}_{x} (\square) .
\end{align}
For the product of two fields as well as the product of a field and a derivative, according to the mixed use of trapezoidal and midpoint rule, we obtain
\begingroup
\allowdisplaybreaks
\begin{align}
\label{eq:applications_definitions_discrete_average_trapezoidal}
v u &\rightarrow
\dfrac{1}{2} \bigg[
\dfrac{v_{\square^{1}} + v_{\square^{4}}}{2}
\dfrac{u_{\square^{1}} + u_{\square^{4}}}{2} +
\dfrac{v_{\square^{2}} + v_{\square^{3}}}{2}
\dfrac{u_{\square^{2}} + u_{\square^{3}}}{2}
\bigg] \equiv \overline{v u} (\square) ,
\\
\label{eq:applications_definitions_discrete_derivative_trapezoidal_t}
v \, \dfrac{\partial u}{\partial t} &\rightarrow
\dfrac{1}{2} \bigg[
\dfrac{v_{\square^{1}} + v_{\square^{4}}}{2}
\dfrac{u_{\square^{4}} - u_{\square^{1}}}{h_{t}} +
\dfrac{v_{\square^{2}} + v_{\square^{3}}}{2}
\dfrac{u_{\square^{3}} - u_{\square^{2}}}{h_{t}}
\bigg] \equiv \overline{v u}_{t} (\square) ,
\\
\label{eq:applications_definitions_discrete_derivative_trapezoidal_x}
v \, \dfrac{\partial u}{\partial x} &\rightarrow
\dfrac{1}{2} \bigg[
\dfrac{v_{\square^{1}} + v_{\square^{2}}}{2}
\dfrac{u_{\square^{2}} - u_{\square^{1}}}{h_{x}} +
\dfrac{v_{\square^{3}} + v_{\square^{4}}}{2}
\dfrac{u_{\square^{3}} - u_{\square^{4}}}{h_{x}}
\bigg] \equiv \overline{v u}_{x} (\square) 
.
\end{align}
\endgroup
This will prove handy in writing the discrete Lagrangians.

\subsection{Linear Advection Equation}
\label{sec:advection}

The initial value problem for the linear advection equation is an instructive example for the derivation of numerical schemes.
Its analytic solution is known and its structure is similar to the systems we will consider in subsequent publications \cite{KrausMajSonnendruecker:2015, KrausMaj:2015, KrausTassi:2015}.
It is given by
\begin{align}\label{eq:linear_advection_equation}
\fu_{t} + c \fu_{\cx} &= 0 ,
\end{align}
with initial condition
\begin{align}
\fu (0, \cx) &= \fu_{0} (\cx) ,
\end{align}
and describes the advection of a field $\fu(\ct, \cx)$ with a constant velocity $c$ that is given as a parameter.
The analytic solution is
\begin{align}
\fu (\ct, \cx) = \fu_{0} (\cx - c \ct) .
\end{align}
Here, $\cx$ denotes the spatial variable only, while the base space $X$ comprises both space and time, i.e., $(\ct, \cx) \in X$.
The field $u$ can be interpreted as a density, such that $\int c \fu \, dx$ is the total momentum and $\tfrac{1}{2} \int c^{2} \fu \, d\cx$ is the total energy of the system. Both are proportional to the ``mass'' $\int \fu \, d\cx$ which is a conserved quantity.
In this case, we have
\begin{align}\label{eq:linear_advection_equation_operator}
\mathcal{F} [\fu] = \fu_{\ct} + c \fu_{\cx} ,
\end{align}
and the Fr\'{e}chet derivative (\ref{eq:formal_lagrangians_frechet_derivative2}),
\begin{align}
D \mathcal{F} [\fu] \fv = \fv_{\ct} + c \fv_{\cx} ,
\end{align}
is not self-adjoint, as under the assumption of appropriate boundary conditions (e.g., periodic or homogeneous Dirichlet), integration by parts gives
\begin{align}
\bracket{D \mathcal{F} [\fu] \fv, \fw} = - \bracket{D \mathcal{F} [\fu] \fw, \fv} .
\end{align}
It follows that there is no functional $\mathcal{A} [u]$ such that $\fr \mathcal{A} = \mathcal{F}$, i.e., $\mathcal{F}$ cannot be written as a variational problem.

\subsubsection{Formal Lagrangian}
\label{sec:advection_formal_lagrangian}

The formal Lagrangian for the advection equation is obtained by multiplying (\ref{eq:linear_advection_equation}) with the auxiliary variable $\fv (\ct, \cx)$.
The solution vector of the extended system is denoted $(\fu,\fv)$ with the corresponding section of the configuration bundle $\tilde{Y}$ denoted by $\tilde{\phy}$. In coordinates, $\tilde{\phy}$ is given explicitly as
\begin{align}
\tilde{\phy} : ( \ct, \cx ) \mapsto ( \ct , \cx , \fu, \fv ) ,
\end{align}
and its jet prolongation as
\begin{align}
\jp^{1} \tilde{\phy} : ( \ct, \cx ) \mapsto ( \ct , \cx , \fu, \fv, \fu_{\ct}, \fu_{\cx}, \fv_{\ct}, \fv_{\cx} ) .
\end{align}
With that, the Lagrangian can be written as
\begin{align}\label{eq:linear_advection_lagrangian}
L = \tilde{\by} \, \big( \bz_{\ct} + c \bz_{\cx} \big) ,
\end{align}
with action functional
\begin{align}\label{eq:linear_advection_action}
\mathcal{A} [\tilde{\phy}] = \int L ( \jp^{1} \tilde{\phy}) \, d\ct \, d\cx = \int \fv \big( \fu_{\ct} + c \fu_{\cx} \big) \, d\ct \, d\cx .
\end{align}
The variation of the action with respect to the adjoint variable $\fv$ gives the advection equation,
\begin{align}
\dfrac{\delta \mathcal{A} [\tilde{\phy}]}{\delta \fv} = \fu_{\ct} + c \fu_{\cx} = 0 ,
\end{align}
while the variation with respect to the original variable $\fu$ yields the adjoint equation,
\begin{align}
\dfrac{\delta \mathcal{A} [\tilde{\phy}]}{\delta \fu} = - \fv_{t} - c \fv_{\cx} = 0 ,
\end{align}
always assuming vanishing variations at the boundaries.
It is immediately observed that the adjoint equation has the same solution as the original equation, such that if $\fu$ is a solution of the advection equation, then $(\fu, \fu)$ solves the Euler-Lagrange equations of the extended Lagrangian (\ref{eq:linear_advection_lagrangian}).
This means, that the advection equation is actually self-adjoint in the sense of Ibragimov, cf. equation (\ref{eq:adjoint_system_self_adjoint}), with $\lambda = -1$.

\subsubsection{Continuous Conservation Laws}
\label{sec:advection_noether_continuous}

We will consider the conservation of mass and the $L^{2}$ norm of $\fu$.
The general form of the Noether charge (\ref{eq:formal_lagrangians_noether_charge_full}) is
\begin{align}\label{eq:linear_advection_conservation_law_global}
\tilde{\mathcal{J}}
= \int \bigg[ \dfrac{\partial L}{\partial \bz_{t}^{a}} \big( \jp^{1} \tilde{\phy} \big) \cdot \eta^{a} (\tilde{\phy}) + \dfrac{\partial L}{\partial \tilde{\bz}_{t}^{a}} \big( \jp^{1} \tilde{\phy} \big) \cdot \tilde{\eta}^{a} (\tilde{\phy}) \bigg] \, d^{n} \bx 
= \int \tilde{J}^{t} (\jp^{1} \tilde{\phy}) \, d^{n} \bx 
= \int \fv \eta \, d^{n} \bx 
,
\end{align}
with $\tilde{\mathcal{J}} = \text{const.}$ for all $\ct$.
Since the advection equation is self-adjoint, we can use the identity as embedding map, i.e., $\phi (\fu) = \fu$, such that $\fv = \fu$, and the restriction of conservation laws for the extended system to the original system is straight forward.

The infinitesimal generator of the transformation that leads to mass conservation and its first jet prolongation are
\begin{align}\label{eq:linear_advection_generator_U_L1}
V = \dfrac{\partial}{\partial \by} ,
\hspace{3em}
\jp^{1} V = \dfrac{\partial}{\partial \by} .
\end{align}
Applying $\jp^{1} V$ to (\ref{eq:linear_advection_equation_operator}), we obtain
\begin{align}
\jp^{1} V (\mathcal{F}) = 0 ,
\end{align}
such that by (\ref{eq:formal_lagrangians_noether_symmetry_equation}) and (\ref{eq:formal_lagrangians_noether_extended_vector_field_component}) we find $\tilde{\eta} = 0$ and the extended vector field (\ref{eq:formal_lagrangians_noether_extended_vector_field}) is
\begin{align}\label{eq:linear_advection_generator_L1}
\tilde{V} &= \eta \dfrac{\partial}{\partial \by} + \tilde{\eta} \dfrac{\partial}{\partial \tilde{\by}} = \dfrac{\partial}{\partial \by} ,
\end{align}
and its prolongation is
\begin{align}
\jp^{1} \tilde{V} &= \dfrac{\partial}{\partial \by} .
\end{align}
The invariance of the Lagrangian (\ref{eq:linear_advection_lagrangian}) is easily confirmed, as
\begin{align}
\jp^{1} \tilde{V} (L) = \dfrac{\partial L}{\partial \by} = 0 .
\end{align}
The resulting conservation law (\ref{eq:linear_advection_conservation_law_global}) reads
\begin{align}
\tilde{\mathcal{J}}
= \int \tilde{J}^{t} (\jp^{1} \tilde{\phy}) \, d\cx
= \int \fv \, d\cx 
= \text{const. for all $\ct$} .
\end{align}
Upon restricting the Noether current $\tilde{J}$ with $\fv = \phi (\fu) = \fu$, this becomes conservation of the total mass in the system,
\begin{align}
\mathcal{J}
= \int \big( \tilde{J}^{t} \circ \jp^{1} \Phi \big) (\jp^{1} \phy) \, d\cx
= \int \fu \, d\cx
= \text{const. for all $\ct$} .
\end{align}
The conservation of momentum and energy follows exactly in the same way for $\eta = c$ and $\eta = \tfrac{1}{2} c^{2}$, respectively.

The infinitesimal generator corresponding to conservation of the $L^{2}$ norm is given by
\begin{align}\label{eq:linear_advection_generator_U_L2}
V &= \by \dfrac{\partial}{\partial \by} .
\end{align}
Applying its prolongation,
\begin{align}
\jp^{1} V = \by \dfrac{\partial}{\partial \by} + \bz_{\ct} \dfrac{\partial}{\partial \bz_{\ct}} + \bz_{\cx} \dfrac{\partial}{\partial \bz_{\cx}} ,
\end{align}
to (\ref{eq:linear_advection_equation}), we obtain
\begin{align}
\jp^{1} V (\mathcal{F}) = \bz_{\ct} + c \bz_{\cx} = \lambda \mathcal{F} ,
\end{align}
with $\lambda = 1$, such that by (\ref{eq:formal_lagrangians_noether_symmetry_equation}) and (\ref{eq:formal_lagrangians_noether_extended_vector_field_component}) the extended vector field (\ref{eq:formal_lagrangians_noether_extended_vector_field}) becomes
\begin{align}\label{eq:linear_advection_generator_L2}
\tilde{V} &= \by \, \dfrac{\partial}{\partial \by} - \tilde{\by} \, \dfrac{\partial}{\partial \tilde{\by}} ,
\end{align}
and its prolongation is
\begin{align}
\jp^{1} \tilde{V}
&= \by \, \dfrac{\partial}{\partial \by}
 + \bz_{\ct} \, \dfrac{\partial}{\partial \by_{\ct}}
 + \bz_{\cx} \, \dfrac{\partial}{\partial \by_{\cx}}
 - \tilde{\by} \, \dfrac{\partial}{\partial \tilde{\by}}
 - \tilde{\bz}_{\ct} \, \dfrac{\partial}{\partial \tilde{\by}_{\ct}}
 - \tilde{\bz}_{\cx} \, \dfrac{\partial}{\partial \tilde{\by}_{\cx}} .
\end{align}
The invariance condition of the Lagrangian (\ref{eq:linear_advection_lagrangian}),
\begin{align}
\jp^{1} \tilde{V} (L) (\jp^{1} \tilde{\phy})
&= - \fv \, \big( \fu_{\ct} + c \fu_{\cx} \big) + \fv \, \big( \fu_{\ct} + c \fu_{\cx} \big) = 0 ,
\end{align}
is again fulfilled.
The corresponding Noether charge (\ref{eq:linear_advection_conservation_law_global}) is
\begin{align}
\tilde{\mathcal{J}}
= \int \tilde{J}^{t} (\jp^{1} \tilde{\phy}) \, d\cx
= \int \fv \fu \, d\cx .
\end{align}
Upon restricting the Noether current $\tilde{J}$ with $\fv = \phi (\fu) = \fu$, this becomes conservation of the $L^{2}$ norm,
\begin{align}
\mathcal{J}
= \int \big( \tilde{J}^{t} \circ \jp^{1} \Phi \big) (\jp^{1} \phy) \, d\cx
= \int \fu^{2} \, d\cx
= \text{const. for all $\ct$} .
\end{align}
In the next sections we will derive the same conservation laws on the discrete level.

\subsubsection{Midpoint Discretisation}
\label{sec:linear_advection_midpoint}

We discretise the Lagrangian (\ref{eq:linear_advection_lagrangian}) as described in section \ref{sec:vi_discrete_action_principle} according to the midpoint rule (\ref{eq:applications_discrete_lagrangian_midpoint}), obtaining
\begin{multline}\label{eq:linear_advection_discrete_lagrangian_midpoint}
L_{d}^{\mathrm{mp}} \big( \jp^{1} \tilde{\phy} (\square) \big) = h_{t} h_{x} \, \dfrac{1}{4} \Big( v_{\square^{1}} + v_{\square^{2}} + v_{\square^{3}} + v_{\square^{4}} \Big) \times \\
\times \bigg[ \dfrac{1}{2} \bigg( \dfrac{u_{\square^{4}} - u_{\square^{1}}}{h_{t}} + \dfrac{u_{\square^{3}} - u_{\square^{2}}}{h_{t}} \bigg) + \dfrac{c}{2} \bigg( \dfrac{u_{\square^{2}} - u_{\square^{1}}}{h_{x}} + \dfrac{u_{\square^{3}} - u_{\square^{4}}}{h_{x}} \bigg) \bigg] .
\end{multline}
In the shorthand notation defined in (\ref{eq:applications_definitions_discrete_average}-\ref{eq:applications_definitions_discrete_derivative_x}), this becomes
\begin{align}\label{eq:linear_advection_discrete_lagrangian_short}
L_{d} \big( \jp^{1} \tilde{\phy} (\square) \big)
= h_{t} h_{x} \, \overline{v} (\square) \, \big[ \overline{u}_{t} (\square) + c \, \overline{u}_{x} (\square) \big] ,
\end{align}
which resembles the continuous Lagrangian.

\subsubsection*{Variational Integrator}

The discrete Euler-Lagrange field equations (\ref{eq:vi_discrete_euler_lagrange_equations}), corresponding to the variation of $v_{d}$, are computed as
\begin{align}\label{eq:linear_advection_integrator}
0
\nonumber
&= \dfrac{1}{4} \bigg[ \dfrac{u_{i+1, j-1} - u_{i-1, j-1}}{2 h_{t}} + 2 \, \dfrac{u_{i+1, j} - u_{i-1, j}}{2 h_{t}} + \dfrac{u_{i+1, j+1} - u_{i-1, j+1}}{2 h_{t}} \bigg] \\
&+ \dfrac{c}{4} \bigg[ \dfrac{u_{i-1, j+1} - u_{i-1, j-1}}{2 h_{x}} + 2 \, \dfrac{u_{i, j+1} - u_{i, j-1}}{2 h_{x}} + \dfrac{u_{i+1, j+1} - u_{i+1, j-1}}{2 h_{x}} \bigg] .
\end{align}
As in the continuous case, the discrete adjoint equation has the exact same form as the discrete advection equation.
We see that the derivatives are approximated by second-order centred finite differences. This means that we need initial data at two consecutive points in time, even though the advection equation is first order in time, so that initial data at one point in time should suffice.
This problem is typical for applying the Veselov discretisation to formal Lagrangians.
A simple solution will be provided in section \ref{sec:vorticity_simplifications}.

We also see that the time derivative is averaged in space and the spatial derivative is averaged in time. Under certain conditions, the spatial average of the time derivative can lead to grid-scale oscillations.
This can be seen as follows. The time derivative in (\ref{eq:linear_advection_integrator}) features a spatial average of the form $\tfrac{1}{4} [ 1 \; 2 \; 1 ]$. So if on top of the actual solution there is some oscillation, e.g., of the form $[ -1 \; +1 \; -1 ]$, this is eliminated by the average.
Similarly, the spatial centred finite difference derivative, $\tfrac{1}{2 h_{x}} [-1 \;\; 0 \; +1 ]$, ``does not see'' such an oscillation. For the case of periodic boundary conditions, the only way to prevent the instability to appear is to use an odd number of grid points in $x$.
This is very similar to the phenomenon of checker-boarding often observed in simulations of the incompressible Euler equations.

In the next section we will prove some important discrete conservation laws of the resulting scheme.
This exemplifies that even though a variational integrator has favourable conservation properties it might be unstable, and that not every discretisation of the Lagrangian leads to a good scheme.

\subsubsection*{Discrete Conservation Laws}

For a formal Lagrangian of a single equation, the discrete symmetry condition (\ref{eq:noether_discrete_symmetry_condition_3}) reads
\begingroup
\allowdisplaybreaks
\begin{alignat}{5}
\label{eq:linear_advection_discrete_symmetry_condition_1}
& 0
\nonumber
&{}={}& \dfrac{\partial L_d}{\partial \by_{\square^1}} \Big( \tilde{\phy}_{i,  j  }, \tilde{\phy}_{i,  j+1}, \tilde{\phy}_{i+1,j+1}, \tilde{\phy}_{i+1,j  } \Big) \cdot \eta_{i,  j  }
&{}+{}& \dfrac{\partial L_d}{\partial \tilde{\by}_{\square^1}} \Big( \tilde{\phy}_{i,  j  }, \tilde{\phy}_{i,  j+1}, \tilde{\phy}_{i+1,j+1}, \tilde{\phy}_{i+1,j  } \Big) \cdot \tilde{\eta}_{i,  j  }
\\
\nonumber &
&{}+{}& \dfrac{\partial L_d}{\partial \by_{\square^2}} \Big( \tilde{\phy}_{i,  j  }, \tilde{\phy}_{i,  j+1}, \tilde{\phy}_{i+1,j+1}, \tilde{\phy}_{i+1,j  } \Big) \cdot \eta_{i,  j+1}
&{}+{}& \dfrac{\partial L_d}{\partial \tilde{\by}_{\square^2}} \Big( \tilde{\phy}_{i,  j  }, \tilde{\phy}_{i,  j+1}, \tilde{\phy}_{i+1,j+1}, \tilde{\phy}_{i+1,j  } \Big) \cdot \tilde{\eta}_{i,  j+1}
\\
\nonumber &
&{}+{}& \dfrac{\partial L_d}{\partial \by_{\square^3}} \Big( \tilde{\phy}_{i,  j  }, \tilde{\phy}_{i,  j+1}, \tilde{\phy}_{i+1,j+1}, \tilde{\phy}_{i+1,j  } \Big) \cdot \eta_{i+1,j+1}
&{}+{}& \dfrac{\partial L_d}{\partial \tilde{\by}_{\square^3}} \Big( \tilde{\phy}_{i,  j  }, \tilde{\phy}_{i,  j+1}, \tilde{\phy}_{i+1,j+1}, \tilde{\phy}_{i+1,j  } \Big) \cdot \tilde{\eta}_{i+1,j+1}
\\ &
&{}+{}& \dfrac{\partial L_d}{\partial \by_{\square^4}} \Big( \tilde{\phy}_{i,  j  }, \tilde{\phy}_{i,  j+1}, \tilde{\phy}_{i+1,j+1}, \tilde{\phy}_{i+1,j  } \Big) \cdot \eta_{i+1,j  }
&{}+{}& \dfrac{\partial L_d}{\partial \tilde{\by}_{\square^4}} \Big( \tilde{\phy}_{i,  j  }, \tilde{\phy}_{i,  j+1}, \tilde{\phy}_{i+1,j+1}, \tilde{\phy}_{i+1,j  } \Big) \cdot \tilde{\eta}_{i+1,j  }
.
\end{alignat}
\endgroup
Here and in the following, evaluation of $\eta_{i,j}$ and $\tilde{\eta}_{i,j}$ at the fields $\phy_{i,j}$ and $\tilde{\phy}_{i,j}$ is implied.
For the Lagrangian (\ref{eq:linear_advection_discrete_lagrangian_midpoint}), this becomes
\begingroup
\allowdisplaybreaks
\begin{align}\label{eq:linear_advection_discrete_symmetry_condition_2}
0
\nonumber
&= h_{t} h_{x} \, \dfrac{1}{8} \bigg[ v_{i,j} + v_{i,j+1} + v_{i+1,j+1} + v_{i+1,j} \bigg]
\bigg[
   \dfrac{\eta_{i+1,j  } - \eta_{i,  j  }}{h_{t}}
 + \dfrac{\eta_{i+1,j+1} - \eta_{i,  j+1}}{h_{t}}
\bigg]
\\
\nonumber
&+ h_{t} h_{x} \, \dfrac{c}{8} \bigg[ v_{i,j} + v_{i,j+1} + v_{i+1,j+1} + v_{i+1,j} \bigg]
\bigg[
   \dfrac{\eta_{i,  j+1} - \eta_{i,  j  }}{h_{x}}
 + \dfrac{\eta_{i+1,j+1} - \eta_{i+1,j  }}{h_{x}}
\bigg]
\\
\nonumber
&+ h_{t} h_{x} \, \dfrac{1}{8} \bigg[
   \tilde{\eta}_{i,  j  }
 + \tilde{\eta}_{i,  j+1}
 + \tilde{\eta}_{i+1,j+1}
 + \tilde{\eta}_{i+1,j  }
\bigg]
\bigg[
   \dfrac{u_{i+1,j  } - u_{i,j  }}{h_{t}}
 + \dfrac{u_{i+1,j+1} - u_{i,j+1}}{h_{t}} 
\bigg]
\\
&+ h_{t} h_{x} \, \dfrac{c}{8} \bigg[ 
   \tilde{\eta}_{i,  j  }
 + \tilde{\eta}_{i,  j+1}
 + \tilde{\eta}_{i+1,j+1}
 + \tilde{\eta}_{i+1,j  }
\bigg]
\bigg[
   \dfrac{u_{i,  j+1} - u_{i,  j}}{h_{x}}
 + \dfrac{u_{i+1,j+1} - u_{i+1,j}}{h_{x}}
\bigg] .
\end{align}
\endgroup
The corresponding discrete conservation law (\ref{eq:noether_discrete_conservation_law_3}) is
\begingroup
\allowdisplaybreaks
\begin{align}\label{eq:linear_advection_discrete_conservation_law_1}
\tilde{\mathcal{J}}_{d} =
\sum \limits_{j=1}^{N_{1}-1} \bigg[
\nonumber
     & \dfrac{\partial L_d}{\partial \by_{\square^3}} \Big( \tilde{\phy}_{i,   j}, \tilde{\phy}_{i,   j+1}, \tilde{\phy}_{i+1, j+1}, \tilde{\phy}_{i+1, j} \Big) \cdot \eta_{i+1,j+1}
\\
\nonumber
+ \, & \dfrac{\partial L_d}{\partial \by_{\square^4}} \Big( \tilde{\phy}_{i,   j}, \tilde{\phy}_{i,   j+1}, \tilde{\phy}_{i+1, j+1}, \tilde{\phy}_{i+1, j} \Big) \cdot \eta_{i+1,j  }
\\
\nonumber
+ \, & \dfrac{\partial L_d}{\partial \tilde{\by}_{\square^3}} \Big( \tilde{\phy}_{i,   j}, \tilde{\phy}_{i,   j+1}, \tilde{\phy}_{i+1, j+1}, \tilde{\phy}_{i+1, j} \Big) \cdot \tilde{\eta}_{i+1,j+1}
\\
+ \, & \dfrac{\partial L_d}{\partial \tilde{\by}_{\square^4}} \Big( \tilde{\phy}_{i,   j}, \tilde{\phy}_{i,   j+1}, \tilde{\phy}_{i+1, j+1}, \tilde{\phy}_{i+1, j} \Big) \cdot \tilde{\eta}_{i+1,j  }
\bigg] = \text{const. for all $i$} .
\end{align}
\endgroup
In contrast to the continuous case (\ref{eq:linear_advection_conservation_law_global}), here we also find contributions due to derivatives of the discrete Lagrangian $L_{d}$ with respect to $\tilde{\by}$.
For the Lagrangian (\ref{eq:linear_advection_discrete_lagrangian_midpoint}), the Noether charge becomes
\begin{align}\label{eq:linear_advection_discrete_conservation_law_2}
\nonumber
\tilde{\mathcal{J}}_{d}
&= h_{t} h_{x} \, \dfrac{1}{8} \sum \limits_{j=1}^{N_{1}-1} \bigg[ v_{i,j} + v_{i,j+1} + v_{i+1,j+1} + v_{i+1,j} \bigg] \bigg[ \dfrac{\eta_{i+1, j+1} + \eta_{i+1, j}}{h_{t}} + c \, \dfrac{\eta_{i+1,j+1} - \eta_{i+1,  j}}{h_{x}} \bigg] \\
\nonumber
&+ h_{t} h_{x} \, \dfrac{1}{4} \sum \limits_{j=1}^{N_{1}-1} \bigg[ \tilde{\eta}_{i+1,j+1} + \tilde{\eta}_{i+1,j} \bigg] \bigg[
  \dfrac{1}{2} \bigg( \dfrac{u_{i+1,j} - u_{i,j}}{h_{t}} + \dfrac{u_{i+1,j+1} - u_{i,j+1}}{h_{t}} \bigg) \\
&\hspace{16em}
+ \dfrac{c}{2} \bigg( \dfrac{u_{i,j+1} - u_{i,j}}{h_{x}} + \dfrac{u_{i+1,j+1} - u_{i+1,j}}{h_{x}} \bigg)
\bigg] .
\end{align}

We begin with mass conservation. The discrete generator (c.f. (\ref{eq:linear_advection_generator_L1})) is
\begin{align}\label{eq:linear_advection_discrete_generator_L1}
\eta_{i,j} &= 1 , &
\tilde{\eta}_{i,j} &= 0 .
\end{align}
The discrete symmetry condition (\ref{eq:linear_advection_discrete_symmetry_condition_2}) is fulfilled for (\ref{eq:linear_advection_discrete_generator_L1}).
The second summation in (\ref{eq:linear_advection_discrete_conservation_law_2}) is zero, hence
\begin{align}\label{eq:noether_discrete_conservation_law_L1_1}
\tilde{\mathcal{J}}_{d} &= h_{x} \, \sum \limits_{j=1}^{N_{1}-1} \dfrac{1}{4} \Big[ v_{i,j} + v_{i,j+1} + v_{i+1,j+1} + v_{i+1,j} \Big]
= \text{const. for all $i$} .
\end{align}
Identifying $\fv_{i,j} = \fu_{i,j}$, we obtain the discrete conservation law for the total mass in the system
\begin{align}\label{eq:noether_discrete_conservation_law_L1_2}
\mathcal{J}_{d} &= \dfrac{h_{x}}{4} \sum \limits_{j=1}^{N_{1}-1} \Big[ u_{i,j} + u_{i,j+1} + u_{i+1,j+1} + u_{i+1,j} \Big]
= \text{const. for all $i$} .
\end{align}
This expression supports grid oscillations.
For the $L^{2}$ norm, the discrete generator (c.f. (\ref{eq:linear_advection_generator_L2})) is
\begin{align}\label{eq:linear_advection_discrete_generator_L2}
\eta_{i,j} &=   \fu_{i,j} , &
\tilde{\eta}_{i,j} &= - \fv_{i,j} .
\end{align}
The discrete symmetry condition (\ref{eq:linear_advection_discrete_symmetry_condition_2}) is again easily checked to be fulfilled.
The corresponding conservation law (\ref{eq:noether_discrete_conservation_law_3}) is $\tilde{\mathcal{J}}_{d} = \text{const. for all $i$}$ with
\begin{align}\label{eq:noether_discrete_conservation_law_L2_1}
\nonumber
\tilde{\mathcal{J}}_{d} 
&= h_{t} h_{x} \, \dfrac{1}{8} \sum \limits_{j=1}^{N_{1}-1}
\bigg[ v_{i,j} + v_{i,j+1} \bigg] 
\bigg[ \dfrac{u_{i+1,j+1} + u_{i+1,  j}}{h_{t}} + c \, \dfrac{u_{i+1,j+1} - u_{i+1,  j}}{h_{x}} \bigg] \\
&+ h_{t} h_{x} \, \dfrac{1}{8} \sum \limits_{j=1}^{N_{1}-1}
\bigg[ v_{i+1,j+1} + v_{i+1,j} \bigg] 
\bigg[ \dfrac{u_{i,j} + u_{i,j+1}}{h_{t}} - c \, \dfrac{u_{i,j+1} - u_{i,j}}{h_{x}} \bigg] .
\end{align}
Upon identifying $\fv_{i,j} = \fu_{i,j}$ and assuming periodic boundary conditions, the discrete Noether charge becomes
\begin{align}\label{eq:noether_discrete_conservation_law_L2_reduced}
\mathcal{J}_{d} &= h_{x} \sum \limits_{j=1}^{N_{1}-1} u_{i,j} \, \bigg[ \dfrac{u_{i+1,j-1} + 2 u_{i+1,j} + u_{i+1,j+1}}{4} + \dfrac{c h_{t}}{2} \dfrac{ u_{i+1,j+1} - u_{i+1,j-1}}{2 h_{x}} \bigg] .
\end{align}
Just as the discrete mass (\ref{eq:noether_discrete_conservation_law_L1_2}), this expression also supports grid oscillations.
It is somewhat unexpected that the discrete conservation law contains terms corresponding to the spatial component $J^{\cx}$ of the Noether current.
We therefore found an example where the discrete conserved quantity differs from the continuous conserved quantity.
Nevertheless, $\mathcal{J}_{d}$ is consistent with the $L^{2}$ norm of $\fu$, as in the limit of $h_{t} \rightarrow 0$, holding $c$ fixed, the additional term vanishes.

\subsubsection{Trapezoidal Discretisation}
\label{sec:linear_advection_leapfrog}

Writing out the discrete Lagrangian (\ref{eq:applications_discrete_lagrangian_leapfrog}) for the linear advection equation, we obtain
\begin{align}\label{eq:linear_advection_discrete_lagrangian_leapfrog_2}
L_{d}^{\mathrm{tr}} \big( \jp^{1} \tilde{\phy} (\square) \big) = h_{t} h_{x} \, \bigg[ &
\nonumber
\dfrac{1}{2} \, \bigg(
  \dfrac{v_{\square^{1}} + v_{\square^{4}}}{2} \dfrac{u_{\square^{4}} - u_{\square^{1}}}{h_{t}}
+ \dfrac{v_{\square^{2}} + v_{\square^{3}}}{2} \dfrac{u_{\square^{3}} - u_{\square^{2}}}{h_{t}}
\bigg) \\
\nonumber
+ &
\dfrac{c}{2} \, \bigg(
  \dfrac{v_{\square^{1}} + v_{\square^{2}}}{2} \dfrac{u_{\square^{2}} - u_{\square^{1}}}{h_{x}}
+ \dfrac{v_{\square^{3}} + v_{\square^{4}}}{2} \dfrac{u_{\square^{3}} - u_{\square^{4}}}{h_{x}}
\bigg)
\bigg],
\end{align}
or in the compact notation of (\ref{eq:applications_definitions_discrete_derivative_trapezoidal_t}-\ref{eq:applications_definitions_discrete_derivative_trapezoidal_x}),
\begin{align}
L_{d}^{\mathrm{tr}} \big( \jp^{1} \tilde{\phy} (\square) \big) = h_{t} h_{x} \, \big[ \overline{vu}_{t} (\square) + c\,\overline{vu}_{x} (\square)
\big] .
\end{align}

\subsubsection*{Variational Integrator}

The discrete Euler-Lagrange equation (\ref{eq:vi_discrete_euler_lagrange_equations}) for the physical variable $\fu_{d}$ reads
\begin{align}\label{eq:linear_advection_integrator_leapfrog}
0 = \dfrac{u_{i+1, j} - u_{i-1, j}}{2 h_{t}} + c \, \dfrac{u_{i, j+1} - u_{i, j-1}}{2 h_{x}} ,
\end{align}
which is the leapfrog scheme.
The discrete Euler-Lagrange equation for the adjoint variable $\fv_{d}$ takes exactly the same form.

\subsubsection*{Discrete Conservation Laws}

The discrete symmetry condition (\ref{eq:linear_advection_discrete_symmetry_condition_1}) applied to (\ref{eq:linear_advection_discrete_lagrangian_leapfrog_2}) becomes
\begingroup
\allowdisplaybreaks
\begin{align}\label{eq:linear_advection_discrete_symmetry_condition_leapfrog}
0 = h_{t} h_{x} \, \bigg[ &
\nonumber
\dfrac{1}{2} \, \bigg(
  \dfrac{v_{i, j  } + v_{i+1, j  }}{2} \dfrac{\eta_{i+1, j  } - \eta_{i, j  }}{h_{t}}
+ \dfrac{v_{i, j+1} + v_{i+1, j+1}}{2} \dfrac{\eta_{i+1, j+1} - \eta_{i, j+1}}{h_{t}}
\bigg) \\
\nonumber
+ \, &
\dfrac{c}{2} \, \bigg(
  \dfrac{v_{i,   j  } + v_{i,   j+1}}{2} \dfrac{\eta_{i,   j+1} - \eta_{i,   j  }}{h_{x}}
+ \dfrac{v_{i+1, j+1} + v_{i+1, j  }}{2} \dfrac{\eta_{i+1, j+1} - \eta_{i+1, j  }}{h_{x}}
\bigg) \\
\nonumber
+ \, &
\dfrac{1}{2} \, \bigg(
  \dfrac{\tilde{\eta}_{i, j  } + \tilde{\eta}_{i+1, j  }}{2} \dfrac{u_{i+1, j  } - u_{i, j  }}{h_{t}}
+ \dfrac{\tilde{\eta}_{i, j+1} + \tilde{\eta}_{i+1, j+1}}{2} \dfrac{u_{i+1, j+1} - u_{i, j+1}}{h_{t}}
\bigg) \\
+ \, &
\dfrac{c}{2} \, \bigg(
  \dfrac{\tilde{\eta}_{i,   j  } + \tilde{\eta}_{i,   j+1}}{2} \dfrac{u_{i,   j+1} - u_{i,   j}}{h_{x}}
+ \dfrac{\tilde{\eta}_{i+1, j+1} + \tilde{\eta}_{i+1, j  }}{2} \dfrac{u_{i+1, j+1} - u_{i+1, j}}{h_{x}}
\bigg)
\bigg] .
\end{align}
\endgroup
The corresponding discrete Noether charge (\ref{eq:linear_advection_discrete_conservation_law_1}) is
\begingroup
\allowdisplaybreaks
\begin{align}\label{eq:linear_advection_discrete_conservation_law_leapfrog}
\tilde{\mathcal{J}}_{d} =
h_{t} h_{x} \, \dfrac{1}{2} \sum \limits_{j=1}^{N_{1}-1} \bigg[
\nonumber
     & \bigg( \dfrac{v_{i,  j+1} + v_{i+1,j+1}}{2 h_{t}} + c \, \dfrac{v_{i+1,j} + v_{i+1,j+1}}{2 h_{x}} \bigg) \cdot \eta_{i+1,j+1} \\
\nonumber
+ \, & \bigg( \dfrac{v_{i,  j  } + v_{i+1,j  }}{2 h_{t}} - c \, \dfrac{v_{i+1,j} + v_{i+1,j+1}}{2 h_{x}} \bigg) \cdot \eta_{i+1,j } \\
\nonumber
+ \, & \bigg( \dfrac{u_{i+1, j+1} - u_{i, j+1}}{h_{t}} + c \, \dfrac{u_{i+1, j+1} - u_{i+1, j}}{h_{x}} \bigg) \cdot \tilde{\eta}_{i+1, j+1} \\
+ \, & \bigg( \dfrac{u_{i+1, j  } - u_{i, j  }}{h_{t}} + c \, \dfrac{u_{i+1, j+1} - u_{i+1, j}}{h_{x}} \bigg) \cdot \tilde{\eta}_{i+1, j  }
\bigg] .
\end{align}
\endgroup
For the discrete generator of mass conservation (\ref{eq:linear_advection_discrete_generator_L1}), the symmetry condition (\ref{eq:linear_advection_discrete_symmetry_condition_leapfrog}) is fulfilled and the discrete Noether charge (\ref{eq:linear_advection_discrete_conservation_law_leapfrog}) amounts to (\ref{eq:noether_discrete_conservation_law_L1_2}).
Similarly, the symmetry condition (\ref{eq:linear_advection_discrete_symmetry_condition_leapfrog}) is fulfilled for the discrete generator (\ref{eq:linear_advection_discrete_generator_L2}). The restriction of the corresponding discrete Noether charge (\ref{eq:linear_advection_discrete_conservation_law_leapfrog}) becomes
\begin{align}\label{eq:noether_discrete_conservation_law_L2_leapfrog}
\mathcal{J}_{d} =
\dfrac{h_{x}}{2} \sum \limits_{j=1}^{N_{1}-1} \bigg[ u_{i,  j+1} \, u_{i+1,j+1} + u_{i,  j  } \, u_{i+1,j } \bigg]
= \text{const. for all $i$} ,
\end{align}
which is indeed different from (\ref{eq:noether_discrete_conservation_law_L2_reduced}) and an immediate discretisation of the $L^{2}$ norm.

\subsubsection{Simplified Implicit Scheme}
\label{sec:linear_advection_simplified_implicit}

The discrete Lagrangian (\ref{eq:applications_discrete_lagrangian_simplified}) reads
\begin{align}\label{eq:linear_advection_discrete_lagrangian_simplified_2}
L_{d}^{\mathrm{mt}} \big( \jp^{1} \tilde{\phy} (\square) \big) = h_{t} h_{x} \, \bigg[ &
\nonumber
\dfrac{1}{2} \, \bigg(
  \dfrac{v_{\square^{1}} + v_{\square^{4}}}{2} \dfrac{u_{\square^{4}} - u_{\square^{1}}}{h_{t}}
+ \dfrac{v_{\square^{2}} + v_{\square^{3}}}{2} \dfrac{u_{\square^{3}} - u_{\square^{2}}}{h_{t}}
\bigg) \\
\nonumber
+ \, &
\dfrac{c}{8} \,
\bigg( v_{\square^{1}} + v_{\square^{2}} + v_{\square^{3}} + v_{\square^{4}} \bigg )
\bigg( \dfrac{u_{\square^{2}} - u_{\square^{1}}}{h_{x}} + \dfrac{u_{\square^{3}} - u_{\square^{4}}}{h_{x}}
\bigg)
\bigg] ,
\end{align}
or equivalently
\begin{align}
L_{d}^{\mathrm{mt}} \big( \jp^{1} \tilde{\phy} (\square) \big) = h_{t} h_{x} \, \big[ &
\overline{vu}_{t} (\square) + c\, \overline{v} (\square) \, u_{x} (\square)
\big] .
\end{align}

\subsubsection*{Variational Integrator}

The discrete Euler-Lagrange equations (\ref{eq:vi_discrete_euler_lagrange_equations}) for the physical variable are obtained in the form
\begin{align}\label{eq:linear_advection_integrator_simplified}
0 = \dfrac{u_{i+1, j} - u_{i-1, j}}{2 h_{t}} + \dfrac{c}{4} \bigg[ \dfrac{u_{i-1, j+1} - u_{i-1, j-1}}{2 h_{x}}
+ 2 \, \dfrac{u_{i, j+1} - u_{i, j-1}}{2 h_{x}} + \dfrac{u_{i+1, j+1} - u_{i+1, j-1}}{2 h_{x}} \bigg] .
\end{align}
As becomes already apparent by comparing the discrete Lagrangian (\ref{eq:linear_advection_discrete_lagrangian_simplified_2}) with (\ref{eq:linear_advection_discrete_lagrangian_leapfrog_2}) and (\ref{eq:linear_advection_discrete_lagrangian_midpoint}), the time derivative is discretised in the same way as in the trapezoidal scheme (\ref{eq:linear_advection_integrator_leapfrog}) and the spatial derivative is discretised in the same way as in the midpoint scheme (\ref{eq:linear_advection_integrator}).
Again, the discrete Euler-Lagrange equation for the adjoint variable takes the same form so that the system is self-adjoint.

\subsubsection*{Discrete Conservation Laws}

In this case, the discrete symmetry condition (\ref{eq:linear_advection_discrete_symmetry_condition_1}) can be written as
\begingroup
\allowdisplaybreaks
\begin{align}\label{eq:linear_advection_discrete_symmetry_condition_simplified_implicit}
0 = h_{t} h_{x} \, \bigg[ &
\nonumber
\dfrac{1}{2} \, \bigg(
  \dfrac{v_{i, j  } + v_{i+1, j  }}{2} \dfrac{\eta_{i+1, j  } - \eta_{i, j  }}{h_{t}}
+ \dfrac{v_{i, j+1} + v_{i+1, j+1}}{2} \dfrac{\eta_{i+1, j+1} - \eta_{i, j+1}}{h_{t}}
\bigg) \\
\nonumber
+ \, &
\dfrac{c}{8} \, \bigg( v_{i,j} + v_{i,j+1} + v_{i+1,j+1} + v_{i+1,j} \bigg)
\bigg(
   \dfrac{\eta_{i,  j+1} - \eta_{i,  j  }}{h_{x}}
 + \dfrac{\eta_{i+1,j+1} - \eta_{i+1,j  }}{h_{x}}
\bigg) \\
\nonumber
+ \, &
\dfrac{1}{2} \, \bigg(
  \dfrac{\tilde{\eta}_{i, j  } + \tilde{\eta}_{i+1, j  }}{2} \dfrac{u_{i+1, j  } - u_{i, j  }}{h_{t}}
+ \dfrac{\tilde{\eta}_{i, j+1} + \tilde{\eta}_{i+1, j+1}}{2} \dfrac{u_{i+1, j+1} - u_{i, j+1}}{h_{t}}
\bigg) \\
+ \, &
\dfrac{c}{8} \, \bigg( \tilde{\eta}_{i,  j  } + \tilde{\eta}_{i,  j+1} + \tilde{\eta}_{i+1,j+1} + \tilde{\eta}_{i+1,j  } \bigg)
\bigg(
  \dfrac{u_{i,  j+1} - u_{i,  j}}{h_{x}}
+ \dfrac{u_{i+1,j+1} - u_{i+1,j}}{h_{x}}
\bigg)
\bigg] ,
\end{align}
\endgroup
which yields the discrete Noether charge, (c.f., (\ref{eq:linear_advection_discrete_conservation_law_1})),
\begingroup
\allowdisplaybreaks
\begin{align}\label{eq:linear_advection_discrete_conservation_law_simplified_implicit}
\tilde{\mathcal{J}}_{d} =
h_{t} h_{x} \, \dfrac{1}{2} \sum \limits_{j=1}^{N_{1}-1} \bigg[
\nonumber
     & \bigg( \dfrac{v_{i,  j+1} + v_{i+1,j+1}}{2 h_{t}} + c \, \dfrac{v_{i,j} + v_{i,j+1} + v_{i+1,j} + v_{i+1,j+1}}{4 h_{x}} \bigg) \cdot \eta_{i+1,j+1} \\
\nonumber
+ \, & \bigg( \dfrac{v_{i,  j  } + v_{i+1,j  }}{2 h_{t}} - c \, \dfrac{v_{i,j} + v_{i,j+1} + v_{i+1,j} + v_{i+1,j+1}}{4 h_{x}} \bigg) \cdot \eta_{i+1,j } \\
\nonumber
+ \, & \bigg( \dfrac{u_{i+1, j+1} - u_{i, j+1}}{h_{t}} + \dfrac{c}{4} \, \dfrac{u_{i,j+1} - u_{i,j}}{h_{x}} + \dfrac{c}{4} \, \dfrac{u_{i+1, j+1} - u_{i+1, j}}{h_{x}} \bigg) \cdot \tilde{\eta}_{i+1, j+1} \\
+ \, & \bigg( \dfrac{u_{i+1, j  } - u_{i, j  }}{h_{t}} + \dfrac{c}{4} \, \dfrac{u_{i,j+1} - u_{i,j}}{h_{x}} + \dfrac{c}{4} \, \dfrac{u_{i+1, j+1} - u_{i+1, j}}{h_{x}} \bigg) \cdot \tilde{\eta}_{i+1, j  }
\bigg] .
\end{align}
\endgroup
The symmetry condition (\ref{eq:linear_advection_discrete_symmetry_condition_simplified_implicit}) is satisfied for both, mass conservation and the $L^{2}$ norm. The corresponding discrete conservation laws reduce to (\ref{eq:noether_discrete_conservation_law_L1_2}) and (\ref{eq:noether_discrete_conservation_law_L2_leapfrog}), respectively.

\subsubsection{Dissipation and Dispersion}

In this section we want to carry out a simple von Neumann analysis of some of the variational integrators we have derived in order to analyse their dissipation and dispersion behaviour.

The oscillatory function $u(t,x) = e^{- \imath \omega t + \imath k x}$ is a solution of the advection equation (\ref{eq:linear_advection_equation}) if and only if $(\omega, k)$ satisfies the exact dispersion relation $\omega = ck$.
When collocated at grid points $t_i = i h_t$, $x_j = j h_x$, this amounts to
\begin{align}
  u_{ij} = e^{- \imath (\tau i - \xi j)},
\end{align}
with $\tau = \omega h_t$ and $\xi = k h_x$. In analogy with the continuous case, the collocated plane wave $u_{ij}$ is a solution of the numerical scheme if and only if $(\tau, \xi)$ satisfies the numerical dispersion relation. 
For the case of the Veselov scheme (\ref{eq:linear_advection_integrator}), the dispersion relation can be written as
\begin{align}\label{eq:linear_advection_dispersion_veselov}
  \sin (\tau) \big[ 1 + \cos(\xi) \big] - \dfrac{c}{c_{\text{grid}}} \big[ 1 + \cos(\tau) \big] \sin(\xi)
  = 0 ,
\end{align}
where $c_{\text{grid}} = h_x / h_t$.
This is the same dispersion relation as the one obtained by \citet{AscherMcLachlan:2004} for the box scheme.
For comparison, we consider the simplified implicit scheme (\ref{eq:linear_advection_integrator_simplified}), for which the dispersion relation is
\begin{align}\label{eq:linear_advection_dispersion_simplified}
  \sin (\tau) - \dfrac{1}{2} \dfrac{c}{c_{\text{grid}}} \big[ 1 + \cos(\tau) \big] \sin(\xi)
  = 0 ,
\end{align}
as well as the leap frog scheme (\ref{eq:linear_advection_integrator_leapfrog}), for which the dispersion relation is
\begin{align}\label{eq:linear_advection_dispersion_leapfrog}
  \sin(\tau) - \dfrac{c}{c_{\text{grid}}} \sin(\xi) = 0.
\end{align}

\begin{figure}[t]
  \includegraphics[width=\textwidth]{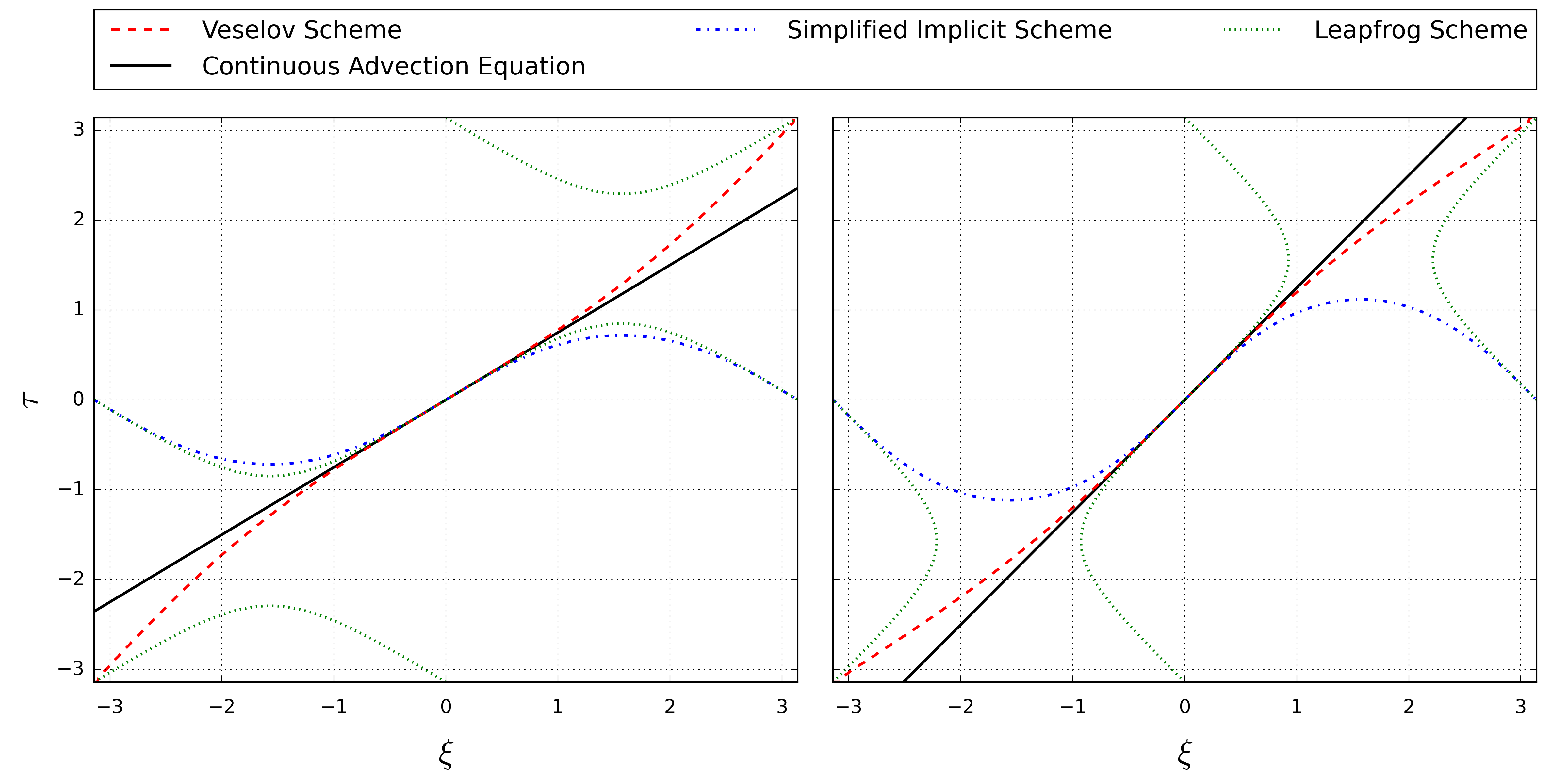}
  \caption{\label{fig:linear_advection_dispersion_relation} Numerical dispersion
    relation of the Veselov scheme (red dashed line), the simplified implicit
    scheme (blue dash-dotted line) and the leap frog scheme (green dotted line)
    applied to the linear advection equation with $c/c_{\text{grid}} = 0.75$
    (left) and $c/c_{\text{grid}} = 1.25$ (right). For comparison the exact
    dispersion relation is also shown (black solid line).} 
\end{figure}

One can observe that in all cases the left-hand side is real valued, i.e., numerical dissipation is exactly zero, but there is numerical dispersion. Figure~\ref{fig:linear_advection_dispersion_relation} shows the solution of the numerical dispersion relation in comparison with the exact linear dispersion of the advection equation in the $(\xi, \tau)$-plane for different values of the ratio $c/c_{\text{grid}}$. 
As usual, for $(\xi,\tau)$ near the origin, the numerical dispersion is a good approximation of the exact linear relation, and thus the corresponding harmonics have a phase velocity close to the exact value $c$. On the other hand, when the wave number is larger than half of the grid wave number $1/h_x$, the frequency of the harmonic deviates from the exact dispersion relation. The corresponding harmonics have a phase velocity different from the exact value $c$, leading to numerical dispersion.

For the simplified implicit scheme, one can also observe local maxima of the numerical dispersion, which correspond to a zero numerical group velocity, while the exact value for the advection equation is equal to the phase velocity $c$.

The left-hand-side panel of Figure~\ref{fig:linear_advection_dispersion_relation} shows the dispersion relation for the case $c/c_{\text{grid}} < 1$. We find that the leapfrog scheme supports parasitic modes, that is to each $\xi \in [-\pi, +\pi]$ there are two values of $\tau$. The high-frequency branches are called parasitic modes. The Veselov and the simplified implicit scheme also feature such modes at the $\tau = -\pi$ and $\tau = +\pi$ lines. As we will see in the numerical experiments, these branches will be populated only if the spectrum of the solution is sufficiently broad to cover the nonlinear region of the dispersion relation, that is for waves which are not well represented by the chosen grid parameters.

For didactical purposes, we have also considered the case $c/c_{\text{grid}} > 1$ on the right-hand-side panel of Figure~\ref{fig:linear_advection_dispersion_relation}. One can see that the numerical group velocity of the leapfrog scheme becomes infinite, e.g., for $\xi = \pm \sin^{-1} (0.8)$, which renders the scheme unstable. This is always observed for $c/c_{\text{grid}} > 1$ which corresponds to the Courant-Friedrichs-Levy (CFL) stability criterion that limits the timestep for explicit schemes. The implicit variational integrators do not have such a limitation. In the following examples, the timestep is chosen small enough so that these instabilities do not appear.

\subsubsection{Numerical Examples}

We consider two test cases, a sum of cosines and a Gaussian. The first is purely meant to verify the dispersion relation, whereas the second is used to judge the quality of the solution.
In both cases, the domain is $\Omega = [-0.5, +0.5)$ with periodic boundary conditions, the timestep is $h_t = 2.5 \times 10^{-3}$, the phase velocity is set to $c=1$, and the number of points in space and time is $n_x = 255$ and $n_t = 4000$, respectively. Hence, the time interval is $[0,10]$, which corresponds to ten passes of the wave packet through the domain.

The following tests are implemented in Python \cite{Python, Langtangen:2014} using NumPy \cite{vanDerWalt:2011} and the sparse linear algebra module of SciPy \cite{SciPy}. Visualisation was done using matplotlib \cite{Matplotlib, Hunter:2007}.

\begin{figure}[htb]
	\centering
	\subfloat[][Veselov Scheme]{
	\includegraphics[width=.3\textwidth]{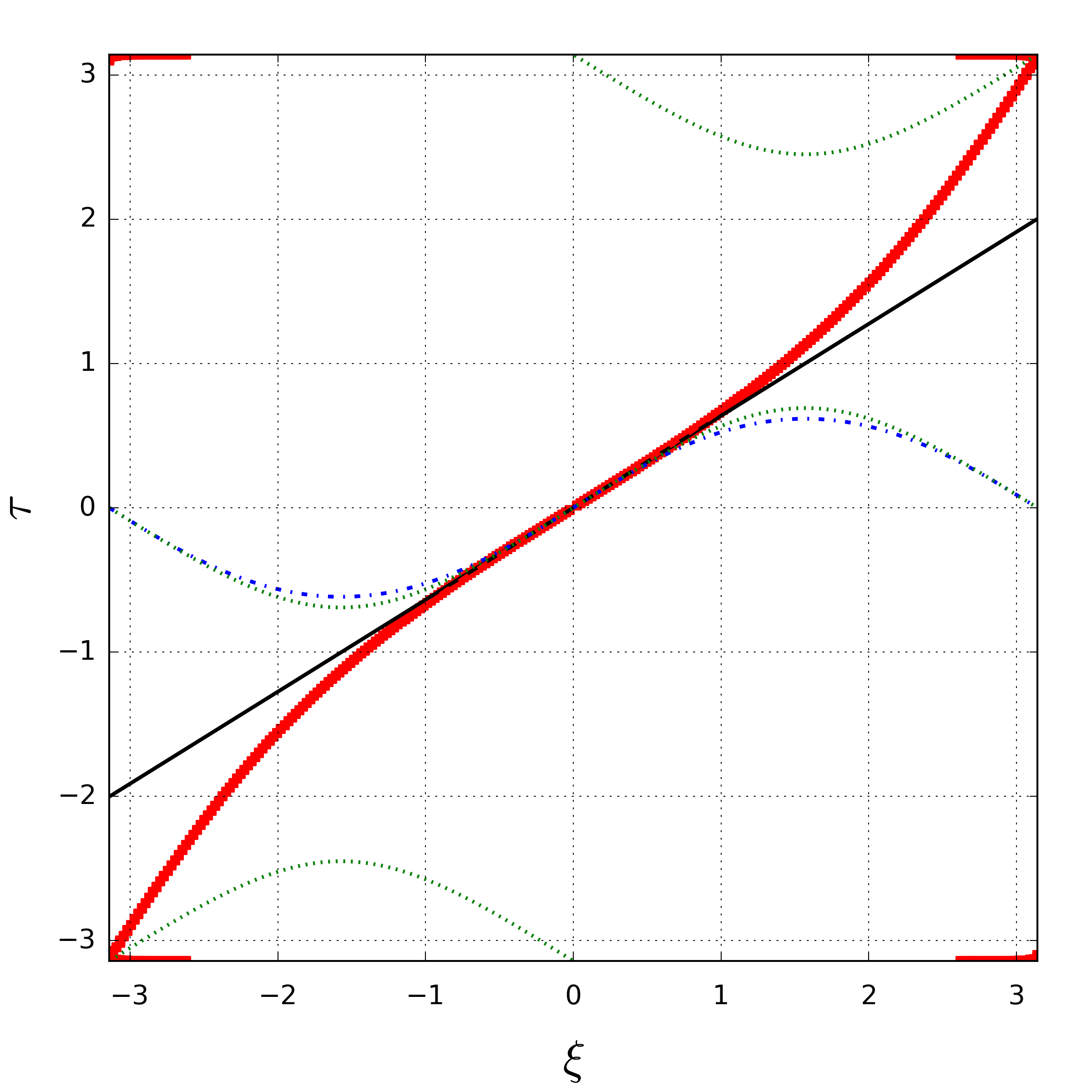}
	}
	\subfloat[][Simplified Implicit Scheme]{
	\includegraphics[width=.3\textwidth]{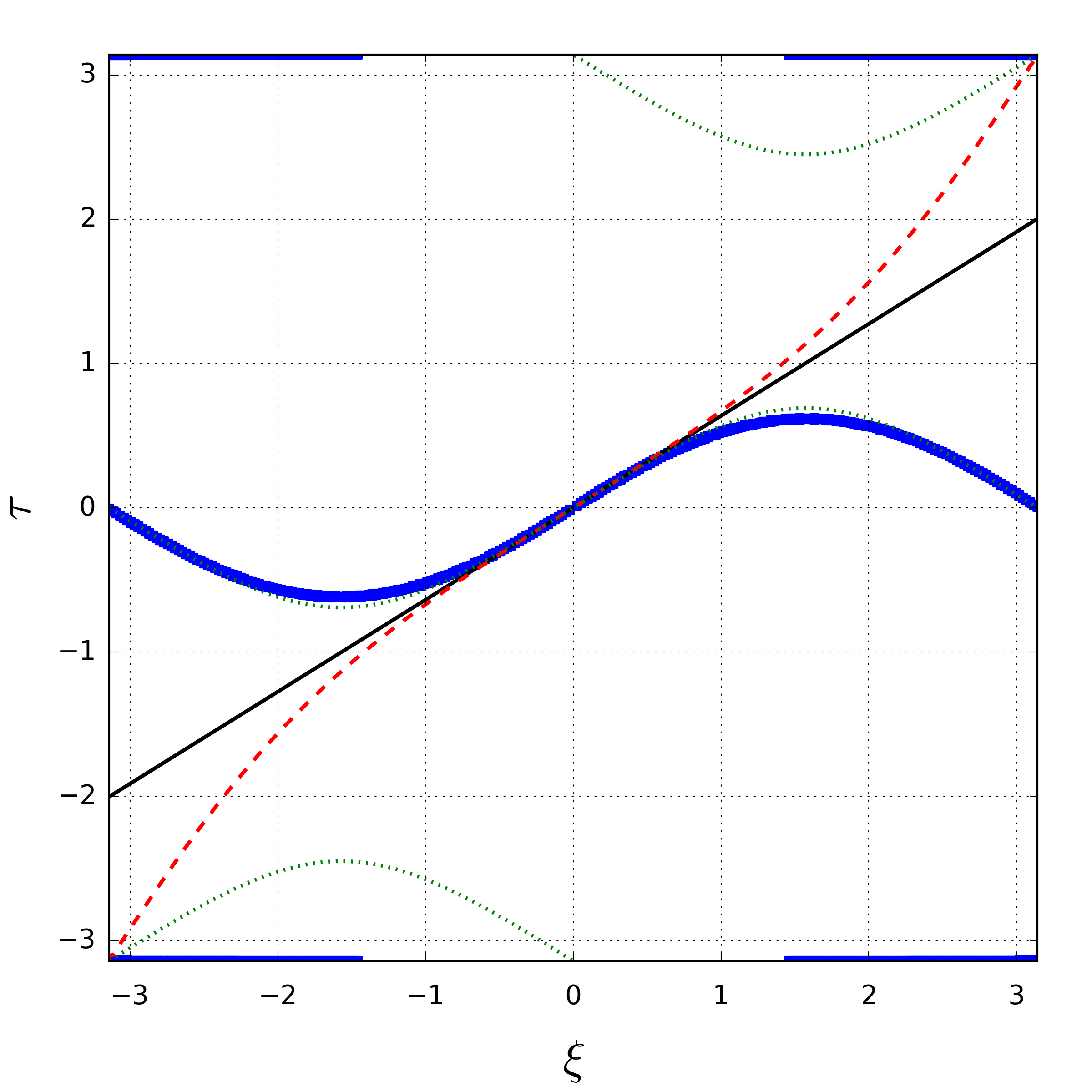}
	}
	\subfloat[][Leapfrog Scheme]{
	\includegraphics[width=.3\textwidth]{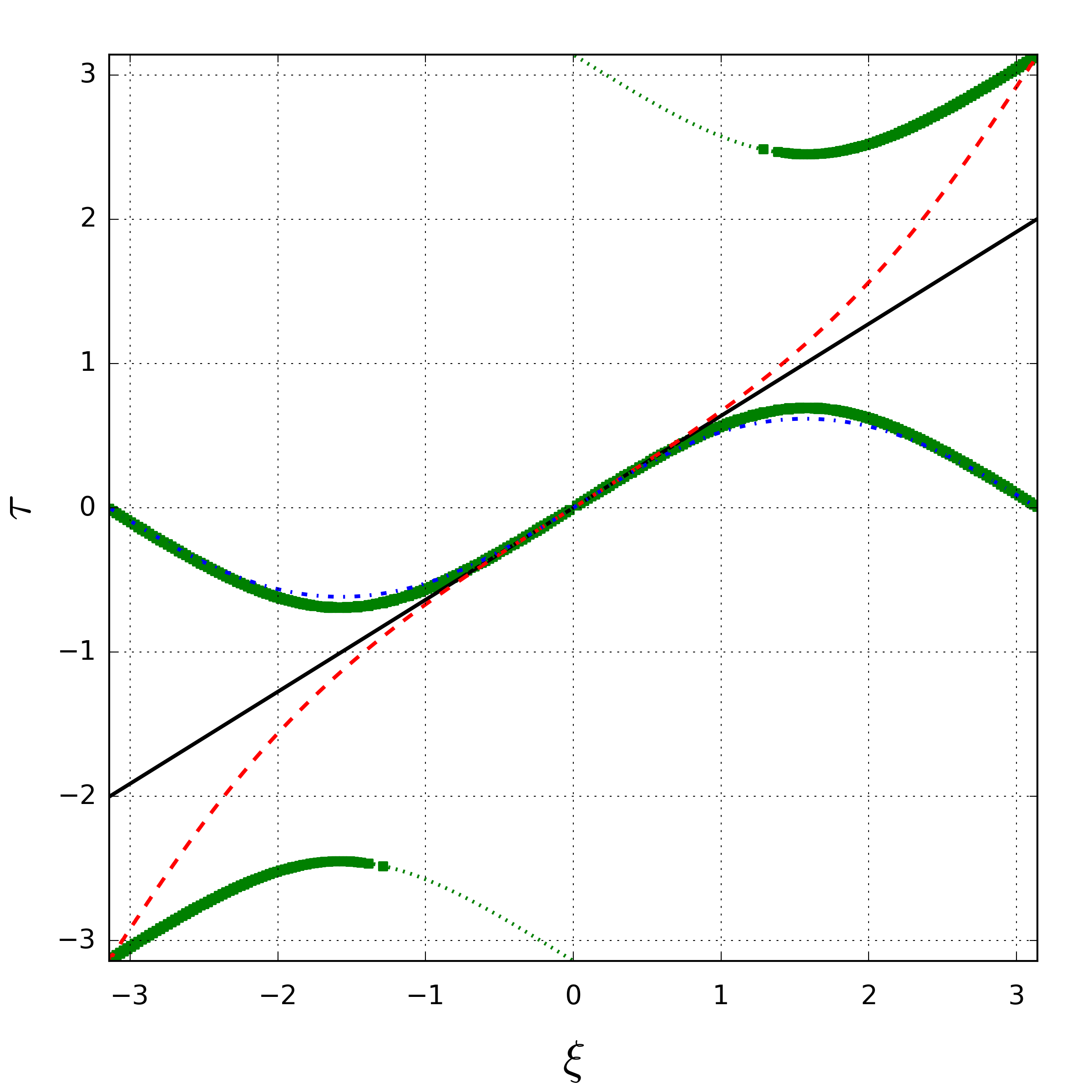}
	}
	\caption{Theoretical and experimental dispersion relation for the sum of
          cosines. The experimental values lie right on top of the theoretical
          curve. Red dashed line: Veselov scheme. Blue dash-dotted line:
          Simplified implicit scheme. Green dotted line: Leapfrog scheme. Black
          solid line: Analytical dispersion relation.}
	\label{fig:linear_advection_cosines_dispersion}
\end{figure}

\subsubsection*{Experimental Dispersion Relation}

In this example, we initialise $u$ with a sum of cosines,
\begin{align}
u_{C,0} (x) = \sum \limits_{i=1}^{n_x/2} \cos (i 2 \pi x) ,
\end{align}
which excites all the modes supported by the grid in order to experimentally verify the dispersion relations (\ref{eq:linear_advection_dispersion_veselov}-\ref{eq:linear_advection_dispersion_leapfrog}) from the previous section.
Figure~\ref{fig:linear_advection_cosines_dispersion} shows the theoretical and experimental dispersion relations for the Veselov scheme, the simplified implicit scheme and the leapfrog scheme. The theoretical dispersion relations are well matched by the simulation, including some of the parasitic modes.

\subsubsection*{Gaussian Wave}

\begin{figure}[p]
	\centering
	\subfloat[][Veselov Scheme]{
 	\includegraphics[width=.63\textwidth]{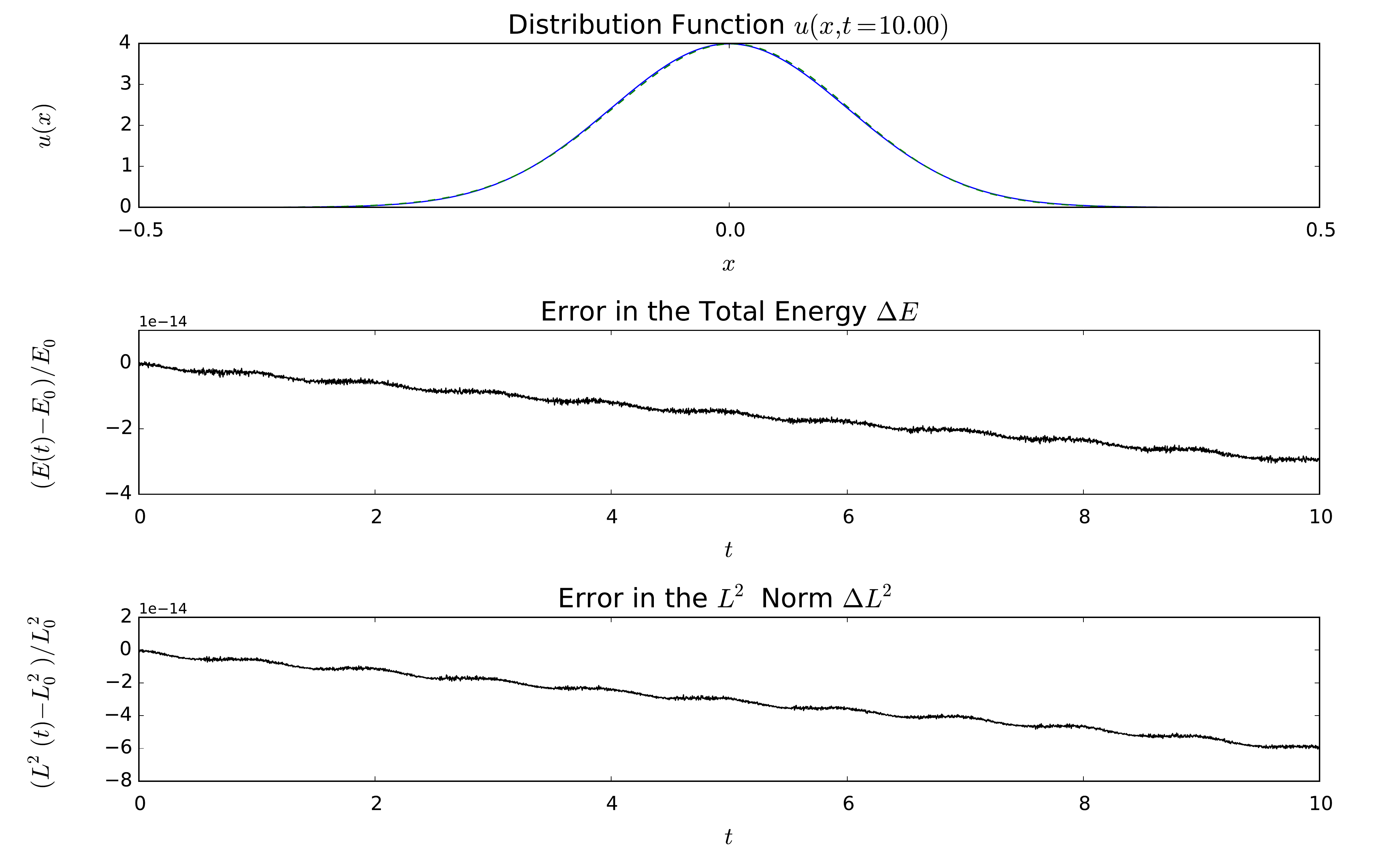}
	}
	
	\subfloat[][Simplified Implicit Scheme]{
 	\includegraphics[width=.63\textwidth]{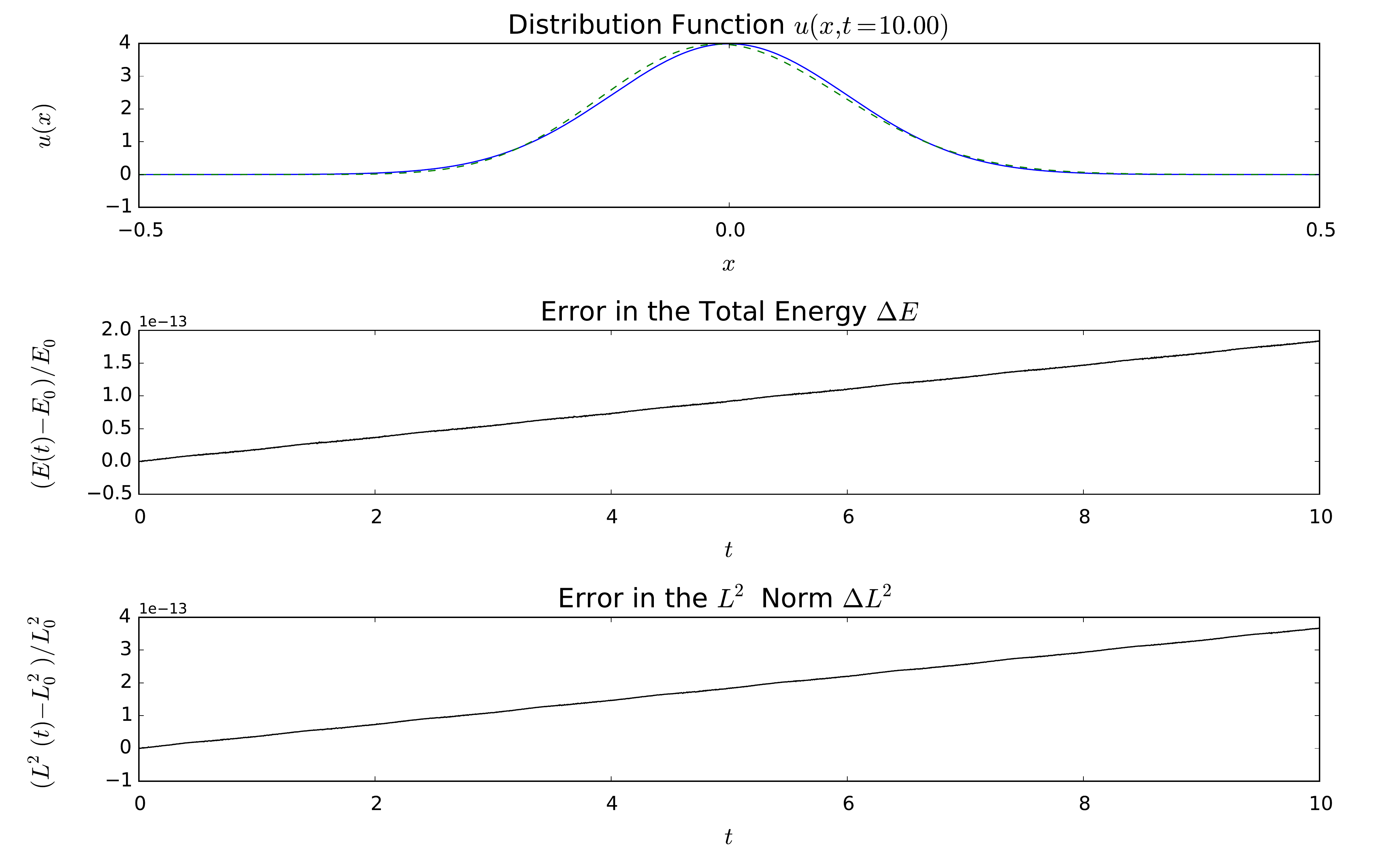}
	}
	
	\subfloat[][Leapfrog Scheme]{
	\includegraphics[width=.63\textwidth]{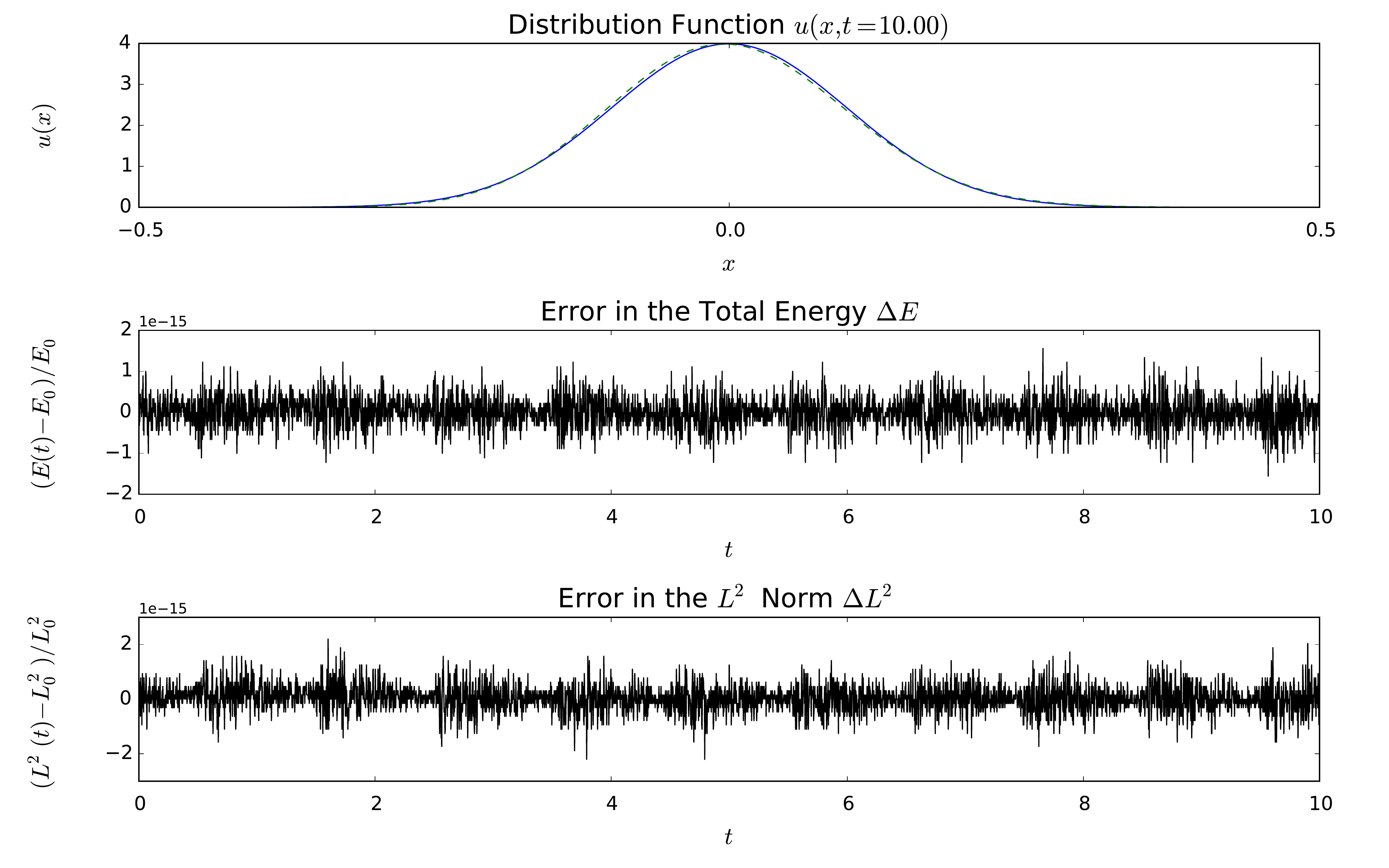}
	}
 	\caption{Gaussian after ten passes with the Veselov scheme, the simplified implicit scheme and the leapfrog scheme. Comparing the solution (green dashed curve) to the initial condition (blue curve) the effect of dispersion becomes visible. Energy and the $L^{2}$ norm are well preserved. For the Veselov scheme as well as the simplified implicit scheme there is some drift in the errors, but this is extremely small.}
	\label{fig:linear_advection_gaussian_timetraces}
\end{figure}

In this section, we consider a Gaussian,
\begin{align}
u_{G,0} (x) = \dfrac{1}{\sigma \sqrt{2 \pi}} \, \exp \bigg( - \dfrac{1}{2} \bigg[ \dfrac{x}{\sigma} \bigg]^2 \bigg) ,
\end{align}
with $\sigma = 0.1$, so that the spectrum is well within the linear region of the dispersion relation.
Figure~\ref{fig:linear_advection_gaussian_timetraces} shows that after 10 passes the form of the Gaussian is close to the initial one for both, the simplified implicit variational integrator (\ref{eq:linear_advection_integrator_simplified}) and the leapfrog scheme (\ref{eq:linear_advection_integrator_leapfrog}), and seems to exactly match the initial conditions for the Veselov scheme (\ref{eq:linear_advection_integrator}). The leapfrog scheme shows slightly less distortion than the simplified implicit scheme, which is explained by their dispersion relations, which is closer to the analytic dispersion relation for the leapfrog scheme and even more so for the Veselov scheme.
For all three integrators, conservation of energy and the $L^{2}$ norm is also shown in Figure~\ref{fig:linear_advection_gaussian_timetraces}.

\subsection{Vorticity Equation}
\label{sec:vorticity}

Next we consider an example with more than one equation and more than one spatial dimension, the vorticity equation,
\begin{align}\label{eq:vorticity_equations}
\omega_{\ct} + \{ \psi , \omega \} &= 0 , &
\omega &= \Delta \psi .
\end{align}
Here, $\psi$ is the streaming function and $\omega$ is the vorticity.
The vorticity equation arises by computing the curl of the incompressible Euler equation in two dimensions,
\begin{align}\label{eq:navier_stokes_equation}
u_{\ct} + ( u \cdot \nabla ) u + \nabla p &= 0 ,
\end{align}
where $u$ is the fluid velocity and $p$ the pressure.
Specifically, $\psi$ and $\omega$ are related to the velocity $u$ by $u = \nabla_{\perp} \psi$ and $\omega = \nabla_{\perp} \cdot u$ with $\nabla_{\perp} = ( - \partial_{y}, \partial_{x})^{T}$.
The Poisson bracket $\{ \cdot , \cdot \}$ is defined by
\begin{align}\label{eq:vorticity_brackets}
\{ \psi , \omega \} &= \psi_{\cx} \omega_{\cy} - \psi_{\cy} \omega_{\cx} ,
\end{align}
which coincides with the determinant of the Jacobian of the transformation $(\cx, \cy) \mapsto (\psi , \omega)$. For this reason, the 2-dimensional Poisson bracket is also referred to as Jacobian.
The variational discretisation of (\ref{eq:vorticity_equations}) leads to a very interesting result, namely Arakawa's famous discretisation of the Poisson bracket \cite{Arakawa:1966}.

\subsubsection{Formal Lagrangian}

The formal Lagrangian of the vorticity equation (\ref{eq:vorticity_equations}) is
\begin{align}\label{eq:vorticity_lagrangian}
L ( \jp^{2} \tilde{\phy} ) = \zeta \, ( \omega_{\ct} + \{ \psi , \omega \} ) + \chi ( \omega - \Delta \psi ) ,
\end{align}
with solution vector $(\omega, \psi, \zeta, \chi)$.
Here, the formal Lagrangian $L$ is defined on the second jet-bundle $\jb^{2} \tilde{Y}$ due to the Laplace operator. In order to avoid second order derivatives, we apply a symmetrisation (corresponding to an integration by parts in the action functional) in the term that produces the Laplacian. This yields an equivalent Lagrangian, which however is defined on the first jet-bundle $\jb^{1} \tilde{Y}$.

Some care has to be taken when discretising the Poisson bracket (see \citeauthor{SalmonTalley:1989} \cite{SalmonTalley:1989}).
In order to retain the antisymmetry property of the continuous bracket at the discrete level, a symmetrisation has to be introduced in the Lagrangian.
Using integration by parts while assuming appropriate boundary conditions, it is seen that the cyclic permutations of the functions in the integrand are all identical\footnote{
This holds assuming boundary conditions such that the surface terms of the integration vanish, e.g., periodic or homogeneous Dirichlet.
},
\begin{align}
  \int \zeta \, \{ \psi , \omega \} \, d\cx \, d\cy
= \int \psi \, \{ \omega , \zeta \} \, d\cx \, d\cy
= \int \omega \, \{ \zeta , \psi \} \, d\cx \, d\cy .
\end{align}
Instead of selecting one of those equivalent forms, a convex combination can be considered, namely,
\begin{align}
\int \zeta \, \{ \psi , \omega \} \, d\cx \, d\cy
&= \int \Big[ \alpha \, \zeta \, \{ \psi , \omega \}
             + \beta \, \psi \, \{ \omega , \zeta \}
            + \gamma \, \omega \, \{ \zeta , \psi \} \Big] \, d\cx \, d\cy ,
\end{align}
with $\alpha + \beta + \gamma = 1$.
One finds that the symmetric case, $\alpha = \beta = \gamma = 1/3$, is the one that retains the properties of the bracket at the discrete level.
We therefore use the modified Lagrangian
\begin{align}\label{eq:vorticity_lagrangian_symmetrised}
L' ( \jp^{1} \tilde{\phy} ) = \zeta \omega_{\ct} + \dfrac{1}{3} \Big[ \zeta \, \{ \psi , \omega \} + \psi \, \{ \omega , \zeta \} + \omega \, \{ \zeta , \psi \} \Big] + \chi \omega + \nabla \chi \cdot \nabla \psi .
\end{align}

Computing the variational derivatives of the action with respect to the adjoint variables, we obtain (\ref{eq:vorticity_equations}). The variational derivatives with respect to the physical variables give
\begin{align}
\zeta_{\ct} + \{ \psi , \zeta \} &= \chi , &
\Delta \chi &= \{ \omega , \zeta \} .
\end{align}
This extended system of equations is obviously not self-adjoint, but we can still find a simple compatible solution of the adjoint equations, namely $(\zeta, \chi) = \phi(\omega, \psi) = (\omega, 0)$.
Assume we can set $\zeta = \omega$, then the Poisson bracket in the adjoint Poisson equation vanishes, which implies that $\chi$ is a harmonic function. In particular, we can choose $\chi = 0$, so that the right-hand side of the equation for $\zeta$ becomes zero and $\zeta$ fulfils the same equation as $\omega$. This justifies the choice of $\phi$.

\subsubsection{Continuous Conservation Laws}
\label{sec:vorticity_noether}

In order to determine the conservation laws for the vorticity equation, it is not enough to just consider vertical transformations as in (\ref{sec:vi_continous_noether_theorem}). 
We have to consider also horizontal transformations, as the interesting symmetries of the vorticity equation are not generated by purely vertical vector fields.
To do so, the discrete Noether theorem as it is presented in section \ref{sec:vi_discrete_noether_theorem} is not sufficiently general.
Unfortunately, it is not straight forward to consider symmetry generators with horizontal components in the framework of finite difference discretisations of the Lagrangian as we use it here.
For that reason, we follow a simplified approach to determine the continuous conservation laws for the vorticity equation which is amenable to discretisation. 
A rigorous analysis of the continuous and discrete conservation laws for the vorticity equation is postponed to a subsequent exposition.

In the following, we will consider the vorticity equation for $\omega$ alone,
\begin{align}\label{eq:vorticity_equation_conservation_laws}
\mathcal{F} [\omega] = \omega_{\ct} + \{ \psi , \omega \} &= 0 ,
\end{align}
without the Poisson equation, assuming a constant-in-time streaming function $\psi$. 
The Lagrangian for this system reads
\begin{align}\label{eq:vorticity_lagrangian_conservation_laws}
\bar{L} = \tilde{\by} \, ( \bz_{\ct} + \psi_{\cx} \bz_{\cy} - \psi_{\cy} \bz_{\cx} ) .
\end{align}
It is easily seen that Euler-Lagrange equations
(\ref{eq:variational_continuous_euler_lagrange}) for
(\ref{eq:vorticity_lagrangian_conservation_laws}) amount to (\ref{eq:vorticity_equation_conservation_laws}) and 
\begin{align}
\zeta_{\ct} + \{ \psi , \zeta \} = 0 ,
\end{align}
such that the extended system is self-adjoint and we can construct a solution for $\zeta$ simply by $\phi(\omega) = \omega$.
The generating vector field for vertical transformations of solutions of the vorticity equation (\ref{eq:vorticity_equation_conservation_laws}) and its first jet prolongation are written as
\begin{align}
V &= \eta \,  \dfrac{\partial}{\partial \by} &
& \text{and} &
j^{1} V &= \eta \dfrac{\partial}{\partial \by} + \eta_{\mu} \dfrac{\partial}{\partial \bz_{\mu}} .
\end{align}
The action of a general vector field $V$ on $\mathcal{F}$ is given by
\begin{align}\label{eq:voriticity_noether_symmetry_equation}
j^{1} V (\mathcal{F}) &= \eta_{\ct} + \psi_{\cx} \eta_{\cy} - \psi_{\cy} \eta_{\cx} .
\end{align}
This is used to determine the component $\tilde{\eta}$ of the corresponding extended vector field $\tilde{V}$ and its first jet prolongation $\jp^{1} \tilde{V}$,
\begin{align}
\tilde{V} &= \eta \dfrac{\partial}{\partial \by} + \tilde{\eta} \dfrac{\partial}{\partial \tilde{\by}} &
& \text{and} &
j^{1} \tilde{V} &= \eta \dfrac{\partial}{\partial \by} + \eta_{\mu} \dfrac{\partial}{\partial \bz_{\mu}} + \tilde{\eta} \dfrac{\partial}{\partial \tilde{\by}} + \tilde{\eta}_{\mu} \dfrac{\partial}{\partial \tilde{\bz}_{\mu}} .
\end{align}
Its action on the Lagrangian (\ref{eq:vorticity_lagrangian_conservation_laws}) is
\begin{align}
\jp^{1} \tilde{V} (\bar{L})
= \tilde{\by}  \, ( \eta_{\ct} + \psi_{\cx} \eta_{\cy} - \psi_{\cy} \eta_{\cx} )
+ \tilde{\eta} \, ( \bz_{\ct}  + \psi_{\cx} \bz_{\cy}  - \psi_{\cy} \bz_{\cx}  ) .
\end{align}

For conservation of circulation (``mass'') we have $\eta = 1$, and therefore
\begin{align}
j^{1} V &= \dfrac{\partial}{\partial \by} .
\end{align}
The action of the prolonged vector field on (\ref{eq:vorticity_equation_conservation_laws}) becomes
\begin{align}
j^{1} V (\mathcal{F}) &= 0 ,
\end{align}
so that $\tilde{\eta} = 0$ and
\begin{align}
j^{1} \tilde{V} &= \dfrac{\partial}{\partial \by} .
\end{align}
The action of the prolonged vector field on the Lagrangian (\ref{eq:vorticity_lagrangian_conservation_laws}) is
\begin{align}
\jp^{1} \tilde{V} (\bar{L} ) = \dfrac{\partial \bar{L}}{\partial \by} = 0 ,
\end{align}
that is we have a symmetry, and the corresponding conservation law (Noether charge) reads
\begin{align}
\mathcal{\tilde{J}}
= \int \tilde{J}^{t} \, d\cx
= \int \dfrac{\partial L}{\partial \bz_{t}} \big( \jp^{1} \tilde{\phy} \big) \cdot \eta (\tilde{\phy}) \, d\cx
= \int \zeta \, d\cx
= \text{const. for all $\ct$} .
\end{align}
Upon restricting the Noether current $\tilde{J}$ with $\zeta = \phi (\omega) = \omega$, this becomes conservation of the total circulation in the system,
\begin{align}
\mathcal{J}
= \int \big( \tilde{J}^{t} \circ \jp^{1} \Phi \big) (\jp^{1} \phy) \, d\cx
= \int \omega \, d\cx
= \text{const. for all $\ct$} .
\end{align}
As here, $\psi$ is treated as constant in time, energy conservation is obtained in the same way for $\eta = \tfrac{1}{2} \psi$.

Next we consider enstrophy, for which $\eta = \by$, so that the first jet prolongation of the generating vector field,
\begin{align}
j^{1} V
&= \by \dfrac{\partial}{\partial \by} 
 + \bz_{\ct} \dfrac{\partial}{\partial \bz_{\ct}} 
 + \bz_{\cx} \dfrac{\partial}{\partial \bz_{\cx}} 
 + \bz_{\cy} \dfrac{\partial}{\partial \bz_{\cy}} ,
\end{align}
acts on the vorticity equation (\ref{eq:vorticity_lagrangian_conservation_laws}) as
\begin{align}
j^{1} V (\mathcal{F}) &= \bz_{\ct} + \psi_{\cx} \bz_{\cy} - \psi_{\cy} \bz_{\cx} = \mathcal{F} ,
\end{align}
hence $\tilde{\eta} = - \tilde{\by}$ and
\begin{align}
j^{1} \tilde{V}
&= \by \dfrac{\partial}{\partial \by} 
 + \bz_{\ct} \dfrac{\partial}{\partial \bz_{\ct}} 
 + \bz_{\cx} \dfrac{\partial}{\partial \bz_{\cx}} 
 + \bz_{\cy} \dfrac{\partial}{\partial \bz_{\cy}}
 - \tilde{\by} \dfrac{\partial}{\partial \tilde{\by}} 
 - \tilde{\bz}_{\ct} \dfrac{\partial}{\partial \tilde{\bz}_{\ct}} 
 - \tilde{\bz}_{\cx} \dfrac{\partial}{\partial \tilde{\bz}_{\cx}} 
 - \tilde{\bz}_{\cy} \dfrac{\partial}{\partial \tilde{\bz}_{\cy}} .
\end{align}
The prolongation of the extended vector field acts on the Lagrangian (\ref{eq:vorticity_lagrangian_conservation_laws}) as
\begin{align}
\jp^{1} \tilde{V} (\bar{L} ) 
= \tilde{\by} \, ( \bz_{\ct} + \psi_{\cx} \bz_{\cy} - \psi_{\cy} \bz_{\cx} ) 
- \tilde{\by} \, ( \bz_{\ct} + \psi_{\cx} \bz_{\cy} - \psi_{\cy} \bz_{\cx} ) = 0 ,
\end{align}
so that we have a symmetry whose Noether charge is computed as
\begin{align}
\mathcal{\tilde{J}}
= \int \tilde{J}^{t} \, d\cx
= \int \dfrac{\partial L}{\partial \bz_{t}} \big( \jp^{1} \tilde{\phy} \big) \cdot \eta (\tilde{\phy}) \, d\cx
= \int \zeta \omega \, d\cx
= \text{const. for all $\ct$} .
\end{align}
Restricting the Noether current $\tilde{J}$ with $\zeta = \phi (\omega) = \omega$, this becomes conservation of enstrophy,
\begin{align}
\mathcal{J}
= \int \big( \tilde{J}^{t} \circ \jp^{1} \Phi \big) (\jp^{1} \phy) \, d\cx
= \int \omega^{2} \, d\cx
= \text{const. for all $\ct$} .
\end{align}

In summary, we obtain conservation of
\begin{enumerate}[(a)]
\begin{subequations}\label{eq:vorticity_conservation_laws}
\item circulation (``mass'')
\begin{align}
\int \omega \, d\cx \, d\cy ,
\end{align}
\item enstrophy ($L^{2}$ norm)
\begin{align}
\int \omega^{2} \, d\cx \, d\cy ,
\end{align}
\item and kinetic energy
\begin{align}
\dfrac{1}{2} \int \psi \omega \, d\cx \, d\cy .
\end{align}
\end{subequations}
\end{enumerate}
For the symmetrised Lagrangian, c.f. equation (\ref{eq:vorticity_lagrangian_symmetrised}),
\begin{align}\label{eq:vorticity_conservation_laws_lagrangian_symmetrised}
\bar{L}' ( \jp^{1} \tilde{\phy} ) = \zeta \omega_{\ct} + \dfrac{1}{3} \Big[ \zeta \, \{ \psi , \omega \} + \psi \, \{ \omega , \zeta \} + \omega \, \{ \zeta , \psi \} \Big] ,
\end{align}
we obtain the same conservation laws, although in that case they correspond to divergence symmetries.

\subsubsection{Variational Integrator}

\begin{figure}[t]
	\centering
\begin{minipage}[c]{.06\textwidth}
$\by_{\ct}$:
\end{minipage}
\begin{minipage}[c]{.22\textwidth}
\includegraphics[width=\textwidth]{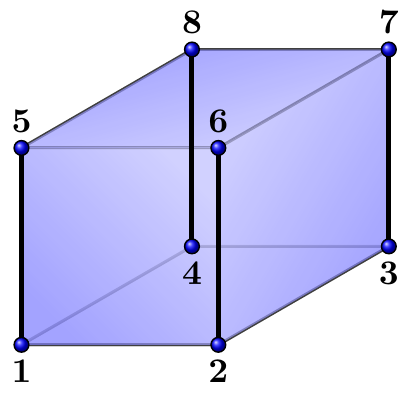}
\end{minipage}
\begin{minipage}[c]{.06\textwidth}
$\by_{\cx}$:
\end{minipage}
\begin{minipage}[c]{.22\textwidth}
\includegraphics[width=\textwidth]{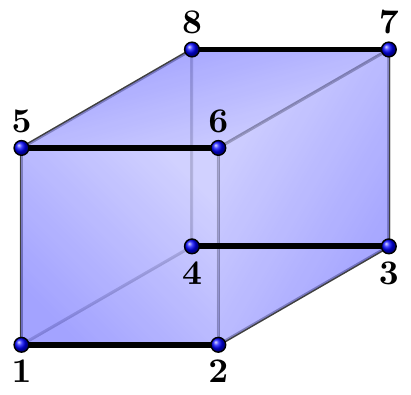}
\end{minipage}
\begin{minipage}[c]{.06\textwidth}
$\by_{\cy}$:
\end{minipage}
\begin{minipage}[c]{.22\textwidth}
\includegraphics[width=\textwidth]{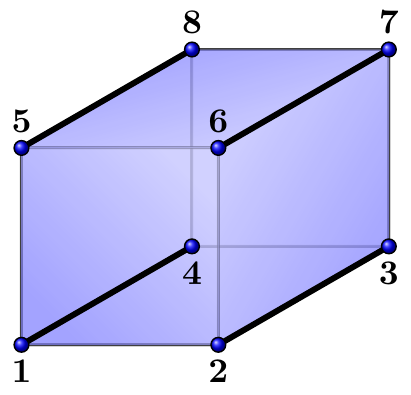}
\end{minipage}
\begin{minipage}[c]{.10\textwidth}
\includegraphics[width=\textwidth]{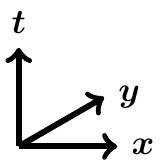}
\end{minipage}

	\caption{For a given spacetime grid cell, there are four possible ways of defining derivatives in the different coordinate directions $(\ct, \cx, \cy)$, highlighted along the black lines.}
	\label{fig:vorticity_derivatives}
\end{figure}

To discretise the Lagrangian $L (\tilde{\phy}, \tilde{\phy}_{\ct}, \tilde{\phy}_{\cx}, \tilde{\phy}_{\cy})$ from (\ref{eq:vorticity_lagrangian_symmetrised}) we adopt the same strategy as in section \ref{sec:linear_advection_simplified_implicit}.
We use simple finite differences to approximate the derivatives (see Figure~\ref{fig:vorticity_derivatives}), the midpoint rule for the time integral and the trapezoidal rule for the two spatial integrals,
\begingroup
\allowdisplaybreaks
\begin{align}\label{eq:vorticity_discrete_lagrangian}
L_{d} \big( \jp^{1} \tilde{\phy} (\square) \big)
\nonumber
&= \dfrac{h_{t} h_{x} h_{y}}{4} L' \Big(
  \tfrac{\tilde{\phy}_{\square^{1}} + \tilde{\phy}_{\square^{5}}}{2} ,
  \tfrac{\tilde{\phy}_{\square^{5}} - \tilde{\phy}_{\square^{1}}}{h_{t}} ,
  \tfrac{1}{2} \Big( \tfrac{\tilde{\phy}_{\square^{2}} - \tilde{\phy}_{\square^{1}}}{h_{x}}
                   + \tfrac{\tilde{\phy}_{\square^{6}} - \tilde{\phy}_{\square^{5}}}{h_{x}} \Big) , 
  \tfrac{1}{2} \Big( \tfrac{\tilde{\phy}_{\square^{4}} - \tilde{\phy}_{\square^{1}}}{h_{y}}
                   + \tfrac{\tilde{\phy}_{\square^{8}} - \tilde{\phy}_{\square^{5}}}{h_{y}} \Big)
\Big) \\
\nonumber
&+ \dfrac{h_{t} h_{x} h_{y}}{4} L' \Big(
  \tfrac{\tilde{\phy}_{\square^{2}} + \tilde{\phy}_{\square^{6}}}{2} ,
  \tfrac{\tilde{\phy}_{\square^{6}} - \tilde{\phy}_{\square^{2}}}{h_{t}} ,
  \tfrac{1}{2} \Big( \tfrac{\tilde{\phy}_{\square^{2}} - \tilde{\phy}_{\square^{1}}}{h_{x}}
                   + \tfrac{\tilde{\phy}_{\square^{6}} - \tilde{\phy}_{\square^{5}}}{h_{x}} \Big) ,
  \tfrac{1}{2} \Big( \tfrac{\tilde{\phy}_{\square^{3}} - \tilde{\phy}_{\square^{2}}}{h_{y}}
                   + \tfrac{\tilde{\phy}_{\square^{7}} - \tilde{\phy}_{\square^{6}}}{h_{y}} \Big)
\Big) \\
\nonumber
&+ \dfrac{h_{t} h_{x} h_{y}}{4} L' \Big(
  \tfrac{\tilde{\phy}_{\square^{3}} + \tilde{\phy}_{\square^{7}}}{2} ,
  \tfrac{\tilde{\phy}_{\square^{7}} - \tilde{\phy}_{\square^{3}}}{h_{t}} ,
  \tfrac{1}{2} \Big( \tfrac{\tilde{\phy}_{\square^{3}} - \tilde{\phy}_{\square^{4}}}{h_{x}}
                   + \tfrac{\tilde{\phy}_{\square^{7}} - \tilde{\phy}_{\square^{8}}}{h_{x}} \Big) ,
  \tfrac{1}{2} \Big( \tfrac{\tilde{\phy}_{\square^{3}} - \tilde{\phy}_{\square^{2}}}{h_{y}}
                   + \tfrac{\tilde{\phy}_{\square^{7}} - \tilde{\phy}_{\square^{6}}}{h_{y}} \Big)
\Big) \\
&+ \dfrac{h_{t} h_{x} h_{y}}{4} L' \Big(
  \tfrac{\tilde{\phy}_{\square^{4}} + \tilde{\phy}_{\square^{8}}}{2} ,
  \tfrac{\tilde{\phy}_{\square^{8}} - \tilde{\phy}_{\square^{4}}}{h_{t}} ,
  \tfrac{1}{2} \Big( \tfrac{\tilde{\phy}_{\square^{3}} - \tilde{\phy}_{\square^{4}}}{h_{x}}
                   + \tfrac{\tilde{\phy}_{\square^{7}} - \tilde{\phy}_{\square^{8}}}{h_{x}} \Big) ,
  \tfrac{1}{2} \Big( \tfrac{\tilde{\phy}_{\square^{4}} - \tilde{\phy}_{\square^{1}}}{h_{y}}
                   + \tfrac{\tilde{\phy}_{\square^{8}} - \tilde{\phy}_{\square^{5}}}{h_{y}} \Big)
\Big) .
\end{align}
\endgroup
In $2+1$ dimensions, the discrete Euler-Lagrange field equations have eight contributions instead of four,
\begin{align}
0
\nonumber
&= \dfrac{\partial L_d}{\partial \tilde{\phy}^{a}_{\square^1}} \Big( \tilde{\phy}_{i,  j,  k  }, \tilde{\phy}_{i,  j+1,k  }, \tilde{\phy}_{i,  j+1,k+1}, \tilde{\phy}_{i,  j,  k+1}, \tilde{\phy}_{i+1,j,  k  }, \tilde{\phy}_{i+1,j+1,k  }, \tilde{\phy}_{i+1,j+1,k+1}, \tilde{\phy}_{i+1,j,  k+1} \Big) \\
\nonumber
&+ \dfrac{\partial L_d}{\partial \tilde{\phy}^{a}_{\square^2}} \Big( \tilde{\phy}_{i,  j-1,k  }, \tilde{\phy}_{i,  j,  k  }, \tilde{\phy}_{i,  j,  k+1}, \tilde{\phy}_{i,  j-1,k+1}, \tilde{\phy}_{i+1,j-1,k  }, \tilde{\phy}_{i+1,j,  k  }, \tilde{\phy}_{i+1,j,  k+1}, \tilde{\phy}_{i+1,j-1,k+1} \Big) \\
\nonumber
&+ \dfrac{\partial L_d}{\partial \tilde{\phy}^{a}_{\square^3}} \Big( \tilde{\phy}_{i,  j-1,k-1}, \tilde{\phy}_{i,  j,  k-1}, \tilde{\phy}_{i,  j,  k  }, \tilde{\phy}_{i,  j-1,k  }, \tilde{\phy}_{i+1,j-1,k-1}, \tilde{\phy}_{i+1,j,  k-1}, \tilde{\phy}_{i+1,j,  k  }, \tilde{\phy}_{i+1,j-1,k  } \Big) \\
\nonumber
&+ \dfrac{\partial L_d}{\partial \tilde{\phy}^{a}_{\square^4}} \Big( \tilde{\phy}_{i,  j  ,k-1}, \tilde{\phy}_{i,  j+1,k-1}, \tilde{\phy}_{i,  j+1,k  }, \tilde{\phy}_{i,  j,  k  }, \tilde{\phy}_{i+1,j,  k-1}, \tilde{\phy}_{i+1,j+1,k-1}, \tilde{\phy}_{i+1,j+1,k  }, \tilde{\phy}_{i+1,j,  k  } \Big) \\
\nonumber
&+ \dfrac{\partial L_d}{\partial \tilde{\phy}^{a}_{\square^5}} \Big( \tilde{\phy}_{i-1,j,  k  }, \tilde{\phy}_{i-1,j+1,k  }, \tilde{\phy}_{i-1,j+1,k+1}, \tilde{\phy}_{i-1,j,  k+1}, \tilde{\phy}_{i,  j,  k  }, \tilde{\phy}_{i,  j+1,k  }, \tilde{\phy}_{i,  j+1,k+1}, \tilde{\phy}_{i,  j,  k+1} \Big) \\
\nonumber
&+ \dfrac{\partial L_d}{\partial \tilde{\phy}^{a}_{\square^6}} \Big( \tilde{\phy}_{i-1,j-1,k  }, \tilde{\phy}_{i-1,j,  k  }, \tilde{\phy}_{i-1,j,  k+1}, \tilde{\phy}_{i-1,j-1,k+1}, \tilde{\phy}_{i,  j-1,k  }, \tilde{\phy}_{i,  j,  k  }, \tilde{\phy}_{i,  j,  k+1}, \tilde{\phy}_{i,  j-1,k+1} \Big) \\
\nonumber
&+ \dfrac{\partial L_d}{\partial \tilde{\phy}^{a}_{\square^7}} \Big( \tilde{\phy}_{i-1,j-1,k-1}, \tilde{\phy}_{i-1,j,  k-1}, \tilde{\phy}_{i-1,j,  k  }, \tilde{\phy}_{i-1,j-1,k  }, \tilde{\phy}_{i,  j-1,k-1}, \tilde{\phy}_{i,  j,  k-1}, \tilde{\phy}_{i,  j,  k  }, \tilde{\phy}_{i,  j-1,k  } \Big) \\
&+ \dfrac{\partial L_d}{\partial \tilde{\phy}^{a}_{\square^8}} \Big( \tilde{\phy}_{i-1,j  ,k-1}, \tilde{\phy}_{i-1,j+1,k-1}, \tilde{\phy}_{i-1,j+1,k  }, \tilde{\phy}_{i-1,j,  k  }, \tilde{\phy}_{i,  j,  k-1}, \tilde{\phy}_{i,  j+1,k-1}, \tilde{\phy}_{i,  j+1,k  }, \tilde{\phy}_{i,  j,  k  } \Big) .
\end{align}
The resulting integrator for the vorticity equation takes the form
\begin{multline}\label{eq:vorticity_integrator}
\dfrac{\omega_{i+1} - \omega_{i-1}}{2 h_{t}}
+ \dfrac{1}{8} \Big[ A (\psi_{i+1}, \omega_{i+1}) + A (\psi_{i}, \omega_{i+1}) + A (\psi_{i+1}, \omega_{i}) \\
+ 2 A(\psi_{i}, \omega_{i}) + A(\psi_{i-1}, \omega_{i}) + A (\psi_{i}, \omega_{i-1}) + A (\psi_{i-1}, \omega_{i-1}) \Big] ,
\end{multline}
where $\omega_{i} = (\omega_{ij})_{j}$ and $\psi_{i} = (\psi_{ij})_{j}$ are the rows of the solution matrix corresponding to constant-time slices, whereas $A$ denotes Arakawa's discretisation of the Poisson bracket \cite{Arakawa:1966}, given by
\begin{align}\label{eq:brackets_poisson_arakawa_5}
A (\psi, \omega) = \dfrac{1}{3} \Big( A^{++} (\psi, \omega) + A^{+ \times} (\psi, \omega) + A^{\times +} (\psi, \omega) \Big) ,
\end{align}
with $\omega$ and $\psi$ being row vectors and
\begingroup
\allowdisplaybreaks
\begin{subequations}\label{eq:brackets_poisson_arakawa_4}
\begin{align}
A^{++} (\psi, \omega)
\nonumber
= \dfrac{1}{4 h_{x} h_{y}} \Big[
     & \big( \psi_{+0} - \psi_{-0} \big) \big( \omega_{0+} - \omega_{0-} \big) \\
- \, & \big( \psi_{0+} - \psi_{0-} \big) \big( \omega_{+0} - \omega_{-0} \big) \Big] , \\
A^{+ \times} (\psi, \omega)
\nonumber
= \dfrac{1}{4 h_{x} h_{y}} \Big[
     & \psi_{+0} \big( \omega_{+-} - \omega_{++} \big) - \psi_{-0} \big( \omega_{--} - \omega_{-+} \big) \\
- \, & \psi_{0+} \big( \omega_{-+} - \omega_{++} \big) + \psi_{0-} \big( \omega_{--} - \omega_{+-} \big) \Big] , \\
A^{\times +} (\psi, \omega)
\nonumber
= \dfrac{1}{4 h_{x} h_{y}} \Big[
     & \psi_{++} \big( \omega_{+0} - \omega_{0+} \big) - \psi_{--} \big( \omega_{0-} - \omega_{-0} \big) \\
- \, & \psi_{-+} \big( \omega_{-0} - \omega_{0+} \big) + \psi_{+-} \big( \omega_{0-} - \omega_{+0} \big) \Big] .
\end{align}
\end{subequations}
\endgroup
The indices $(-,0,+)$ indicate the increment of the corresponding index on the grid relative to the point where the bracket is computed.
The integrator for the Poisson equation is
\begin{align}\label{eq:vorticity_integrator_poisson}
& \bracket{\Delta_{x} \psi_{i,j,k}}_{t}
+ \bracket{\Delta_{y} \psi_{i,j,k}}_{t}
= \bracket{\omega_{i,j,k}}_{t} .
\end{align}
Here, $\Delta_{x}$ and $\Delta_{y}$ denote the standard centred finite difference second-order derivative with respect to $x$ and $y$, i.e.,
\begin{align*}
\Delta_{x} \psi_{i,j,k} &= \dfrac{\psi_{i,j-1,k} - 2 \psi_{i,j,k} + \psi_{i,j+1,k}}{h_{x}^{2}} , \\
\Delta_{y} \psi_{i,j,k} &= \dfrac{\psi_{i,j,k-1} - 2 \psi_{i,j,k} + \psi_{i,j,k+1}}{h_{y}^{2}} ,
\end{align*}
and the angle brackets $\bracket{\cdot}_{t}$ denote an average in time of the form $\tfrac{1}{4} \big[ 1 \; 2 \; 1 \big]$, namely,
\begin{align*}
\bracket{\omega_{i,j,k}}_{t} = \dfrac{\omega_{i-1,j,k} + 2 \omega_{i,j,k} + \omega_{i+1,j,k}}{4} .
\end{align*}

It is remarkable that Arakawa's discretisation of the Poisson bracket is recovered. Indeed, a similar derivation was proposed by \citeauthor{SalmonTalley:1989} \cite{SalmonTalley:1989}.
Our approach, however, is fully covariant, leading to a complete spacetime discretisation, that is a combination of Arakawa's scheme in space with a symplectic integrator in time.
Whereas the Arakawa bracket alone guarantees energy conservation only for the spatial discretisation, i.e., up to errors due to the discretisation of the time derivative, the variational integrator is exactly energy preserving.

\subsubsection{Simplifications}\label{sec:vorticity_simplifications}

The integrator (\ref{eq:vorticity_integrator}) for the vorticity equation has one drawback. Even though we are considering a partial differential equation that is first order in time, we need to initialise the fields at two consecutive points in the discrete time domain, that is (\ref{eq:vorticity_integrator}) is a multistep integrator.
Nevertheless, we can reduce it to a single-step integrator.
Observe that we can rewrite (\ref{eq:vorticity_integrator}) in the following way,
\begin{align}\label{eq:vorticity_integrator_rewritten}
\nonumber
0
&= \dfrac{\psi_{i+1} - \psi_{i}}{2 h_{t}}
 + \dfrac{1}{8} \Big[ A (\psi_{i+1}, \omega_{i+1}) + A (\psi_{i+1}, \omega_{i}) + A (\psi_{i}, \omega_{i+1}) + A(\psi_{i}, \omega_{i}) \Big] \\
&+ \dfrac{\psi_{i} - \psi_{i-1}}{2 h_{t}}
 + \dfrac{1}{8} \Big[ A(\psi_{i}, \omega_{i}) + A(\psi_{i}, \omega_{i-1}) + A (\psi_{i-1}, \omega_{i}) + A (\psi_{i-1}, \omega_{i-1}) \Big] ,
\end{align}
which is the average over two points in time of the single-step integrator
\begin{align}\label{eq:vorticity_integrator_simplified}
0
&= \dfrac{\psi_{i} - \psi_{i-1}}{h_{t}}
 + \dfrac{1}{4} \Big[ A(\psi_{i}, \omega_{i}) + A(\psi_{i}, \omega_{i-1}) + A (\psi_{i-1}, \omega_{i}) + A (\psi_{i-1}, \omega_{i-1}) \Big] .
\end{align}
Analogously, the Poisson equation (\ref{eq:vorticity_integrator_poisson}) is the average over three points in time of
\begin{align}\label{eq:vorticity_integrator_poisson_simplified}
& \Delta_{x} \psi_{i,j,k} + \Delta_{y} \psi_{i,j,k} = \omega_{i,j,k} .
\end{align}
The solution of the single-step integrator (\ref{eq:vorticity_integrator_simplified}-\ref{eq:vorticity_integrator_poisson_simplified}) is fully determined by specifying the vorticity $\omega = \omega_{0}$ at $t=t_{0}$ on the spatial grid.
Then $\psi_{0}$ can be obtained by (\ref{eq:vorticity_integrator_poisson_simplified}) and the vorticity can be advected by solving (\ref{eq:vorticity_integrator_simplified}).
We observe that if the sequence $(\omega_{i}, \psi_{i})$ solves (\ref{eq:vorticity_integrator_simplified}-\ref{eq:vorticity_integrator_poisson_simplified}) with initial conditions $\omega_{0}$, then it is also a solution of (\ref{eq:vorticity_integrator}-\ref{eq:vorticity_integrator_poisson}) with initialisation $(\omega_{0}, \psi_{0}, \omega_{1}, \psi_{1})$.
Vice versa, if the sequence $(\omega_{i}, \psi_{i})$ is a solution of (\ref{eq:vorticity_integrator_simplified}-\ref{eq:vorticity_integrator_poisson_simplified}) initialised by using (\ref{eq:vorticity_integrator_simplified}-\ref{eq:vorticity_integrator_poisson_simplified}) with the initial data $\omega_{0}$, then it is also a solution of (\ref{eq:vorticity_integrator_simplified}-\ref{eq:vorticity_integrator_poisson_simplified}).
Equations (\ref{eq:vorticity_integrator_simplified}-\ref{eq:vorticity_integrator_poisson_simplified}) are called the underlying one-step method of (\ref{eq:vorticity_integrator}-\ref{eq:vorticity_integrator_poisson}).

The simplified equations (\ref{eq:vorticity_integrator_simplified}) and (\ref{eq:vorticity_integrator_poisson_simplified}) can also be obtained directly from a discrete Lagrangian as is shown in \cite{KrausMajSonnendruecker:2015}.
This is of course preferable as it allows us to check for symmetries and to compute the discrete conserved quantities more easily.
As solutions of the simplified integrator (\ref{eq:vorticity_integrator_simplified}-\ref{eq:vorticity_integrator_poisson_simplified}) will also be solutions of (\ref{eq:vorticity_integrator}-\ref{eq:vorticity_integrator_poisson}), they will satisfy the same conservation laws.
Nevertheless, it is possible that (\ref{eq:vorticity_integrator_simplified}-\ref{eq:vorticity_integrator_poisson_simplified}) will have additional symmetries and therefore additional invariants.

\subsubsection{Discrete Conservation Laws}

The discrete conservation laws of the variational integrator (\ref{eq:vorticity_integrator}) for the vorticity equation are computed in the same way as for the linear advection equation (section \ref{sec:advection}) except that we have to account for discrete divergence symmetries as introduced in section \ref{sec:vi_discrete_noether_theorem} and that we are using the generators discussed in section \ref{sec:vorticity_noether}.
As in the continuous case, we restrict our analysis to the linear case, that is the Lagrangian (\ref{eq:vorticity_conservation_laws_lagrangian_symmetrised}) discretised according to (\ref{eq:vorticity_discrete_lagrangian}).

Because the vorticity equation is defined on three-dimensional spacetime, the actual computations are tediously long and hardly possible to carry out by hand.
We therefore provide only the results, obtained via computer aided calculations, namely
\begin{enumerate}[(a)]
\begin{subequations}\label{eq:vorticity_discrete_conservation_laws}
\item circulation (``mass'')
\begin{align}
h_x h_y \sum \limits_{j,k} \omega_{j,k}  = \mathrm{const.} ,
\end{align}
\item enstrophy ($L^{2}$ norm)
\begin{align}
h_x h_y \sum \limits_{j,k} \omega_{j,k}^{2} = \mathrm{const.} ,
\end{align}
\item and kinetic energy
\begin{align}
\dfrac{h_x h_y}{2} \sum \limits_{j,k} \omega_{j,k} \phi_{j,k} = \mathrm{const.} .
\end{align}
\end{subequations}
\end{enumerate}
It is well known that these quantities are preserved by Arakawa's discrete bracket \cite{Arakawa:1966}. In numerical simulations we indeed find that our integrator conserves them to machine accuracy or at least the tolerance of the nonlinear solver, not only in the linear case (when $\psi$ is constant in time), but even in the nonlinear case (with self-consistent streaming function $\psi$).

\subsubsection{Numerical Examples}

We implemented the simplified variational integrator (\ref{eq:vorticity_integrator_simplified}-\ref{eq:vorticity_integrator_poisson_simplified}) using Python \cite{Python, Langtangen:2014}, Cython \cite{Cython, Behnel:2010}, PETSc \cite{petsc-web-page, petsc-user-ref} and petsc4py \cite{Dalcin:2011}. Visualisation was done using NumPy \cite{vanDerWalt:2011}, SciPy \cite{SciPy} and matplotlib \cite{Matplotlib, Hunter:2007}.
The nonlinear system is solved via Picard iteration. Within each nonlinear step the two linear systems corresponding to the vorticity equation and the Poisson equation are solved separately. 
The vorticity equation is solved with GMRES and the Poisson equation via LU decomposition with SuperLU \cite{superlu, Li:2005}. The tolerance of the nonlinear solver is set to $10^{-10}$ or smaller, which is reached after $5-10$ iterations. Usually, the GMRES solver needs between $5$ and $25$ iterations to converge with a relative tolerance of $10^{-8}$ or an absolute tolerance of $10^{-15}$.

\subsubsection*{The Linear Case}

At first we consider the linear case, where we prescribe the streaming function
$\psi$ and only solve the vorticity equation. The streaming potential is set to
\begin{align}
\psi (x,y) = \tfrac{1}{2} y^2 + 1 - \cos (x) ,
\end{align}
while the vorticity is initialised with a localised, symmetric Gaussian, 
\begin{align}
\omega_{0} (x,y) = f(x, x_{0}, \sigma_{x}) \, f(y, y_{0}, \sigma_{y}) ,
\end{align}
with
\begin{align}
f(z, z_{0}, \sigma) = \dfrac{1}{\sigma \, \sqrt{2 \pi}} \, \exp \bigg( - \dfrac{1}{2} \bigg[ \dfrac{z-z_{0}}{\sigma} \bigg]^2 \bigg) ,
\hspace{3em}
z \in \mathbb{R} .
\end{align}
The domain is $\Omega = [-2\pi, +2\pi) \times [-2\pi, +2\pi)$ with periodic boundaries, the spatial resolution is $1024 \times 1024$ and the timestep is $10^{-2}$.

\begin{figure}[p]
	\centering
	\subfloat[][$t=0.0$]{
	\includegraphics[width=.42\textwidth]{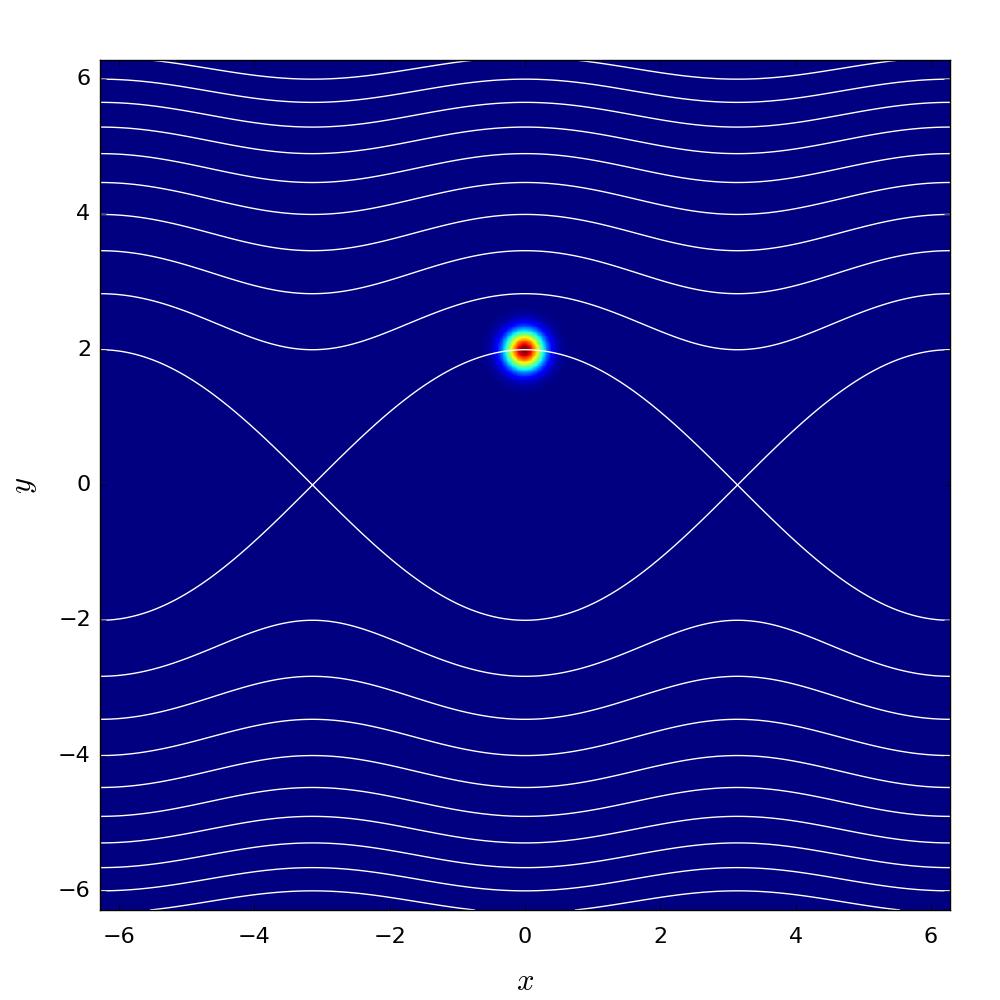}
	}
	\subfloat[][$t=0.5$]{
	\includegraphics[width=.42\textwidth]{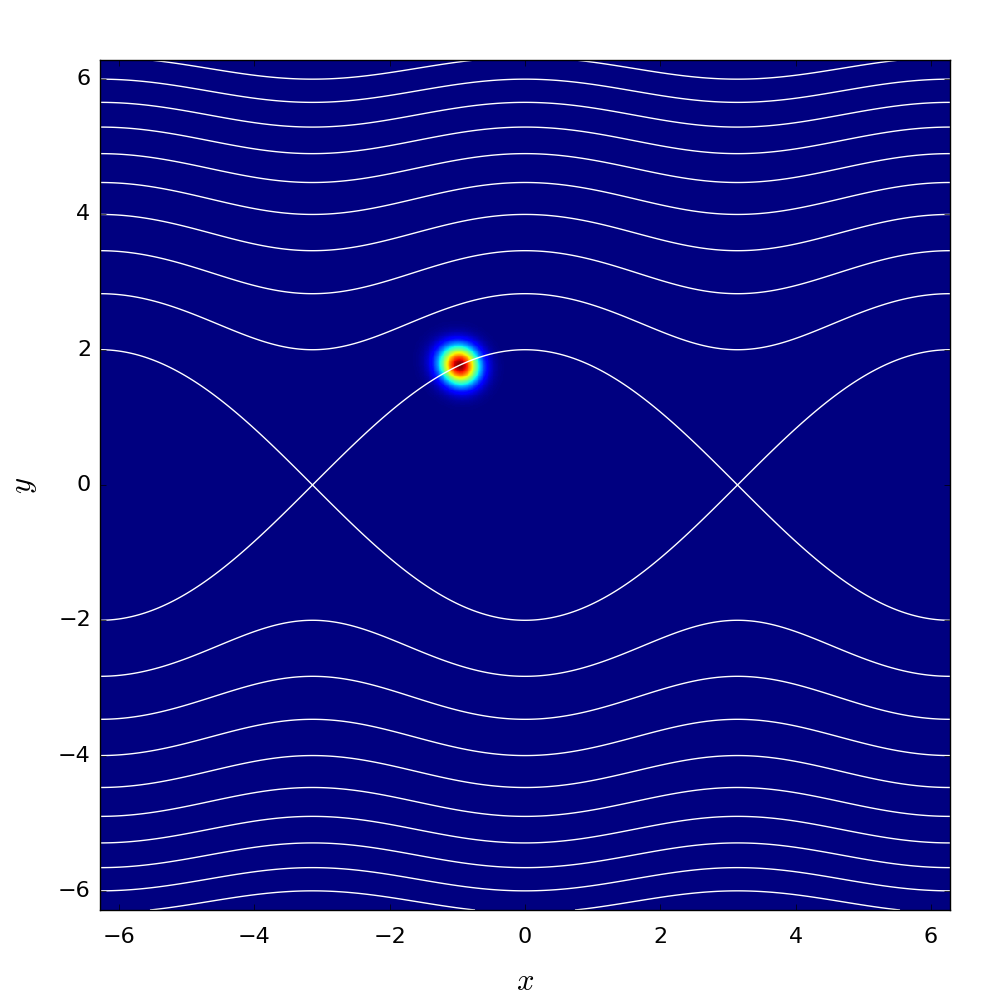}
	}
	
	\subfloat[][$t=1.0$]{
	\includegraphics[width=.42\textwidth]{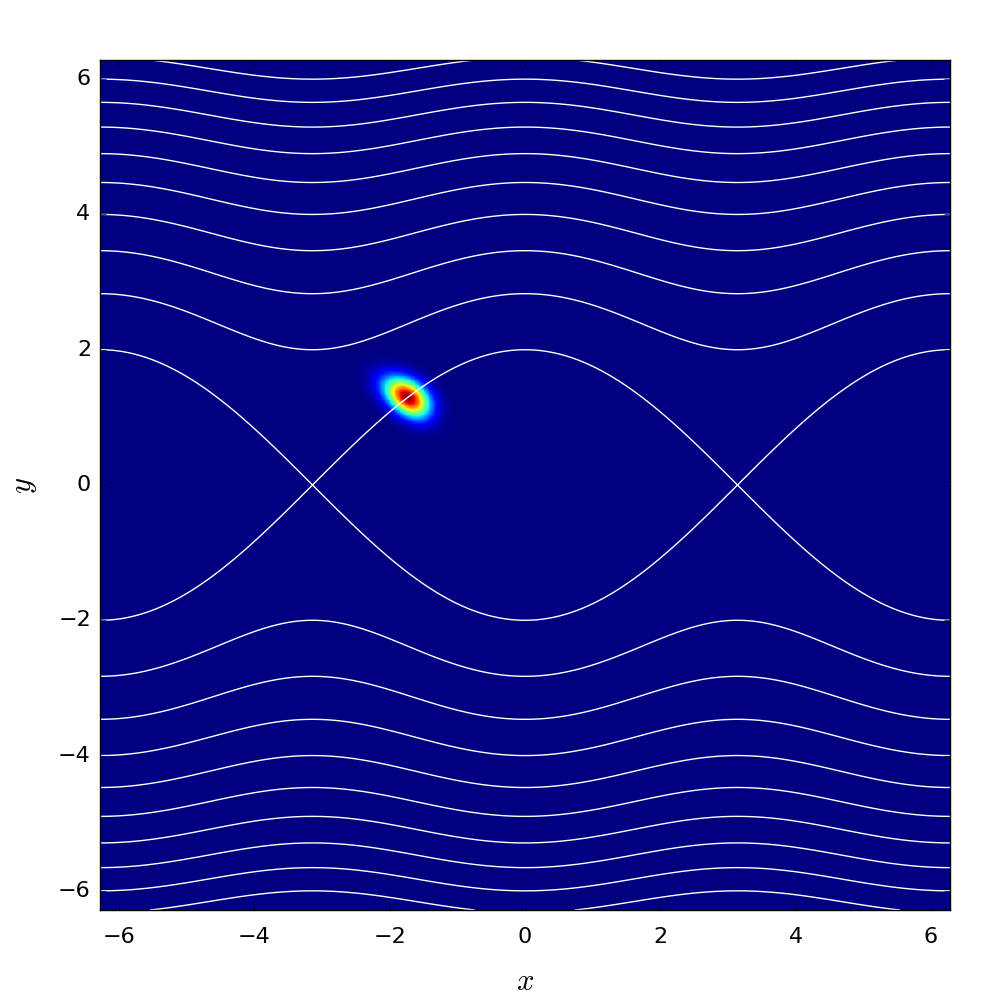}
	}
	\subfloat[][$t=2.0$]{
	\includegraphics[width=.42\textwidth]{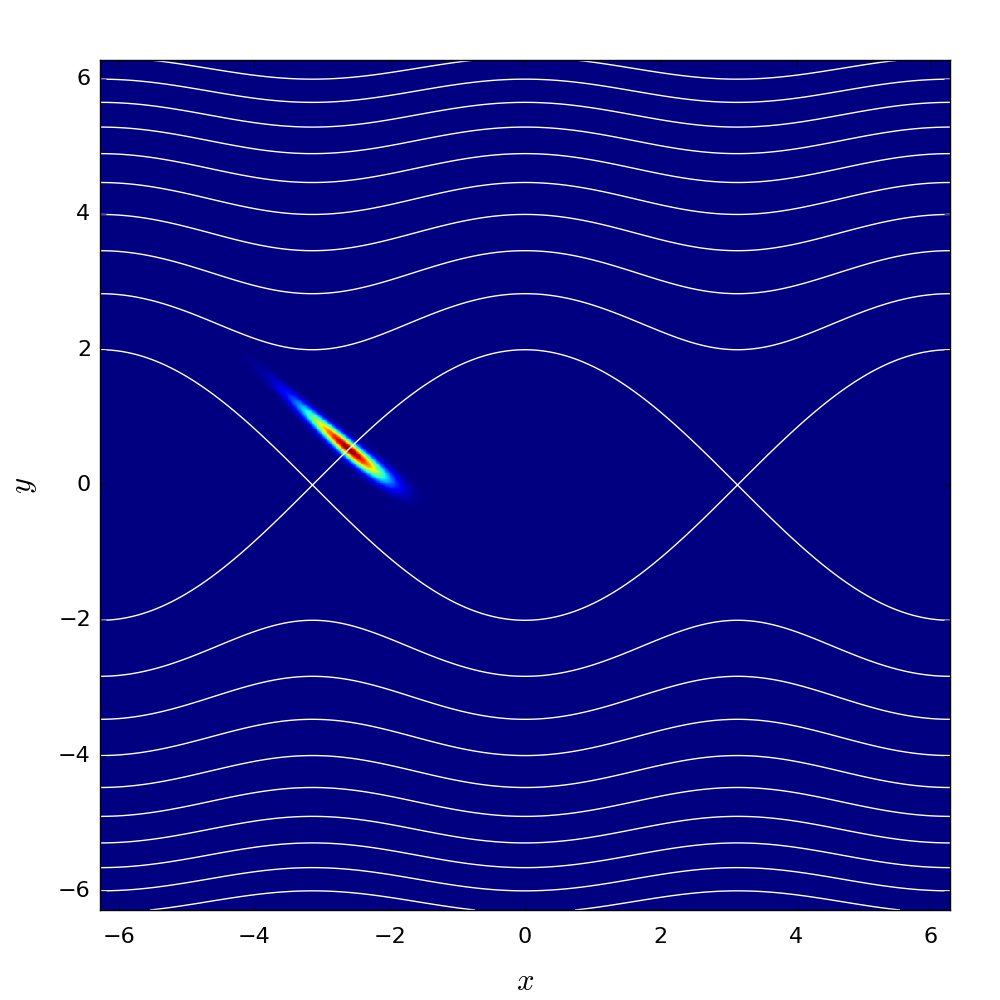}
	}
	
	\subfloat[][$t=3.0$]{
	\includegraphics[width=.42\textwidth]{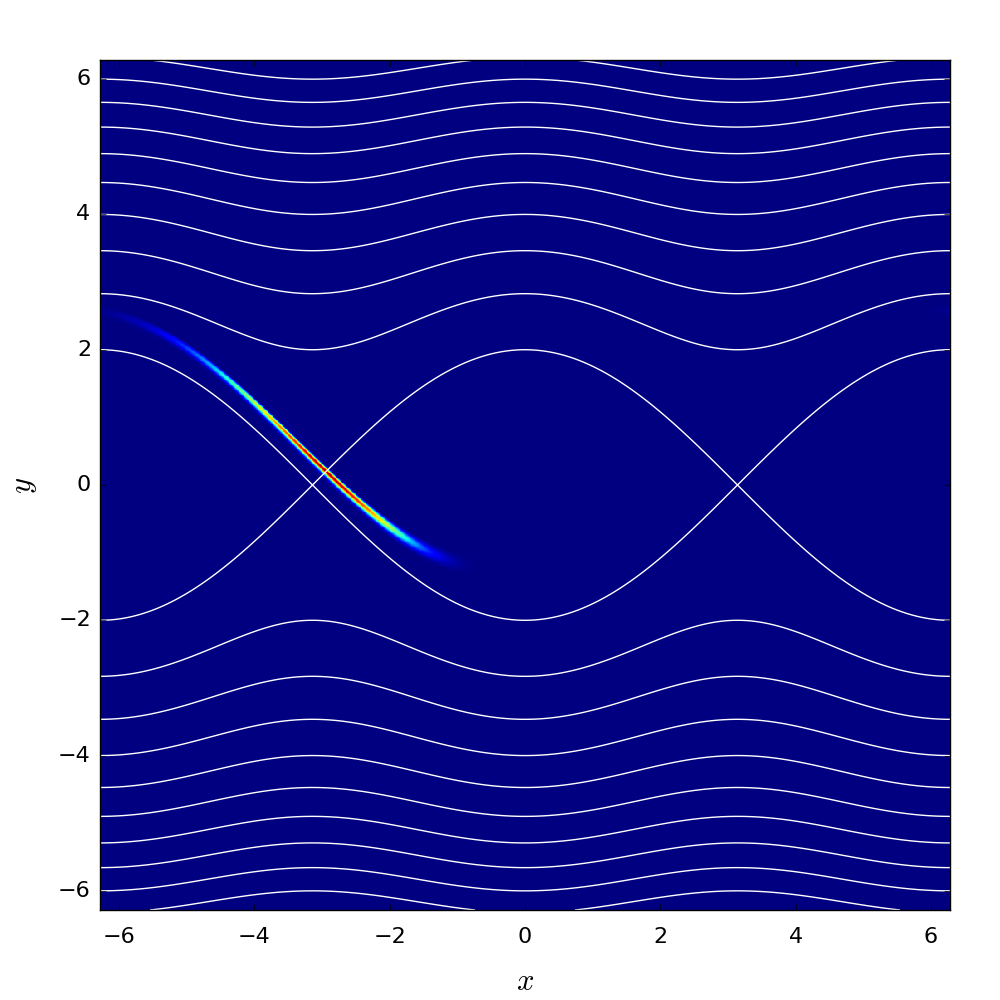}
	}
	\subfloat[][$t=6.0$]{
	\includegraphics[width=.42\textwidth]{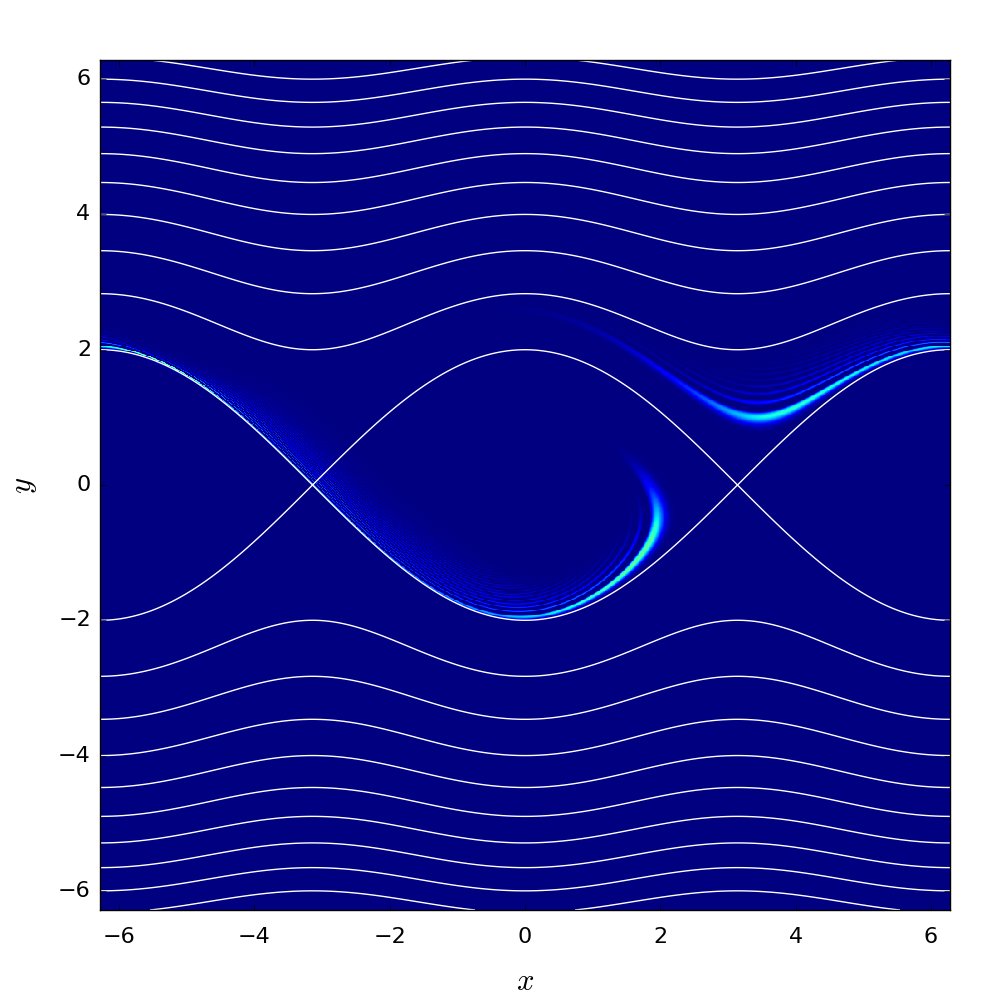}
	}
	
	\caption{Linear case with a separatrix. Vorticity $\omega$ at $t=0.0$, $0.5$, $1.0$, $2.0$, $3.0$ and $6.0$. Contours of the streaming function $\psi$ are also shown.}
	\label{fig:vorticity_linear_separatrix_vorticity}
\end{figure}

\begin{figure}[p]
	\centering
	\subfloat[][]{\label{fig:vorticity_linear_separatrix_contours}
	\label{fig:linear_separatrix_contours}
	\includegraphics[width=\textwidth]{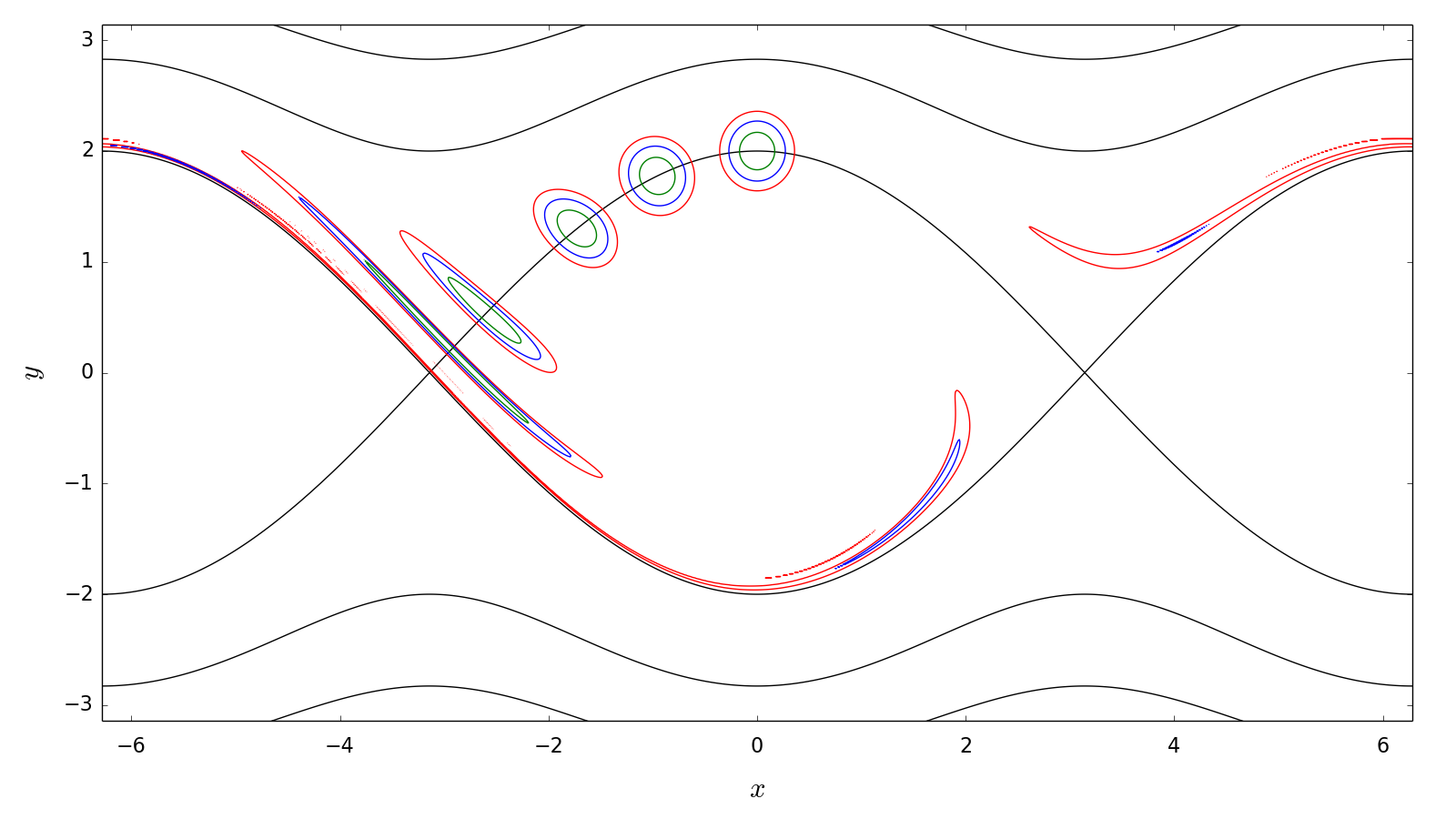}
	}

	\subfloat[][]{\label{fig:vorticity_linear_separatrix_traces}
	\label{fig:linear_separatrix_timetraces}
	\includegraphics[width=\textwidth]{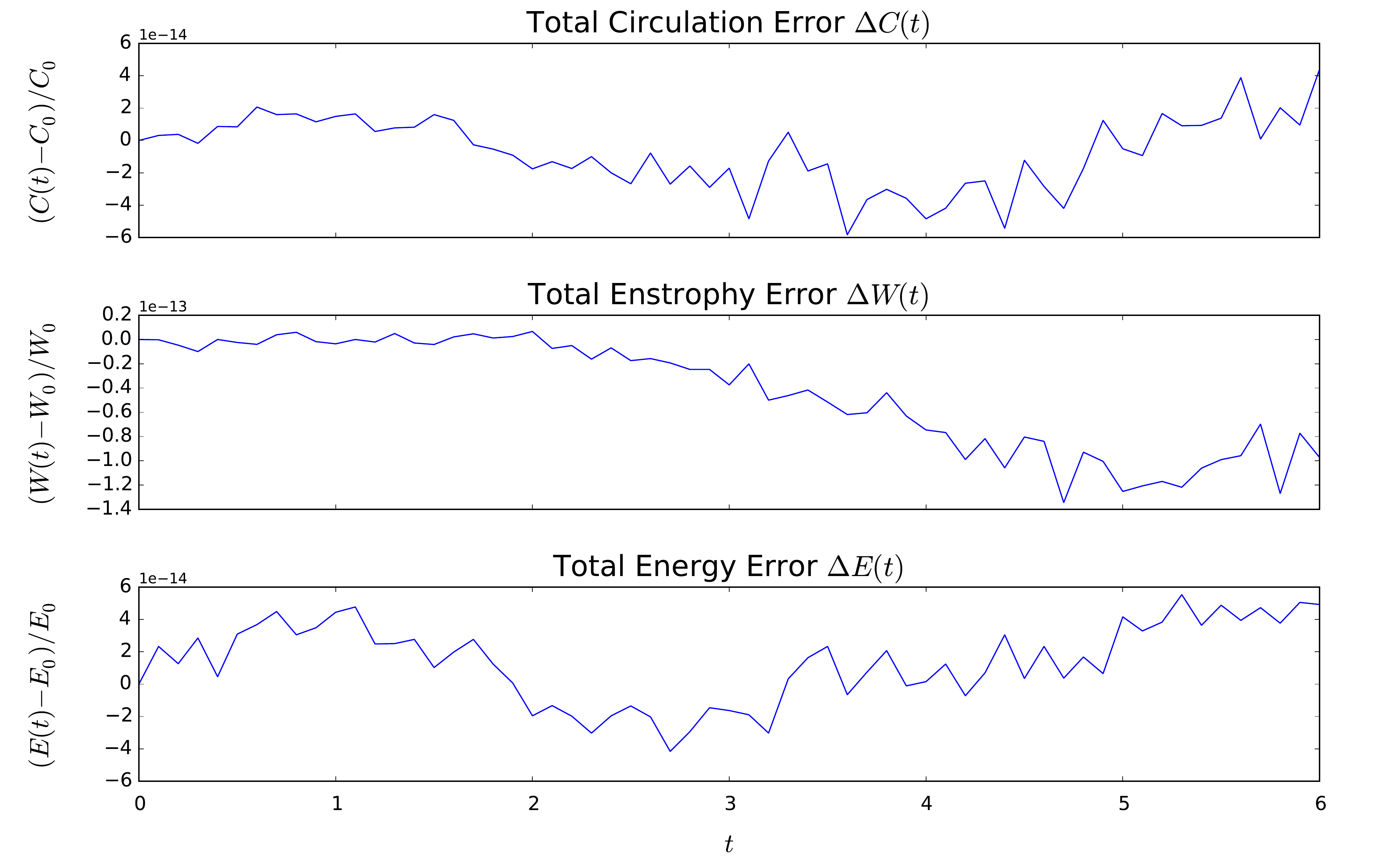}
	}

	\caption{Linear case with a separatrix. (a) Contours of the streaming function $\psi$ (black) and of the vorticity $\omega$ (at $0.2 \max \omega_{0}$, $0.4 \max \omega_{0}$ and $0.7 \max \omega_{0}$) at times $t=0.0$, $0.5$, $1.0$, $2.0$, $3.0$ and $6.0$. (b) Conservation laws.}
\end{figure}

The parameters are set to $x_{0} = 0$, $y_{0} = 2$ and $\sigma_{x} = \sigma_{y}
= 0.2$, so that the Gaussian is placed on the separatrix of the streaming
function (c.f., Figure~\ref{fig:vorticity_linear_separatrix_vorticity}). As the
Gaussian is moving along the separatrix, it is stretched while the contours of
the vorticity are preserved. In fact, if a field $\omega$ is advected by a smooth flow, the topology of the contours of $\omega$ should be preserved. This behaviour appears to be respected by the integrator. In Figure~\ref{fig:vorticity_linear_separatrix_contours}, the contours for different values of the vorticity are shown. Until $t=3$ all three contours stay intact. At about $t=6$ the stretching of the contours is so strong that resolution becomes insufficient and only the outermost contour is still preserved. Towards the centre and the boundaries in $x$ the effect of dispersion becomes visible.
Figure~\ref{fig:vorticity_linear_separatrix_contours} shows the evolution of the errors of circulation, enstrophy and energy.

\subsubsection*{Lamb Dipole}

The lamb dipole \cite{NielsenRasmussen:1997} is a stationary solution of the vorticity equation which is at rest in its frame of reference. 
The vorticity is initialised as
\begin{align}
\omega_{0} (x,y) = 
\begin{cases}
2 \lambda U \, \cos \theta \, \dfrac{J_{1} (\lambda r)}{J_{0} (\lambda R)} & r \leq R , \\
0 & r > R ,
\end{cases}
\end{align}
which leads to a dipole of radius $R$ moving in the $y$-direction with velocity $U$.
Here, $r$ and $\theta$ denote polar coordinates, i.e.,
\begin{align}
r = \sqrt{ x^2 + y^2 }
\hspace{3em}
\text{and}
\hspace{3em}
\theta = \arctan (y/x) ,
\end{align}
$J_{i}$ is the $i$-th Bessel function of first kind and $\lambda R$ is the first root of $J_{1}$. In the example, we use $R = 0.2$ and $U = 1$. The domain is $\Omega = [-1, +1) \times [-1, +1)$ with periodic boundaries, the spatial resolution is $1024 \times 1024$ and the timestep is $10^{-3}$. The tolerance of the nonlinear solver is set to $10^{-12}$.

\begin{figure}[p]
	\centering
	\subfloat[][$t=0.0$]{\label{fig:vorticity_lamb_dipole_vorticity_initial}
	\includegraphics[width=.45\textwidth]{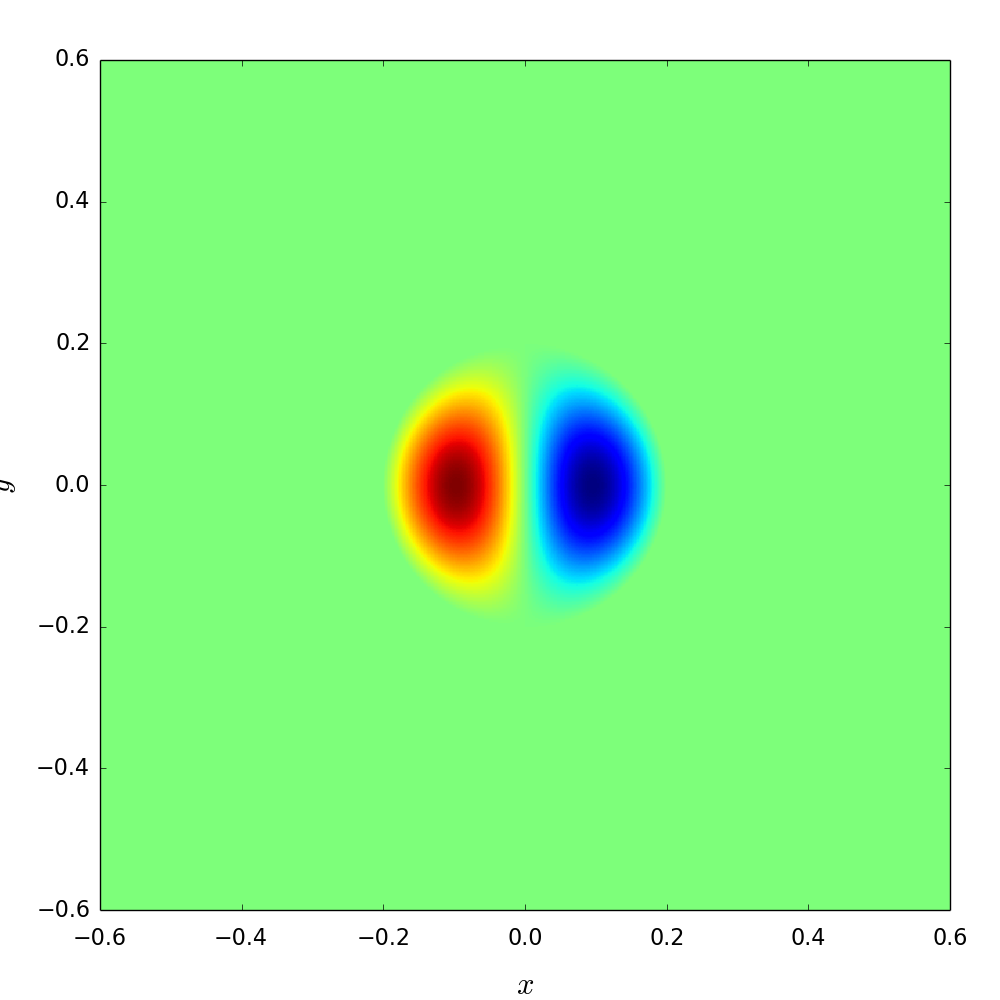}
	}
	\subfloat[][$t=10.0$]{\label{fig:vorticity_lamb_dipole_vorticity_final}
	\includegraphics[width=.45\textwidth]{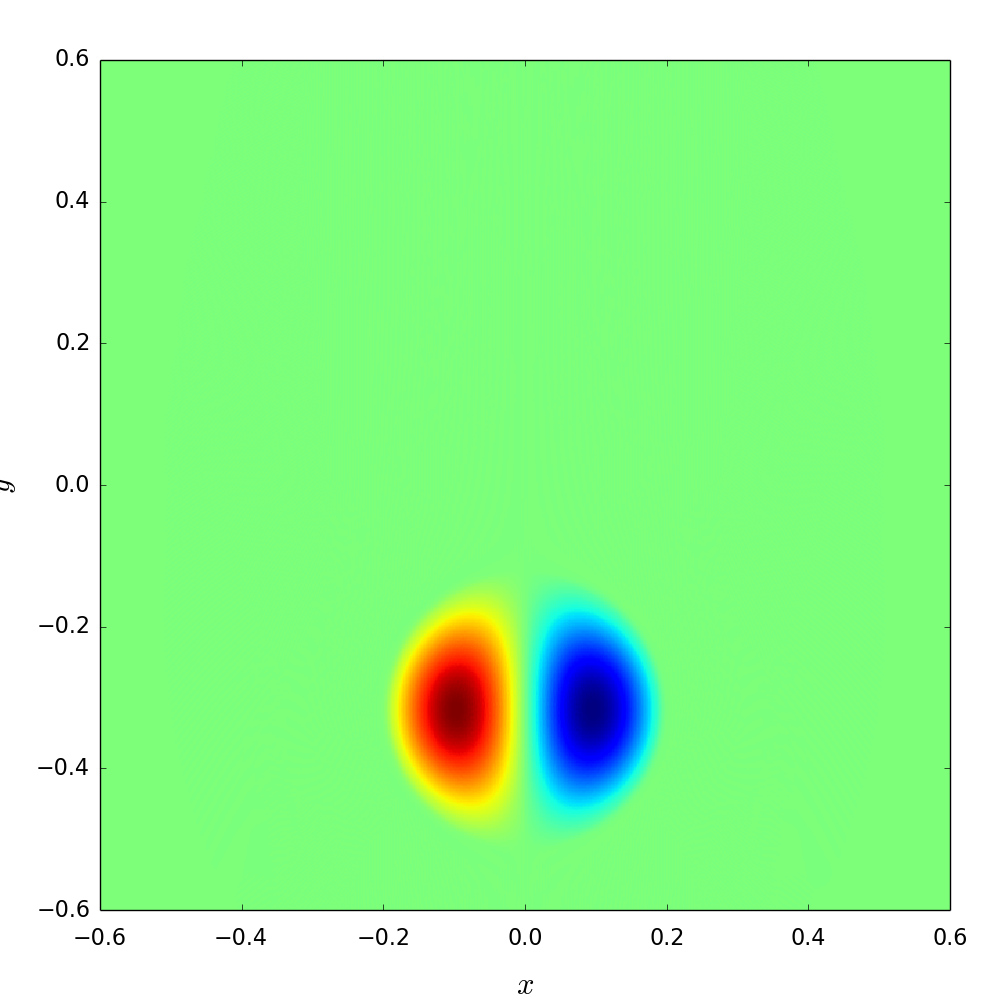}
	}
	
	\subfloat[][]{
	\label{fig:vorticity_lamb_dipole_timetraces}
	\includegraphics[width=\textwidth]{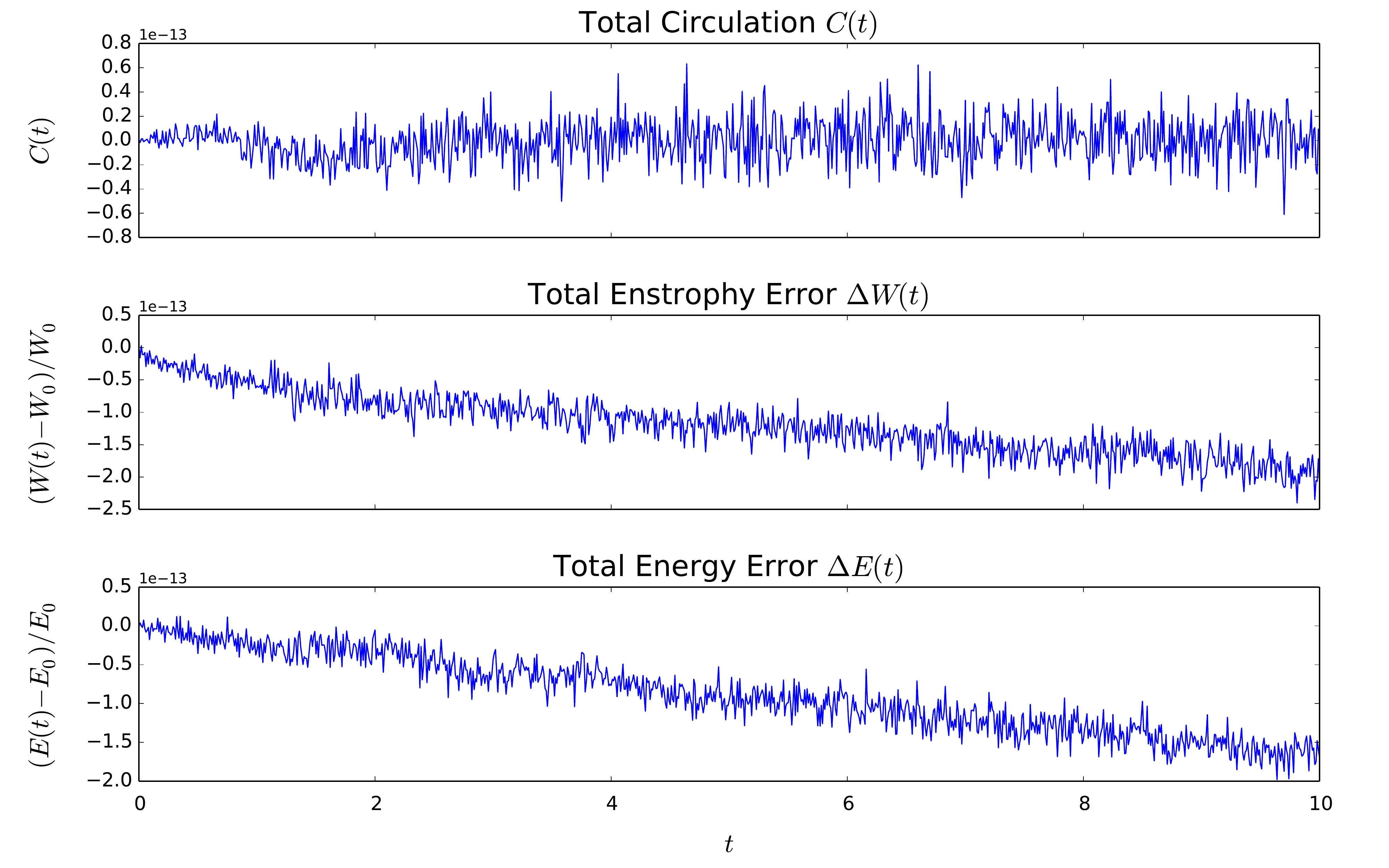}
	}

	\caption{Lamb Dipole. Top: Vorticity $\omega$ at $t=0.0$ and $10.0$. Bottom: Conservation Laws.}
	\label{fig:lamb_dipole}
\end{figure}

After $10,000$ timesteps (at $t=10$) no distortion of the lamb dipole is visible (compare figures~\subref*{fig:vorticity_lamb_dipole_vorticity_initial} and \subref*{fig:vorticity_lamb_dipole_vorticity_final}). The dipole keeps its shape just as it is supposed to. It does not end up exactly at the initial position due to the finite size of the box which reduces the velocity of the reference frame \cite[section IV.B]{NielsenRasmussen:1997}. 
 Figure~\subref*{fig:vorticity_lamb_dipole_timetraces} shows that circulation, enstrophy and energy are preserved to about machine accuracy. There is a small drift in the errors of enstrophy and energy, but this is of $\mathcal{O} (10^{-13})$ and therefore comparable to the tolerance of the nonlinear solver.

\subsubsection*{Gaussian Vortex}

\begin{figure}[p]
	\centering
	\subfloat[][$t=0$]{
	\includegraphics[width=.42\textwidth]{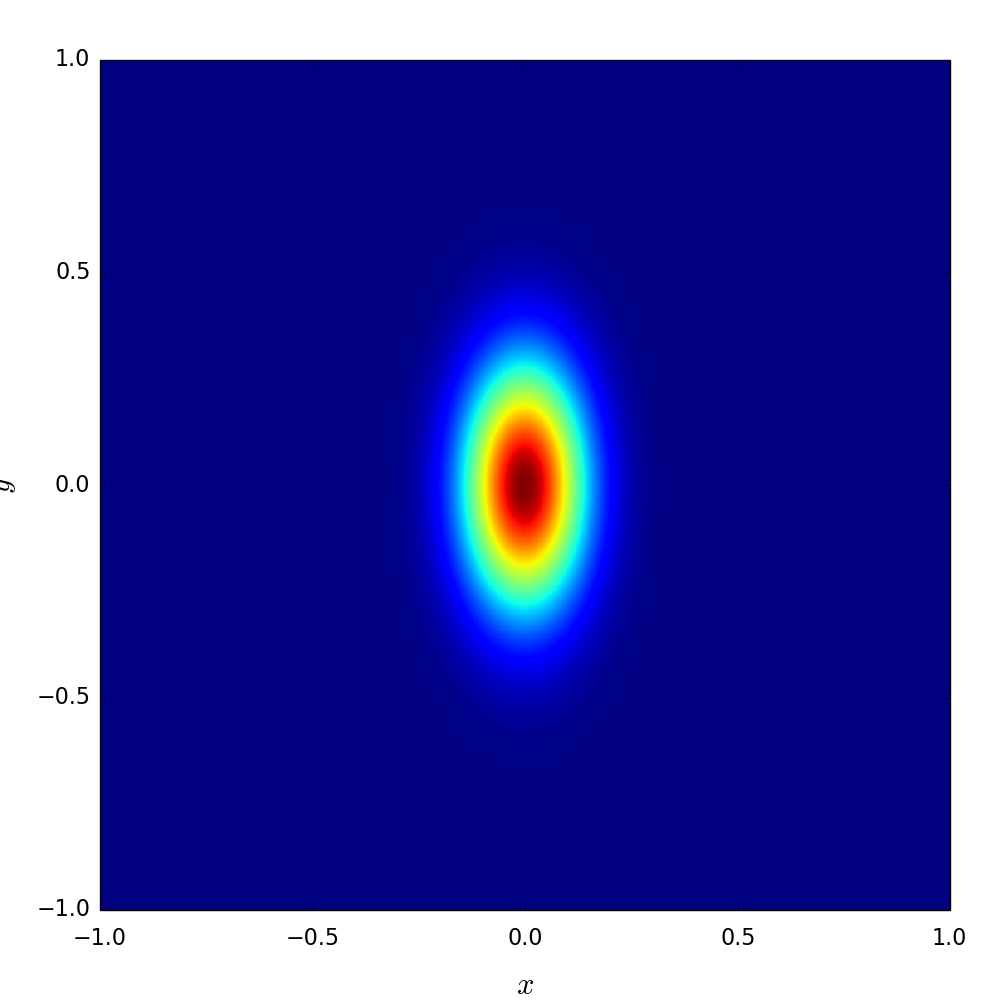}
	}
	\subfloat[][$t=2$]{
	\includegraphics[width=.42\textwidth]{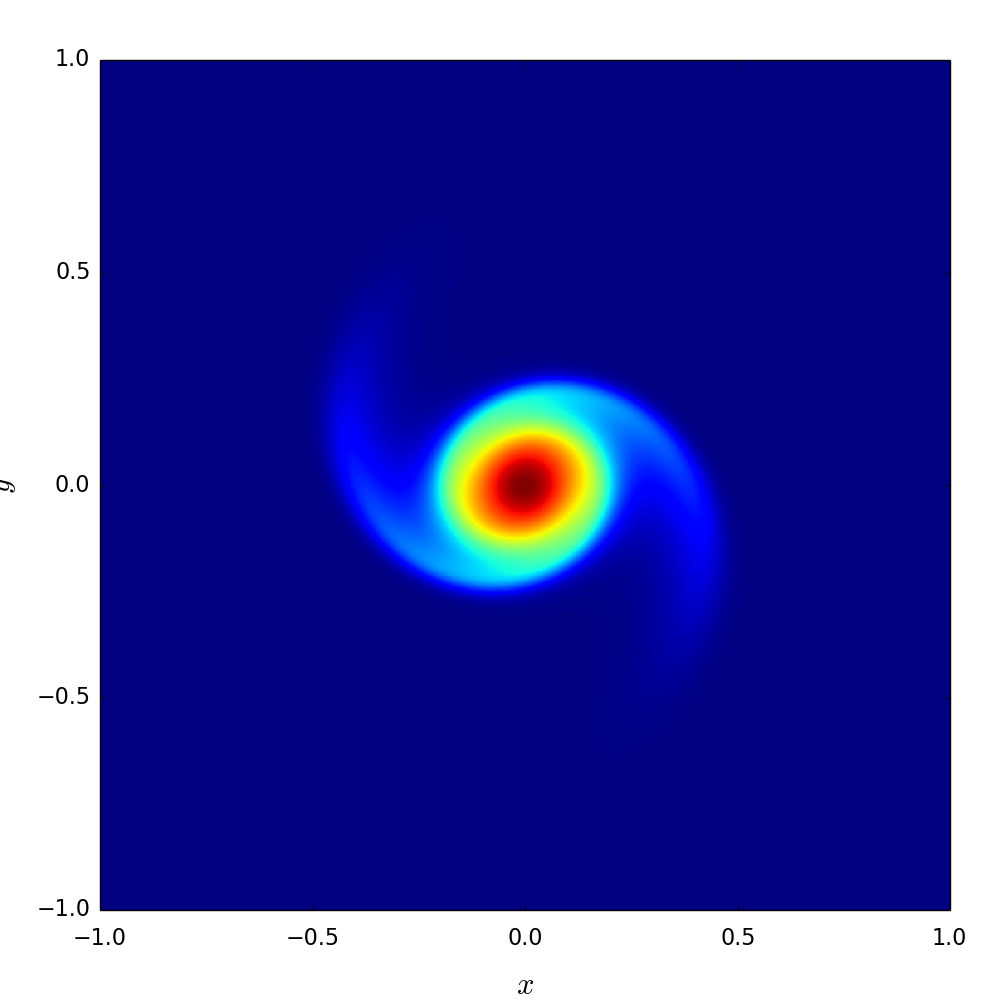}
	}
	
	\subfloat[][$t=4$]{
	\includegraphics[width=.42\textwidth]{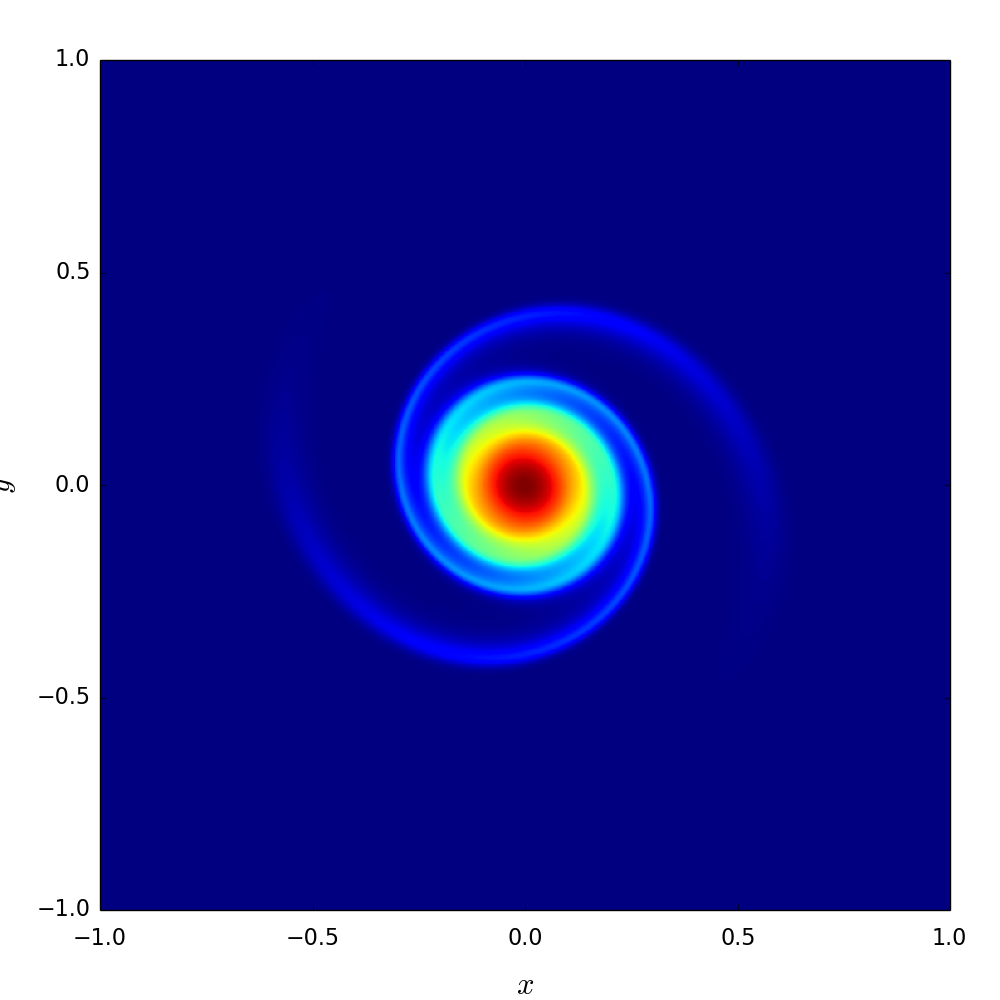}
	}
	\subfloat[][$t=6$]{
	\includegraphics[width=.42\textwidth]{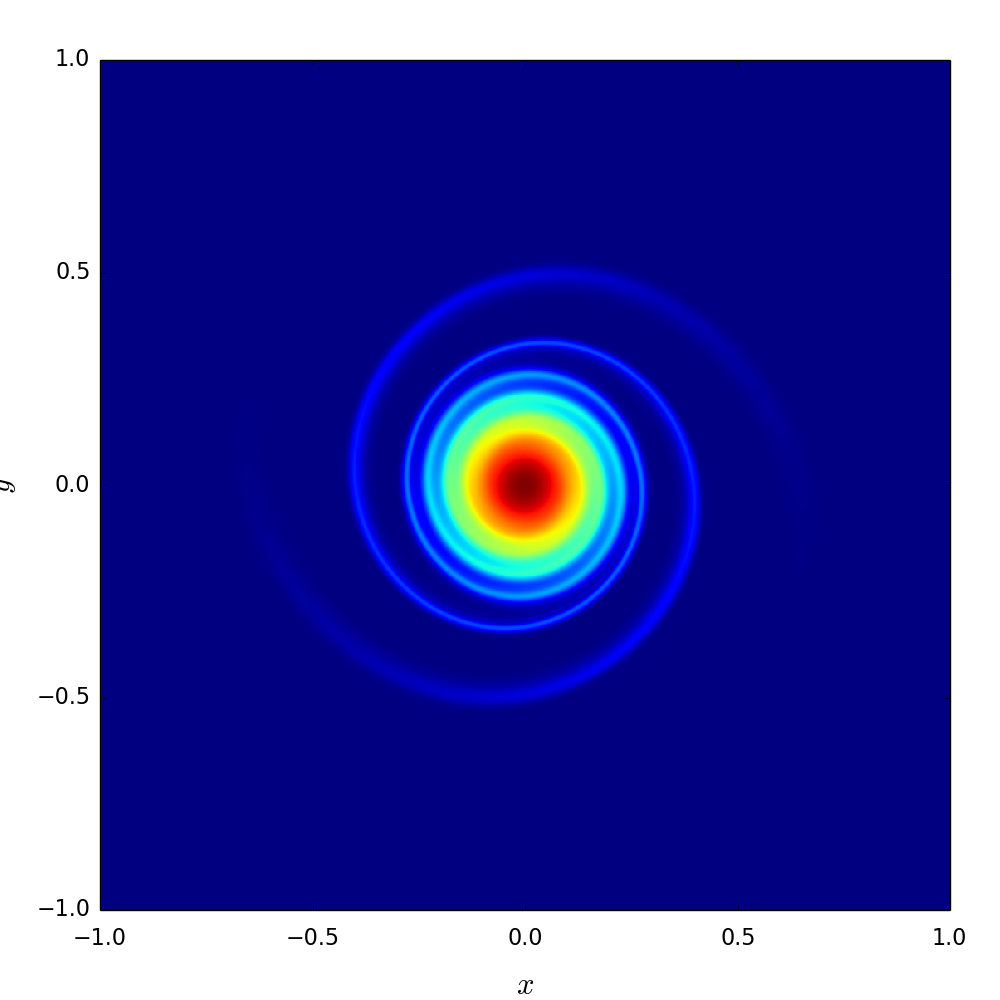}
	}
	
	\subfloat[][$t=8$]{
	\includegraphics[width=.42\textwidth]{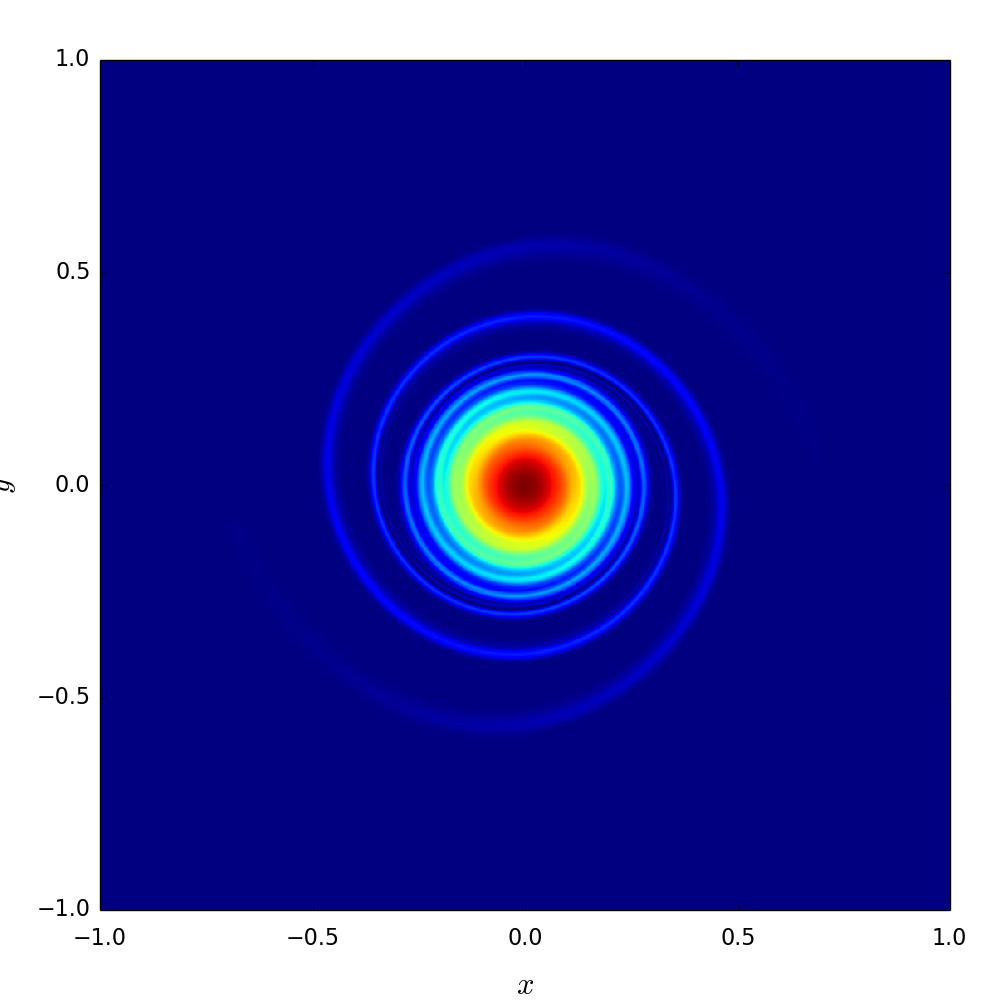}
	}
	\subfloat[][$t=10$]{
	\includegraphics[width=.42\textwidth]{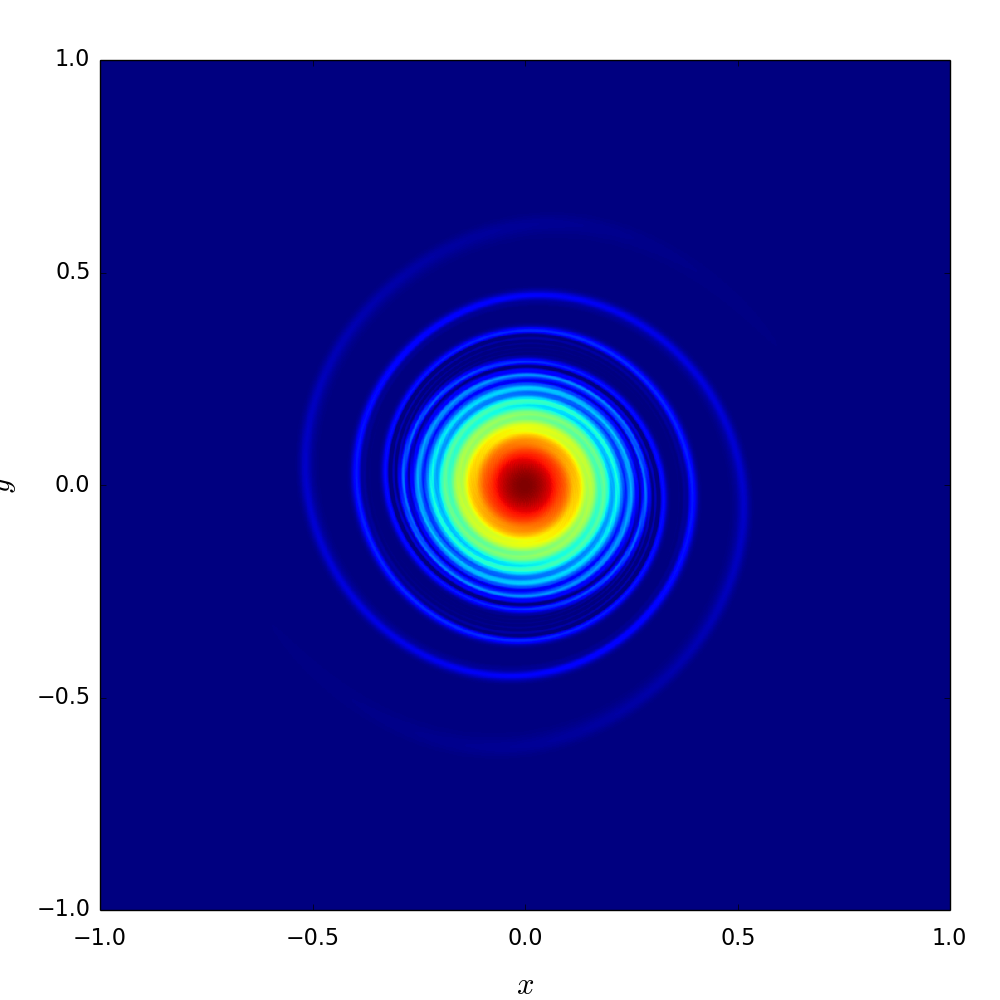}
	}
	
	\caption{Gaussian Vortex. Vorticity $\omega$ at $t=0$, $2$, $4$, $6$, $8$, $10$.}
	\label{fig:vorticity_gaussian_vortex_vorticity}
\end{figure}

Vorticity is initialised by a nonsymmetric Gaussian,
\begin{align}
\omega_{0} (x,y) = f(x, x_{0}, \sigma_{x}) \, f(y, y_{0}, \sigma_{y}) ,
\end{align}
with
\begin{align}
f(x, x_{0}, \sigma) = \dfrac{1}{\sigma \, \sqrt{2 \pi}} \, \exp \bigg( - \dfrac{1}{2} \bigg[ \dfrac{x-x_{0}}{\sigma} \bigg]^2 \bigg) ,
\end{align}
which develops into a spiral structure.
The domain is $\Omega = [-1, +1) \times [-1, +1)$ with periodic boundaries, the spatial resolution is $1024 \times 1024$ and the timestep is $10^{-2}$. The tolerance of the nonlinear solver is $10^{-11}$. We set $\sigma_{x} = 0.1$, $\sigma_{y} = 0.2$ and $x_{0} = y_{0} = 0$.
The solution follows the expected dynamics (c.f., Figure~\ref{fig:vorticity_gaussian_vortex_vorticity}) and conservation of circulation, enstrophy and energy is respected to machine accuracy (c.f., Figure~\subref*{fig:vorticity_gaussian_vortex_timetraces}).

\subsubsection*{Vortex Sheet Rollup}

\begin{figure}[p]
	\centering
	\subfloat[][$t=0.00$]{
	\includegraphics[width=.42\textwidth]{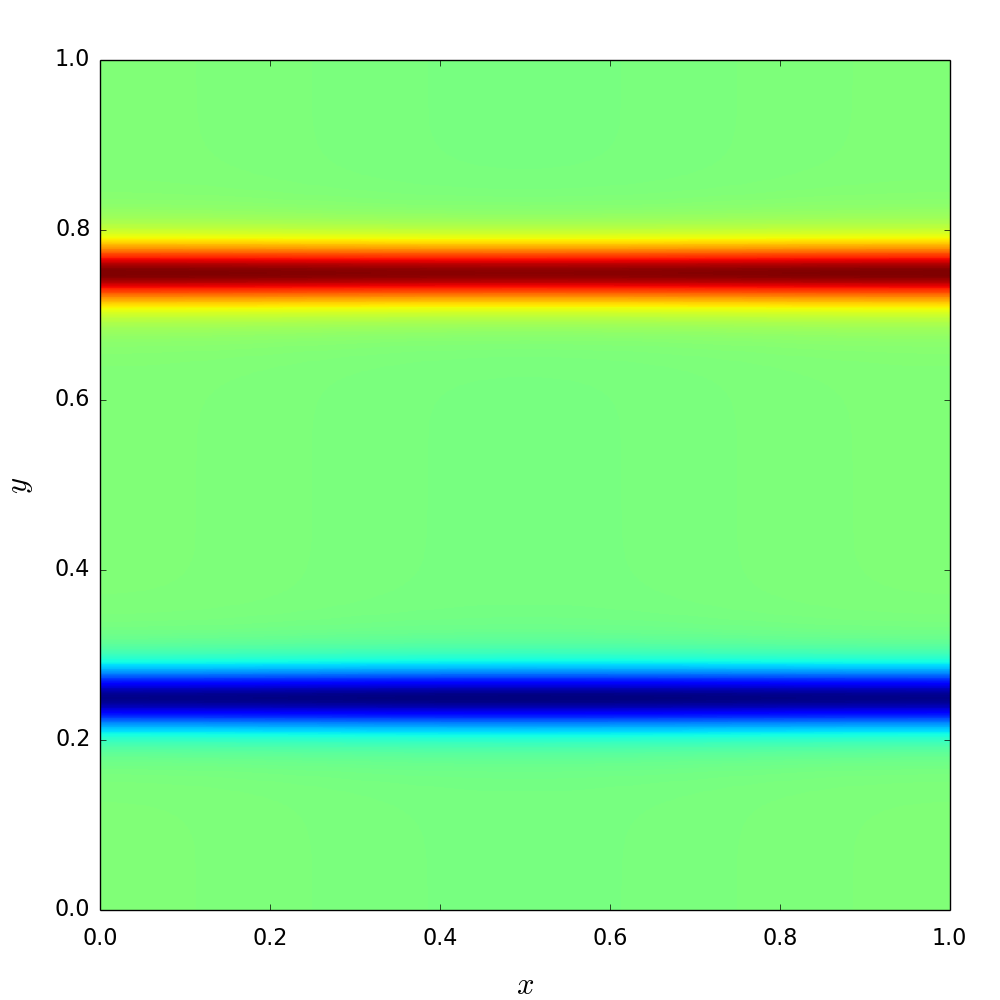}
	}
	\subfloat[][$t=1.00$]{
	\includegraphics[width=.42\textwidth]{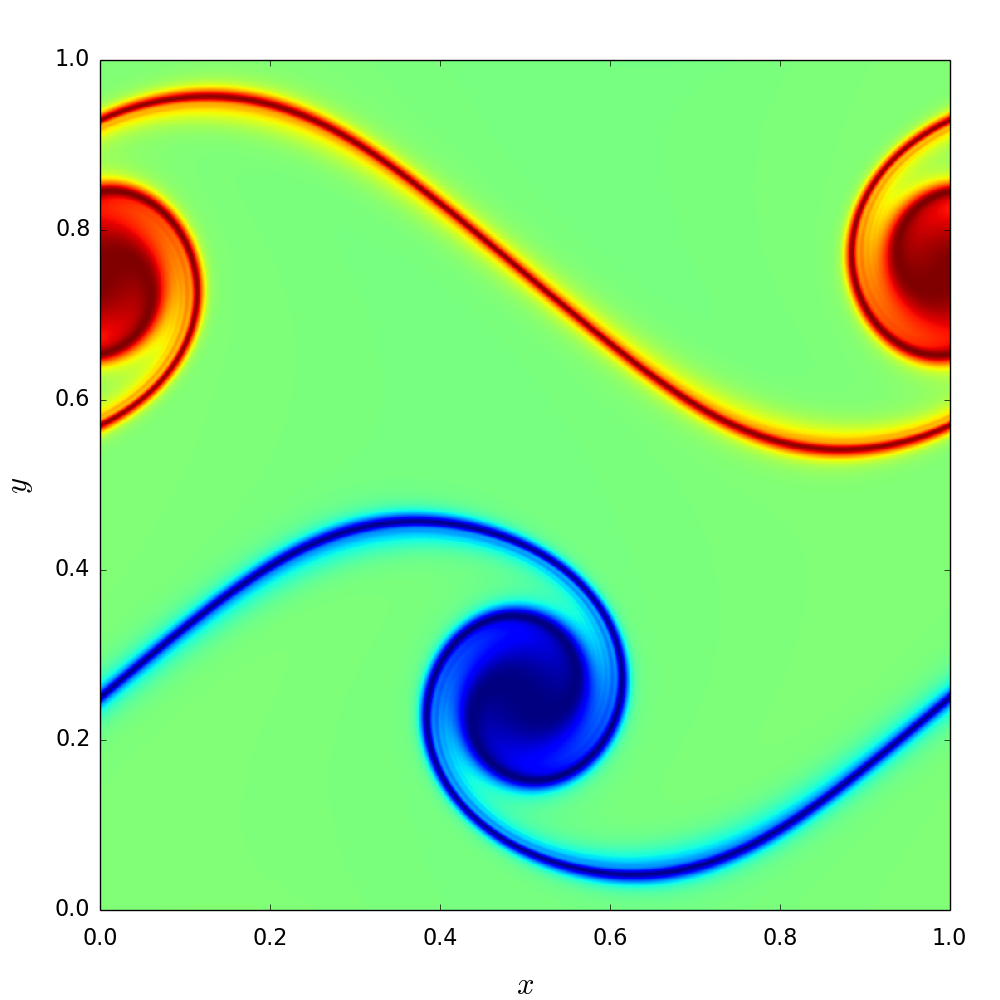}
	}
	
	\subfloat[][$t=1.25$]{
	\includegraphics[width=.42\textwidth]{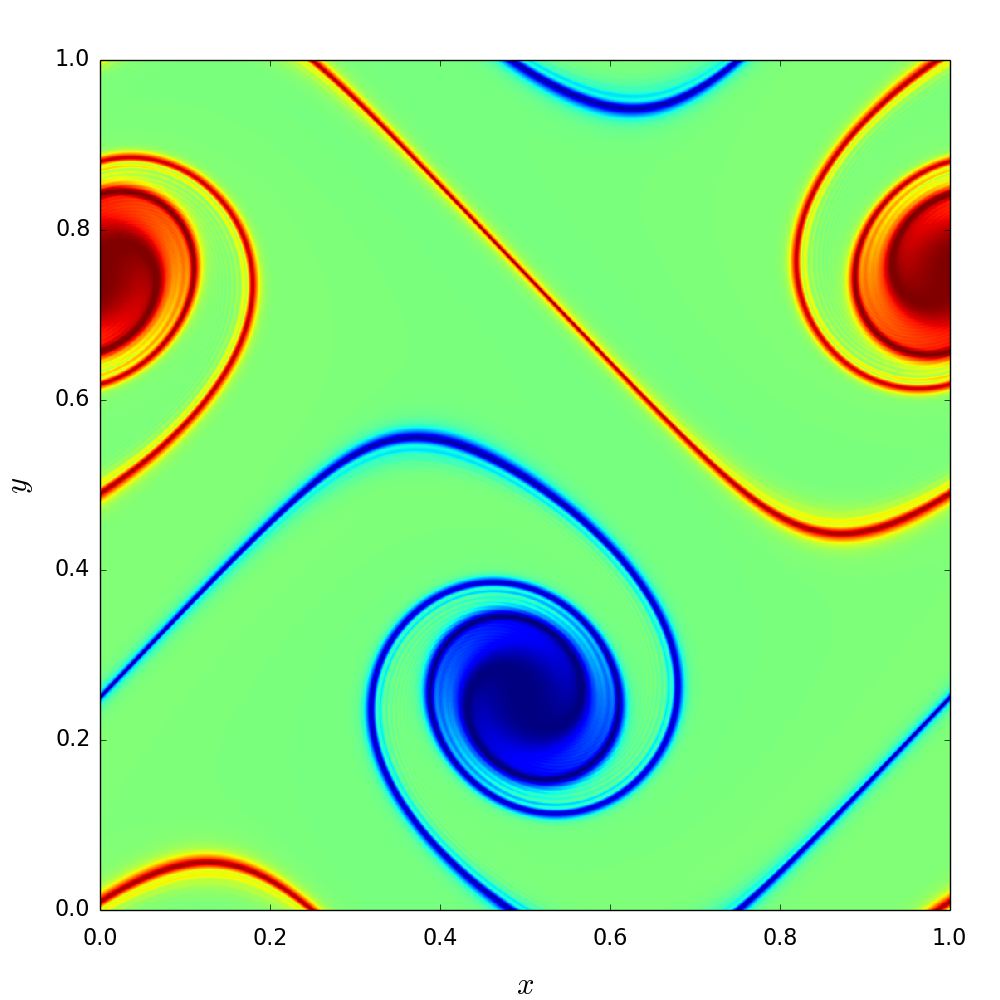}
	}
	\subfloat[][$t=1.50$]{
	\includegraphics[width=.42\textwidth]{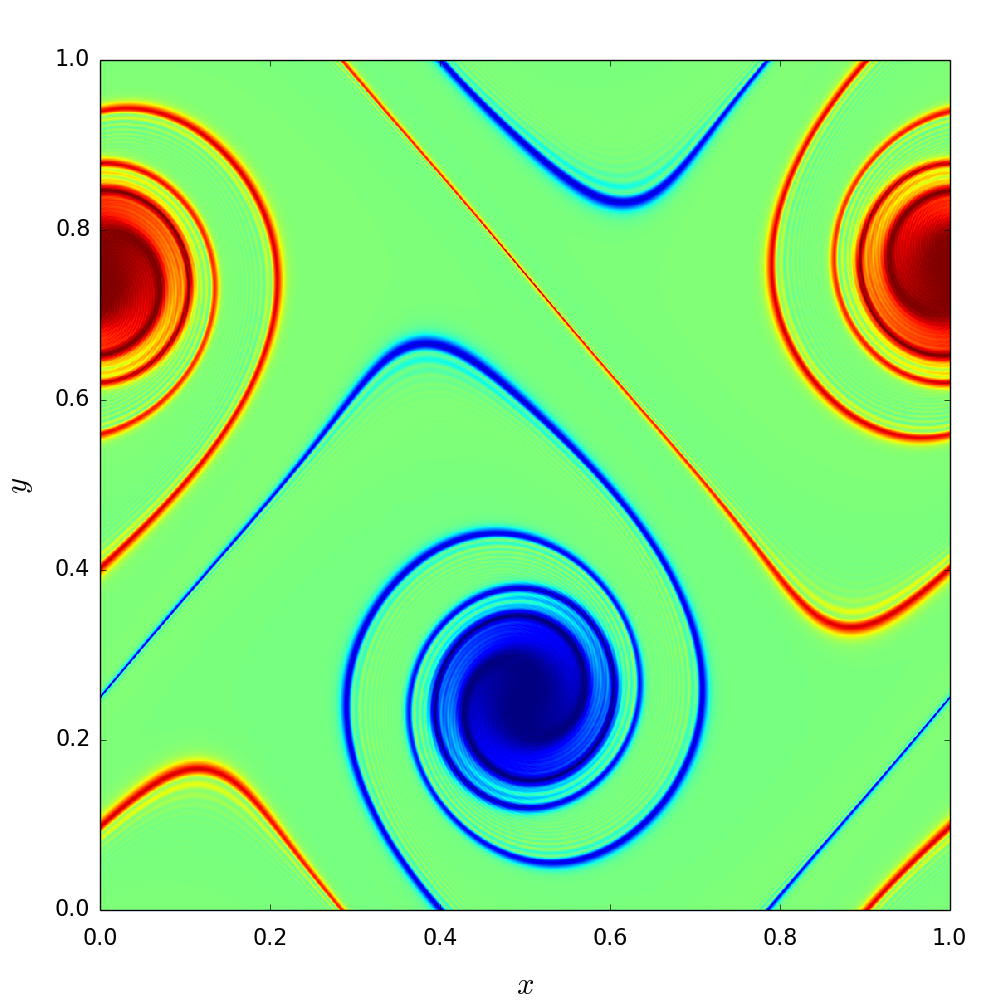}
	}
	
	\subfloat[][$t=1.75$]{
	\includegraphics[width=.42\textwidth]{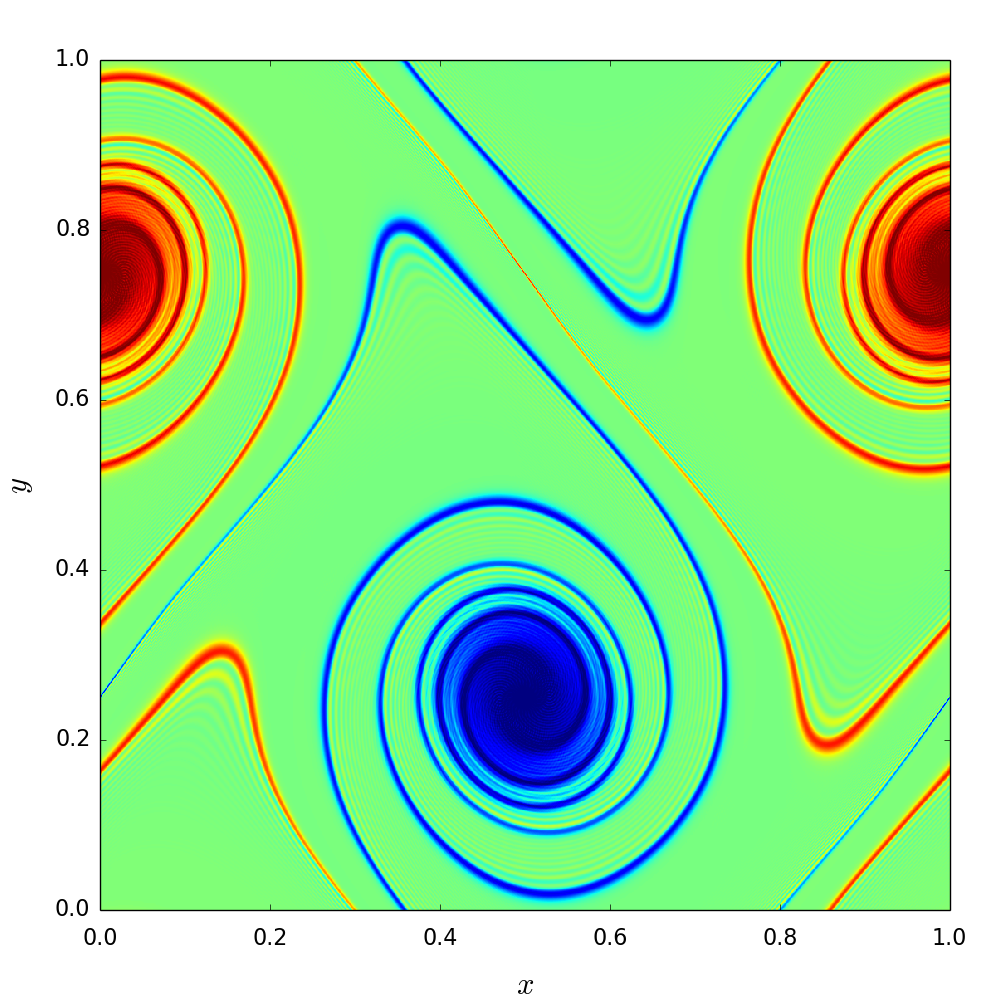}
	}
	\subfloat[][$t=2.00$]{
	\includegraphics[width=.42\textwidth]{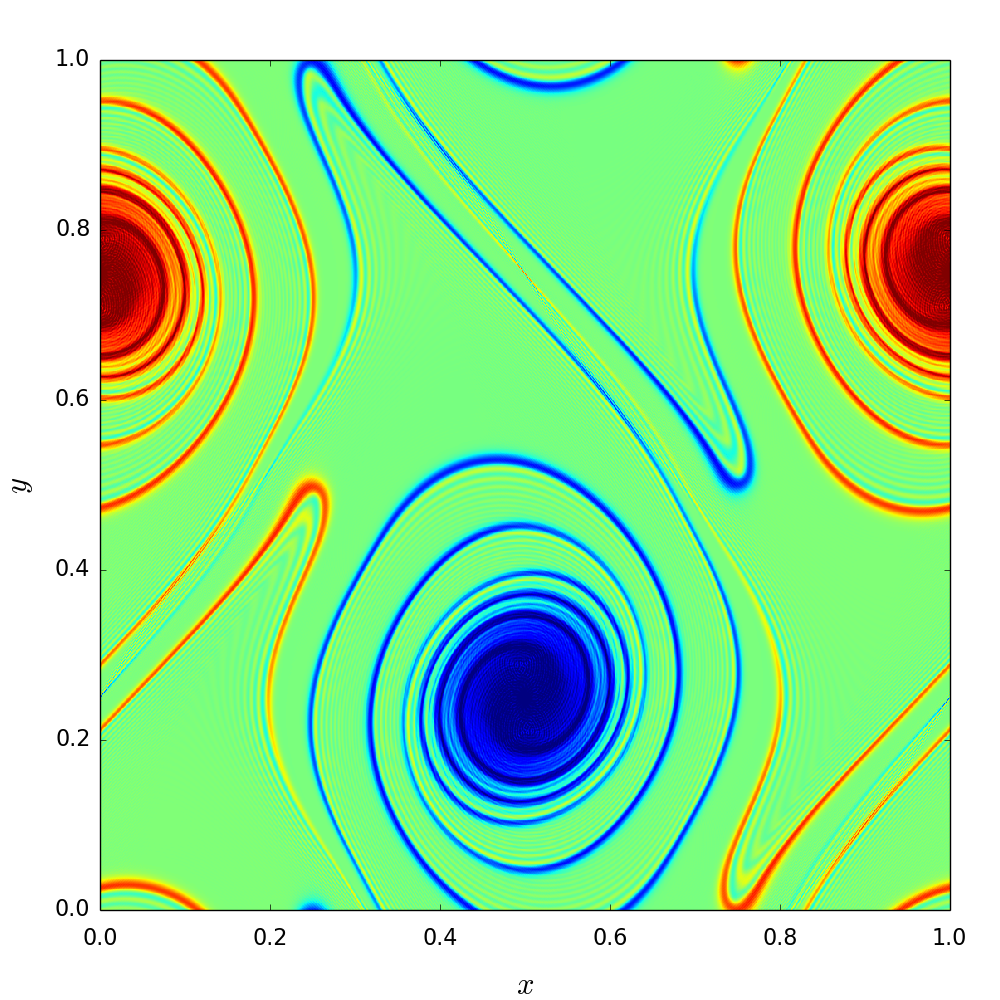}
	}
	
	\caption{Vortex sheet rollup. Vorticity $\omega$ at $t=0.00$, $1.00$, $1.25$, $1.50$, $1.75$, $2.00$.}
	\label{fig:vorticity_vortex_sheet_rollup_vorticity}
\end{figure}

\begin{figure}[p]
	\centering
	\subfloat[][Gaussian Vortex]{
	\label{fig:vorticity_gaussian_vortex_timetraces}
	\includegraphics[width=\textwidth]{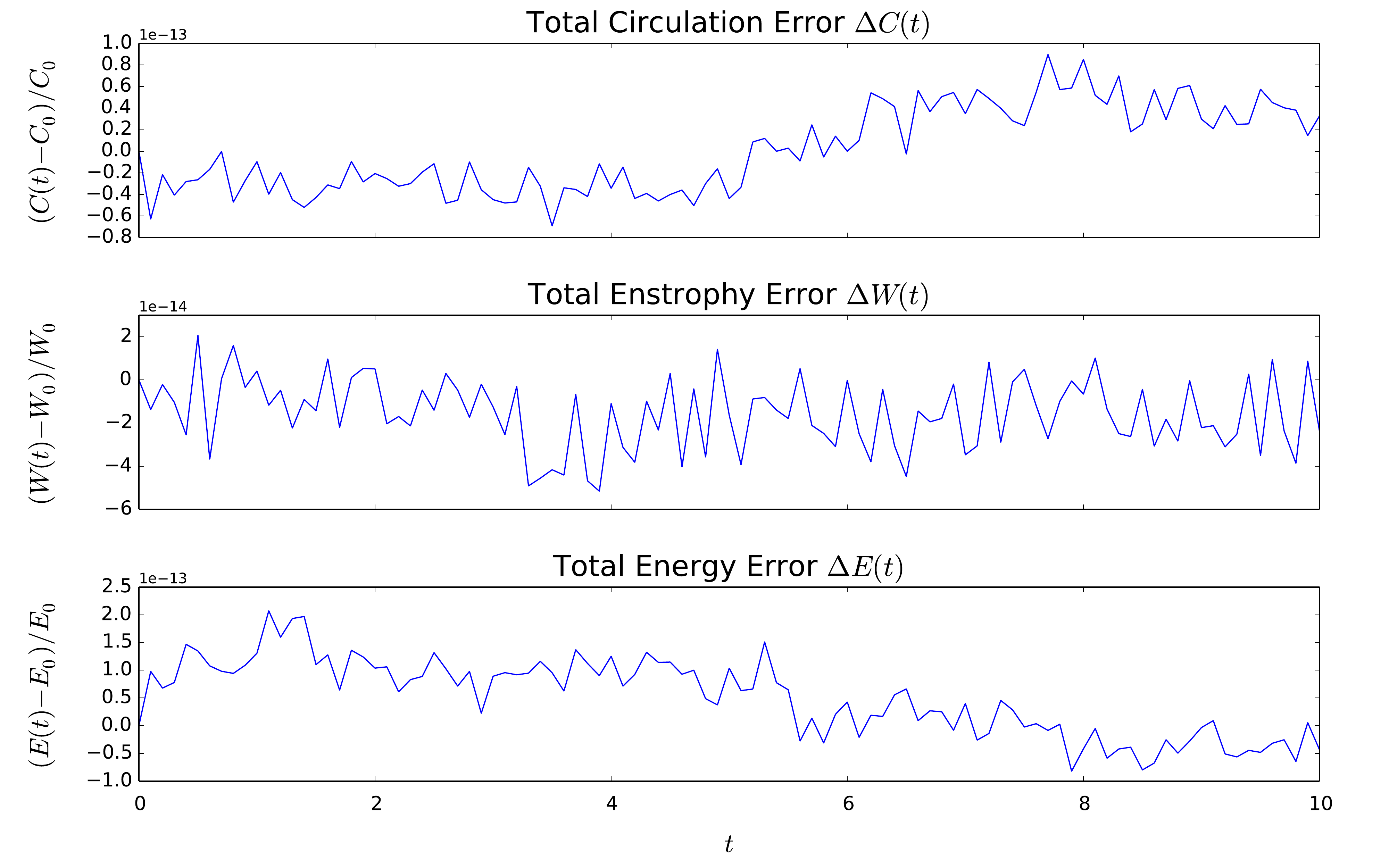}
	}

	\subfloat[][Vortex Sheet Rollup]{
	\label{fig:vorticity_vortex_sheet_rollup_timetraces}
	\includegraphics[width=\textwidth]{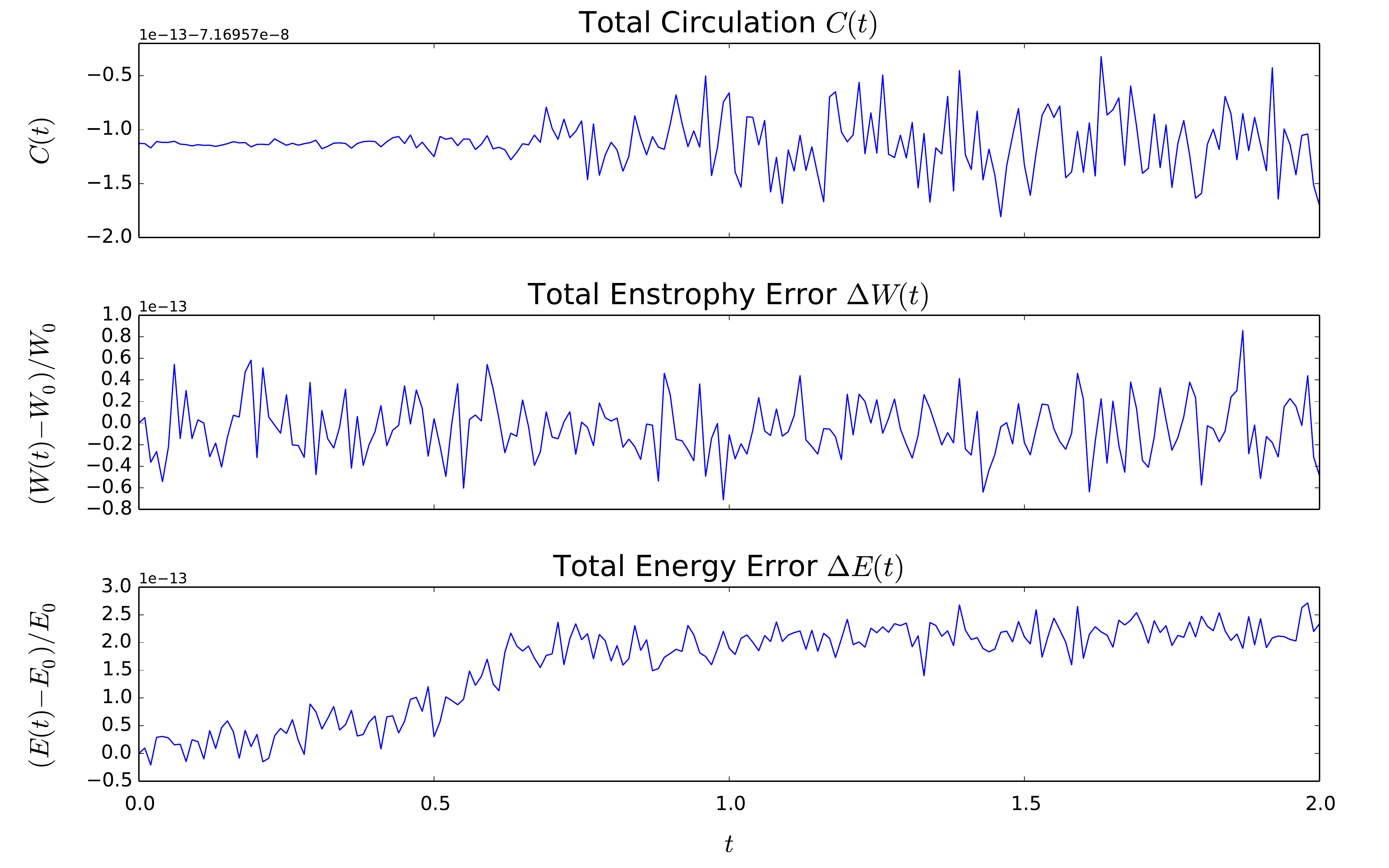}
	}

	\caption{Conservation Laws. The errors in the circulation, enstrophy and energy are all of the order of the machine accuracy.}
\end{figure}

The vortex sheet is a popular model in fluid dynamics used to approximate shear layers.
Following \cite{Deville:2002}, we use initial conditions
\begin{align}
u_{0} = \bigg( 
\begin{array}{ll}
\tanh (\rho ( y - 0.25 )) & \text{for} \; y \leq 0.5 , \\
\tanh (\rho ( 0.75 - y )) & \text{for} \; y > 0.5 ,
\end{array} , \quad
0.05 \, \sin (2 \pi x)
 \bigg)^{T} ,
\end{align}
which corresponds to the initial vorticity
\begin{align}
\omega_{0} (x,y) = \nabla_{\perp} \cdot u_{0} =
0.1 \pi \cos(2 \pi x) + 
\begin{cases}
- \dfrac{\rho}{\cosh^2 (\rho (y - 0.25))} & y \leq r , \\
+ \dfrac{\rho}{\cosh^2 (\rho (0.75 - y))} & y >    r .
\end{cases}
\end{align}
The domain is $\Omega = [0, 1) \times [0,1)$ with periodic boundaries, discretised by a grid of $1024 \times 1024$ points, and the timestep is $h_t = 10^{-3}$. The tolerance of the nonlinear solver is $10^{-10}$ and the parameter $\rho = 30$.
We observe good qualitative agreement of the solution in Figure~\ref{fig:vorticity_vortex_sheet_rollup_vorticity} with \cite[p. 328]{Deville:2002} together with excellent conservation properties (c.f., Figure~\subref*{fig:vorticity_vortex_sheet_rollup_timetraces}).

\clearpage

\section{Summary and Outlook}
\label{sec:summary}

We created a link between variational integrators, formal Lagrangians and Ibragimov's theory of conservation laws for arbitrary differential equations.
The proposed method allows us to derive variational discretisation schemes for potentially arbitrary systems of differential equations, even the ones that do not possess a classical variational formulation.
Thereby we were able to extend the applicability of the variational integrator method to a class of systems much larger than originally envisaged.
The main strength of our method is that it allows for the straight forward design of numerical schemes that respect certain conservation laws of the system at hand.

We presented an introduction to the theory of variational integrators that tries to make the geometric framework of the variational integrators available also to non specialists. Thereby we hope to make this class of methods accessible to a wider audience.
We extended the discrete Noether theorem to include discrete divergence symmetries.

The power of the method was demonstrated by several examples that are prototypical for problems arising in fluid dynamics and plasma physics. We emphasised the analysis of discrete conservation laws, which is seldom found in the variational integrator literature, and verified these theoretical properties in numerical experiments.
In follow-up papers we will present numerical results for the Vlasov-Poisson system  \cite{KrausMajSonnendruecker:2015} as well as ideal and reduced magnetohydrodynamics \cite{KrausMaj:2015, KrausTassi:2015}. There, explicit computations will demonstrate the favourable properties of the variational integrators for more elaborate and challenging applications.

Remarkably, we recovered Arakawa's discretisation of the Poisson bracket combined with a symplectic integrator in time. That is, we constructed a spacetime generalisation of Arakawa's discretisation. With our method, it is also straight-forward to generalise Arakawa's method to higher spatial dimensions and to higher order schemes.

An open question within our framework is the meaning of the discrete multisymplectic form arising from the boundary terms of the action sum (see \cite{MarsdenPatrick:1998} for further details) and its restriction to the original system. We observed that the multisymplectic form vanishes identically if the extended system is self-adjoint in the sense of Ibragimov. This is a topic we have not included in our discussion but one that certainly deserves attention.
 
The limits of the finite-difference approach to the discretisation of the Lagrangian became obvious in several places in this work. The treatment of discrete divergence symmetries is not straight-forward so that we are limited to global conservation laws. We cannot treat horizontal transformations in the symmetry generator but only vertical transformations. This is a problem for more complicated systems like the vorticity equations or the nonlinear advection equation (inviscid Burgers equation).
For nonlinear systems there is the additional complication that even if the system is self-adjoint on the continuous level, it will in general not be self-adjoint on the discrete level due to the absence of a discrete Leibniz rule.
Moreover, we are limited to low-order schemes. While it is in principle possible to design higher-order methods, this quickly becomes very cumbersome and confusing, especially in the analysis of the discrete conservation laws.
All those issues appear easier to deal with in a Galerkin framework with finite elements or splines as basis functions. The former has already been considered to a certain extend \cite{Chen:2008}, that is for the discretisation of the Lagrangian and the approximation of the action integral, but not for the analysis of conservation laws. The latter is a topic under active development \cite{Kraus:2015:Splines}.

So far, we only considered examples of advection-diffusion type. A detailed analysis of the applicability of the proposed method to other types of equations would be most interesting.

\section*{Acknowledgements}

The first author would like to thank Bruce D. Scott for his guidance and the freedom to pursue the path that seemed most interesting, as well as Eric Sonnendr\"ucker for supporting this project.
We are thankful to Emanuele Tassi and Jonathan Squire for raising and discussing important questions that added a lot to the understanding of the matters presented as well as reading various versions of the manuscript.
Furthermore, we wish to thank Claudio Dappiaggi for the discussions on the inverse problem of calculus of variations and for suggesting the paper by \citet{BampiMorro:1982}.

\bibliographystyle{plainnat}
\bibliography{vi_nonvar}

\end{document}